\documentclass[11pt,a4paper]{article}

\usepackage[english,francais]{babel}
\usepackage[utf8]{inputenc}
\usepackage[T1]{fontenc}

\usepackage{hyperref} 

\usepackage{amsmath,amssymb,amsfonts,latexsym,xspace,amsthm,array}
\usepackage{cases} 
\usepackage{stmaryrd} 
\usepackage{adjustbox} 
\usepackage{blkarray} 

\usepackage{tikz}
\usepackage{tikz-cd}
\usetikzlibrary{patterns}

\usepackage{geometry}
\newenvironment{pf}{{\em \noindent Preuve:}}{ \hfill \qed\smallskip}

\newtheorem{theo}{Théorème}[section]
\newtheorem{lem}{Lemme}[section]
\newtheorem{prop}{Proposition}[section]
\newtheorem*{prop*}{Proposition} 

\newtheorem{defi}{Définition}[section]
\newtheorem*{defi*}{Définition} 

\theoremstyle{definition} 
\newtheorem{rem}{Remarque}[section]

\newtheorem*{notat*}{Notation}
\newtheorem{notat}{Notation}[section]
\newtheorem{conv}{Conventions}[section]

\def\Spec{\mathrm{Spec}}
\def\Frac{\mathrm{Frac}}

\def\deg{\mathrm{deg}}

\renewcommand{\textbf}[1]{\begingroup\bfseries\mathversion{bold}#1\endgroup} 

\begin{document}

\title{Calcul du Frobenius divisé modulo $p$ sur la cohomologie cristalline de certains revêtements de la droite projective}
\author{Amandine Pierrot\\
   Université de Strasbourg,\\
   France,\\
   \texttt{pierrot@math.unistra.fr}}
\date{\today}

\maketitle

\begin{abstract} Dans cet article nous décrivons une famille de revêtements modérément ramifiés de la droite projective sur un corps fini, pour lesquels nous effectuons le calcul de la matrice du Frobenius divisé cristallin. Les formules que nous obtenons généralisent les formules de Hasse-Witt classiques dans le cas des courbes hyperelliptiques. Un des outils est un résultat récent de Huyghe-Wach qui démontre que le Frobenius divisé cristallin coïncide avec le morphisme explicite, construit par Deligne-Illusie en 1987, pour établir la dégenerescence de la suite spectrale de Hodge vers de Rham dans le cas algébrique. \end{abstract}
\selectlanguage{english}
\begin{abstract} In this paper we describe a family of tamely ramified coverings of the projective line over a finite field, for which we compute the matrix of the divided crystalline Frobenius. The formulas we obtain generalize the classical Hasse-Witt formulas in the case of hyperelliptic curves. Our result relies on a result of Huyghe-Wach which shows that the divided crystalline Frobenius coincides with the explicit morphism, constructed by Deligne-Illusie in 1987, for their proof of the degeneration of the Hodge-de Rham spectral sequence in the algebraic case. \end{abstract}

\section{Introduction}

Soit $k$ un corps fini et $X_0$ une courbe projective lisse sur $\Spec(k)$. Hasse a démontré \cite{Hasse1935} que le rang de la matrice de Hasse-Witt est un invariant important de $X_0$ et donc de sa jacobienne. Lorsque le rang est minimal, c'est-à-dire nul, on dit que la courbe est supersingulière; lorsqu'il est maximal et donc égal à $g$, genre de la courbe, elle est dite ordinaire. Dans le cas d'une courbe elliptique, définie par l'équation $y^2=t(t-1)(t-\lambda)$, le calcul de cet invariant est déjà fait par Hasse en loc. cit. et il est donné, à un inversible modulo $p$ près.  Dans le cas des courbes hyperelliptiques et de certaines courbes de Fermat, le calcul a été fait par Gonzales \cite{MR1447179}. Par ailleurs pour les courbes hyperelliptiques on dispose depuis peu pour le calculer d'un algorithme en temps polynomial \cite{MR3240808}.

\vspace{+5mm}

Parallèlement, Kedlaya \cite{Memoire} s'est rendu compte que, lorsque $X_0$ est une courbe hyperelliptique, il est possible de calculer la cohomologie rigide de l'ouvert complémentaire des points de ramification de $X_0$. On peut en déduire la fonction Zêta de cette courbe et donc le nombre de points de la courbe. Ce calcul a été généralisé par Tuitman dans le cas de revêtements séparables de degré $n$ de la droite projective sur un corps fini \cite{MR3434890}\cite{MR3631366}. Même si le contexte de notre travail est différent de celui de ces algorithmes, nous avons réutilisé certaines idées de Kedlaya, notamment celle qui consiste à se placer sur l'ouvert complémentaire du lieu de ramification pour trouver un relèvement explicite du Frobenius (cela sera explicité dans la section \ref{paszero}).

\vspace{+5mm}

Considérons $k$ un corps fini, $W=W(k)$ l'anneau des vecteurs de Witt de $k$ et $K=\Frac(W)$. On se propose dans cet article de calculer un invariant plus complet que l'invariant de Hasse-Witt pour une certaine classe de courbes projectives lisses qui sont des revêtements séparables de degré $n$ de la droite projective sur $k$. Décrivons plus précisément ce dont il s'agit. Soit $X$ une courbe projective lisse sur $\Spec(W)$ qui est un revêtement séparable de degré $n$ de la droite projective relative sur $\Spec(W)$, notons $X_0$ sa fibre spéciale et $X'_0$ la courbe obtenue après changement de base par le Frobenius sur $k$. Commençons par quelques rappels sur l'invariant de $X_0$ que nous souhaitons calculer. En 1987 Fontaine et Messing ont donné une construction du Frobenius divisé cristallin $\Psi$, qui est à valeurs dans le groupe de cohomologie cristalline, $H^1_{cris}(X_0)$ \cite{MR902593}. Or il se trouve que pour les courbes que nous considérons ici, $H^1_{cris}(X_0)$ est un $W$-module libre de rang fini égal à $2g$. Sous nos hypothèses on sait, grâce au théorème de comparaison de Berthelot \cite{Berth}, que le groupe  $H^1_{cris}(X_0)$ s'identifie à $H^1_{DR}(X)$, et cet isomorphisme donne un isomorphisme modulo $p$
$$H^1_{cris}(X_0)/pH^1_{cris}(X_0) \simeq H^1_{DR}(X_0).$$
Finalement, modulo $p$, le Frobenius divisé $\Psi$ donne un morphisme de $k$-espaces vectoriels
\begin{center} \label{isoDI}
$\Phi' : H^0(X'_0,\Omega^1_{X'_0})\bigoplus H^1(X'_0,\mathcal{O}_{X'_0})  \rightarrow H^1_{DR}(X_0)$
\end{center}
qui, sous nos hypothèses, est un isomorphisme. Par ailleurs, puisque le Frobenius sur $k$ est un isomorphisme, la projection 
$X'_0 \rightarrow X_0$ est un isomorphisme de schémas, et on en déduit des isomorphismes 
semi-linéaires par rapport au Frobenius sur $k$
\begin{center}
$H^0(X_0,\Omega^1_{X_0})\simeq H^0(X_0',\Omega^1_{X_0'})$ et $ H^1(X_0,\mathcal{O}_{X_0})\simeq H^1(X'_0,\mathcal{O}_{X'_0})$.
\end{center}
En composant l'isomorphisme précédent avec ces isomorphismes semi-linéaires, on en déduit un isomorphisme semi-linéaire de $k$- espaces vectoriels
\begin{center}
$\Phi : H^0(X_0,\Omega^1_{X_0})\bigoplus H^1(X_0,\mathcal{O}_{X_0}) \rightarrow H^1_{DR}(X_0)$.
\end{center}
C'est cet isomorphisme dont on cherche à calculer la matrice dans une base adaptée au scindage. Par définition, dans la base que nous avons choisie pour la représenter, le quart inférieur droit de la matrice redonne la matrice de Hasse-Witt, tandis que le quart supérieur gauche redonne la matrice de l'opérateur de Cartier \cite{DT}.  Il était déjà connu de Cartier que ces deux opérateurs se correspondent par la dualité de Serre sur $X_0$. L'intérêt de calculer ce nouvel invariant (la matrice du Frobenius divisé) est que cela permet, après application de l'algorithme de Wach \cite{QT}\cite{MR1468834}, de calculer le $(\varphi,\Gamma)$-module modulo $p$ associé à la représentation galoisienne $H^1_{et}(X_{\overline{K}},\mathbf{Q}_p)$,  où $\overline{K}$ est une clôture algébrique de $K$. Ce dernier point ne sera pas abordé ici.

\vspace{+5mm}

Pour mener à bien le calcul de la matrice du Frobenius divisé, nous utilisons un résultat récent de Huyghe-Wach \cite{Encours} qui démontre que $\Phi'$ coïncide avec le morphisme construit par Deligne-Illusie  en 1987 \cite{MR894379} pour démontrer la dégénérescence de la suite spectrale de Hodge vers de Rham. Le point remarquable de la construction de Deligne-Illusie est qu'elle est complètement algorithmique, au moins dans certains cas, comme l'illustre notre travail.
Signalons enfin qu'il est standard que $\Phi$ coïncide avec le Frobenius divisé de Mazur modulo $p$ \cite{MR0321932}, qui est construit à partir d'un scindage de la filtration de Hodge, donnant un isomorphisme
$$H^1_{DR}(X)\simeq H^0(X,\Omega^1_X)\bigoplus H^1(X,\mathcal{O}_X).$$ 
Enfin nous obtenons au passage une formule générale pour calculer la matrice de Hasse-Witt pour la famille de courbes considérée, formule analogue à celle obtenue par Elkin pour l'opérateur de Cartier \cite{TT}.

\vspace{+5mm}

Je remercie Christine Huyghe et Nathalie Wach qui ont encadré la thèse dont est issue cet article. Je remercie également Annette Huber-Klawitter pour m'avoir invitée à présenter mes premiers résultats en conférence. Enfin je remercie chaleureusement Xavier Caruso pour son astuce d'inversion qui a permis l'obtention de formules explicites et pour son aide pour les calculs de complexité.

\section{Définition de la courbe}\label{cas}

\subsection{Cadre de travail}

Dans ce texte on va considérer $p$ un nombre premier, $n$ un nombre entier premier à $p$ et $k$ un corps fini de caractéristique $p>0$. On notera $W=W(k)$ l'anneau des vecteurs de Witt de $k$, $K=\Frac(W)$ et $W_1=W/p^2W$. On écrit la droite projective $\mathbb{P}^1_W$ comme $\mathbb{P}^1_W=A\cup B$ avec $A=\Spec(W[t])$ et $B=\Spec(W[s])$ où $s=1/t$. \\

Considérons $l$ un nombre entier congru à $-1$ modulo $n$ tel que $l$ n'est pas multiple de $p$ et un polynôme \textbf{$f$} de $W[t]$ de degré $l$ tel que \textbf{$f$} et \textbf{$f'$} engendrent l'idéal unité de $W[t]$. On notera $r$ le nombre entier $(l+1)/n$ et \textbf{$f_2(s)$} la fraction rationnelle
\begin{center}
\textbf{$f_2$}$(s)=$\textbf{$f$}$(1/s)s^{rn}=$\textbf{$f$}$(1/s)s^{l+1}$
\end{center}
qui est un polynôme de $W[s]$ de degré $l+1$. Pour finir posons $z=s^ry$. Dans tout ce qui suit on notera $f$, $f'$, $f_2$ et $f'_2$ les classes de ces polynômes modulo $p$ qui sont donc des éléments de $k[t]$ et $k[s]$ respectivement. Nous faisons enfin l'hypothèse que $t\wedge f(t)=1$ dans $k[t]$. Nous verrons au lemme \ref{camarche} que cette hypothèse n'est en fait pas restrictive.

\begin{rem}
Remarquons qu'avec ces hypothèses les polynômes $f(t)$ et $f'(t)$ sont premiers entre eux dans $k[t]$.
\end{rem}

\begin{defi}
Avec les notations précédentes, on définit un schéma affine $U$ sur $A$ par 
\begin{center}
$U=\Spec(W[t,y]/(y^n-$\textbf{$f$}$(t)))$
\end{center}
et
\begin{center}
$W[t]\rightarrow \mathcal{O}_U(U), \ t\mapsto t$
\end{center}
ainsi qu'un ouvert affine $V$ sur $B$ par
\begin{center}
$V=\Spec(W[s,z]/(z^n-$\textbf{$f_2$}$(s)))$
\end{center}
et
\begin{center}
$W[s]\rightarrow \mathcal{O}_V(V), \ s\mapsto s$.
\end{center}
\end{defi}

\begin{prop}
Les courbes $U$ et $V$ ainsi définies se recollent au dessus de $A\cap B$. Les conditions de recollement sont données par $t=s^{-1}$ et $z=s^ry$. 
\end{prop}

\begin{pf}
Montrons que l'on peut définir des isomorphismes 
$$\rho : W[t,t^{-1}]\otimes_{W[t]}\mathcal{O}_U(U) \rightarrow W[s,s^{-1}]\otimes_{W[s]}\mathcal{O}_V(V)$$
$$\tilde{\rho} : W[s,s^{-1}]\otimes_{W[s]} \mathcal{O}_V(V)\rightarrow W[t,t^{-1}]\otimes_{W[t]} \mathcal{O}_U(U)$$
inverses l'un de l'autre. Commençons par rappeler que $W[t,t^{-1}]\otimes_{W[t]} \mathcal{O}_U(U)\simeq W[t,t^{-1},y]/(y^n-$\textbf{$f$}$(t))$ et $W[s,s^{-1}]\otimes_{W[s]} \mathcal{O}_V(V)\simeq W[s,s^{-1},z]/(z^n-$\textbf{$f_2$}$(s))$. Par ailleurs il existe un morphisme d'algèbres $W[t,t^{-1},y] \rightarrow W[s,s^{-1},z]$ qui envoie $t$ sur $s^{-1}$ et $y$ sur $s^{-r}z=s^{-(l+1)/n}z$. Ce morphisme passe au quotient en un morphisme d'algèbres $$\rho : W[t,t^{-1}]\otimes_{W[t]}\mathcal{O}_U(U) \rightarrow W[s,s^{-1}]\otimes_{W[s]}\mathcal{O}_V(V).$$ On définit de même le morphisme $\tilde{\rho} : W[s,s^{-1}]\otimes_{W[s]} \mathcal{O}_V(V)\rightarrow W[t,t^{-1}]\otimes_{W[t]} \mathcal{O}_U(U)$ avec $\tilde{\rho}(s)=t^{-1}$ et $\tilde{\rho}(z)=t^{-r}y$. Il est alors facile de vérifier que $\tilde{\rho}\circ \rho (t)=t$, $\tilde{\rho}\circ \rho (y)=y$ et $\rho \circ \tilde{\rho}(s)=s$, $\rho \circ \tilde{\rho}(z)=z$. Ainsi les ouverts $U$ et $V$ se recollent au dessus de l'ouvert $A\cap B$. \\
\end{pf}

\begin{defi}
En conservant les notations précédentes, on définit $X$ comme le recollement des deux schémas $U$ et $V$. C'est une courbe relative sur $\Spec(W)$. \\
\end{defi}

\begin{prop} \label{pik}
Par construction on a un morphisme de schémas $\Pi:X\rightarrow \mathbb{P}^1_W$. Le morphisme $\Pi$ est fini et plat. De plus les applications induites par $\Pi$ sur les fibres, $\Pi_0:X_0\rightarrow \mathbb{P}^1_k$ et $\Pi_K:X_K\rightarrow \mathbb{P}^1_K$, sont des morphismes finis séparables de degré $n$. En particulier $X$ est propre sur $\Spec(W)$.
\end{prop}

\begin{pf}
Le fait que le morphisme $\Pi:X\rightarrow \mathbb{P}^1_W$ soit fini et plat est une propriété locale sur $\mathbb{P}^1_W$, il suffit de montrer que $U\rightarrow \Spec(W[t])$ et $V\rightarrow \Spec(W[s])$ sont finis et plats, ce qui est équivalent à montrer que $W[t]\rightarrow \mathcal{O}_U(U)$ et $W[s]\rightarrow \mathcal{O}_V(V)$ sont finis et plats. Or $\mathcal{O}_U(U)$ est un $W[t]$-module libre de rang $n$ et de base $1,y,...,y^{n-1}$. De même $\mathcal{O}_V(V)$ est un $W[s]$-module libre de rang $n$ et de base $1,z,...,z^{n-1}$ ce qui permet de conclure. En ce qui concerne le morphisme sur la fibre générique c'est alors une conséquence de ce qui précède puisque platitude et finitude sont stables par changement de base. L'argument est le même pour le morphisme sur la fibre spéciale car $k$ est de caractéristique $p$ et qu'on a pris soin que $n$ soit premier à $p$. \\
\end{pf}

\begin{notat} \label{defisurfibrespe}
Notons $X_0$ la fibre spéciale de $X$, alors $X_0$ est la courbe sur $\Spec(k)$ définie comme étant le recollement des deux schémas $\mathcal{U}$ et $\mathcal{V}$ où
\[\mathcal{U}=\Spec(k[t,y]/(y^n-f(t))) \] 
et
\[\mathcal{V}=\Spec(k[s,z]/(z^n-f_2(s))) \]\\
\end{notat}

\begin{lem} \label{camarche}
Supposons que $X_0$ soit telle que $t\wedge f(t)=t$, alors quitte à faire une extension finie de $k$, il existe $u$ dans $W$ tel que \textbf{$g$}$(t)=$\textbf{$f$}$(t+u)$ vérifie $g(t)\wedge t=1$ $mod$ $p$.
\end{lem}

\begin{pf}
Considérons $k'$ une extension finie de $k$ de cardinal supérieur ou égal à $l+1$ telle que $f$ est scindé (à racines simples) dans $k'$. Notons $0,\alpha_2,...,\alpha_l$ les racines de $f$ et prenons $v$ dans $k'\backslash \{0,\alpha_2,...,\alpha_l\}$. Alors les racines de $\tau_v(f)$ sont $-v,\alpha_2-v,...,\alpha_l-v$ et donc $0$ n'est pas racine de $\tau_v(f)$. Considérons alors $u$ un relevé de $v$ dans $W(k')$, \textbf{$g$}$(t)=$\textbf{$f$}$(t+u)$ et $Y$ la courbe sur $W(k')$ correspondant à \textbf{$g$}, alors $Y_0\simeq \Spec(k')\times_{\Spec(k)}X_0$. \\
\end{pf}

Ainsi quitte à faire une extension finie de $k$ on suppose dans toute la suite que $t\wedge f(t)=1$.\\

\begin{notat*}
Afin d'alléger les notations du faisceau des formes différentielles de degré $1$ de $X_0$ sur $\mathbb{P}^1_k$ on notera désormais $\Omega^1_{X_0}=\Omega^1_{X_0/\Spec(k)}$. De même on notera $\Omega^1_{X_K}=\Omega^1_{X_K/\Spec(K)}$ et $\Omega^1_{X}=\Omega^1_{X/\Spec(V)}$. \\
\end{notat*}

\begin{conv} \label{conventions}
Soient $C_L$ une courbe lisse irréductible sur un corps $L$ et $Q$ un point de $C_L$. On suppose que $Q$ correspond à un idéal premier, encore noté $Q$, d'un ouvert affine $\Spec(D)\subset C_L$ et on désigne par $D_Q=\mathcal{O}_{C_L,Q}$ l'anneau local associé. Alors $D_Q$ est une anneau de valuation discrète d'uniformisante $\xi$, un paramètre local de $C_L$ en $Q$. Notons $L(Q)=\mathcal{O}_{C_L,Q}/\xi \mathcal{O}_{C_L,Q}$ le corps résiduel de $C_L$ en $Q$, on a un morphisme canonique $$\lambda_Q : D\rightarrow D/Q\rightarrow L(Q).$$
Soit $f\in D$ on utilisera le vocabulaire suivant
\begin{itemize}
\item $f(Q)$ pour $\lambda_Q(f)$ 
\item $f$ est inversible au voisinage de $Q$ pour $f$ est inversible dans l'anneau local $D_Q$
\item $v_{Q}$ pour la valuation $\xi$-adique de $D_Q$
\end{itemize}
par ailleurs on peut étendre $v_{Q}$ à $\Frac(D)$ le corps des fractions de la courbe $C_L$ et on se permettra de la noter encore $v_{Q}$. Enfin on note $\Omega^1_{C_L,Q}$ le localisé en $Q$ du module des différentielles $\Omega^1_{C_L}$ i.e. $\Omega^1_{C_L,Q}=\mathcal{O}_{C_L,Q}.d\xi$ et si $\omega\in \Omega^1_{C_L,Q}$ avec $\omega=h.d\xi$ on posera $v_{Q}(\omega)=v_{Q}(h)$. \\
\end{conv}

\begin{prop} \label{lisse}
La courbe relative $X$ est un schéma lisse sur $\Spec(W)$. De plus $U=D(y)\cup D($\textbf{$f'$}$)$.
\end{prop}

\begin{pf}
La démonstration se réalise en deux temps. On vérifie d'abord la lissité sur $U$ puis sur $V$. Soit donc $Q$ un idéal premier de $W[t,y]$ contenant $y^n-$\textbf{$f$}$(t)$ et notons $L(Q)$ le corps résiduel de $Q$ i.e. $L(Q)=\Frac(W[t,y]/Q)$. Considérons $\lambda : W[t,y]\rightarrow L(Q)$, définie via la surjection et l'injection canonique $W[t,y]\rightarrow W[t,y]/Q \rightarrow L(Q)$. En vertu du critère Jacobien (cf \cite{Neron} II.2 prop 7(d)), pour vérifier que la courbe est lisse, il suffit de vérifier que $\lambda(ny^{n-1}) \neq 0$ ou $\lambda($\textbf{$f'$}$(t))  \neq  0$. Comme $n$ est inversible sur $W$ il suffit donc de vérifier que $\lambda(y) \neq 0$ ou $\lambda($\textbf{$f'$}$(t))\neq  0$. Supposons que $\lambda(y)= 0$ alors $\lambda($\textbf{$f$}$(t))= 0$. Or \textbf{$f$} et \textbf{$f'$} sont premiers entre eux dans $W[t]$ par hypothèse donc il existe $a$ et $b$ dans $W[t]$ tels que $a$\textbf{$f$}$+b$\textbf{$f'$}$=1$ et alors on en déduit que $\lambda($\textbf{$f'$}$(t))  \neq  0$. \\
En ce qui concerne la lissité sur $V$, il suffit de la tester sur les points $Q$ de l'ouvert $V$ tels que $s=0$. Soit $L(Q)$ le corps résiduel de $V$ en $Q$, alors \textbf{$f_2$}$(s)=$\textbf{$f$}$(1/s)s^{l+1}=s$\textbf{$g$}$(s)$ où \textbf{$g$} est un élément de $W[s]$ qui vérifie \textbf{$g$}$(0)\neq 0$, puisque \textbf{$f$}$(0)$ est non nul, et donc la classe de \textbf{$g$} dans $L(Q)$ est non nulle. D'autre part, on a \textbf{$f'_2$}$(s)=s $\textbf{$g'$}$(s) + $\textbf{$g$}$(s) \in W[s]$
et donc dans le corps résiduel $L(Q)$ la classe de \textbf{$f'_2$} est égale à la classe de \textbf{$g$} dans $L(Q)$ et est non nulle. Ainsi la courbe $V$ est lisse en tout point $Q$ de $V(s)$ par le critère jacobien et $V$ est lisse.
\end{pf}

\begin{rem}
Les polynômes en deux variables $F(Y,T)$ et $F_2(Y,T)$ définis par $F(T,Y)=Y^n-f(T)$ et $F_2(Y,T)=Y^n-f_2(T)$, sont irréductibles. \\
\end{rem}

\begin{prop}
Pour $i$ inférieur ou égal à $2$ on a 
\[H^i_{DR}(X)\simeq H^i_{cris}(X_0).\]
\end{prop}

\begin{pf}
Comme $X$ est un schéma propre et lisse sur $\Spec(W)$ et que $X_0$ est sa fibre spéciale, ce résultat est une conséquence directe du théorème de comparaison de Berthelot \cite{Berth}. \\
\end{pf}

\subsection{Calcul du genre de la courbe $X_0$}

\begin{rem}
La proposition \ref{lisse} permet d'affirmer que la courbe $X_0$ définie précédemment est lisse sur $\Spec(k)$ et que $\mathcal{U}=D(y)\cup D(f')$. \\
\end{rem}

Pour calculer le genre, quitte à réaliser une extension de corps, on peut supposer que le polynôme $f$ est scindé sur $k$. Considérons donc des éléments $\alpha_i$ non nuls de $k$, tous distincts dans $k$ tels que 
\[f(t)=\prod_{i=1}^l (t-\alpha_i).\]
Rappelons que le polynôme $f_2(s)$ est défini par 
\[f_2(s)=f(1/s)s^{rn}=f(1/s)s^{l+1}=s\prod_{i=1}^l (1-\alpha_is)\]
et qu'il est de degré $l+1$. \\

\begin{prop}\label{ram}
La courbe $X_0$ est ramifiées en les $l+1$ points fermés correspondants aux idéaux maximaux 
\[M_i=(t-\alpha_i)\mathcal{O}_\mathcal{U}(\mathcal{U})+y\mathcal{O}_\mathcal{U}(\mathcal{U}) \ et \ \infty = s\mathcal{O}_\mathcal{V}(\mathcal{V})+z\mathcal{O}_{\mathcal{V}}(\mathcal{V})\]
de degré de ramification $n$. De plus les corps résiduels $k(M_i)$ et $k(\infty)$ sont égaux à $k$. 
\end{prop}

\begin{pf}
Nous ne démontrons pas ce résultat très classique.\\
\end{pf}

\begin{prop} \label{genre}
En conservant les notations précédentes, la courbe $X_0$ est de genre $g$ où \[g=\dfrac{(l-1)(n-1)}{2}=\dfrac{rn(n-1)}{2}-(n-1).\]
\end{prop}

\begin{pf}
Grâce à ce qui précède, on sait que la courbe $X_0$ est lisse (voir proposition \ref{lisse}) or une courbe lisse est normale. De plus on a vu que $\Pi:X_0\rightarrow \mathbb{P}^1_k$ est un morphisme fini séparable de degré $n$ (voir proposition \ref{pik}). Les conditions d'application du théorème d'Hurwitz sont donc vérifiées. On a vu à la proposition \ref{ram} que les degrés de ramifications des points de $X_0$ étaient $1$ ou $n$ mais $k$ est de caractéristique $p$ et $n$ est premier à $p$ il n'y a donc pas de ramification sauvage pour $X_0$. Pour un point $x$ de $X_0$ on a donc $e_x=1$ ou $e_x=n$ mais lorsque $e_x=n$ on a $[k(x):k]=1$ (toujours en vertu de la proposition \ref{ram}). Ainsi la formule d'Hurwitz donne 
\[2g(X_0)-2=n(2g(\mathbb{P}^1_k)-2)+ \sum_{\alpha_i} (e_{\alpha_i}-1).\] 
Or $\mathbb{P}^1_k$ est de genre nul et la somme contient $l+1$ termes, on a donc $2g(X_0)=(l+1)(n-1)-2(n-1)$ d'où les deux égalités voulues puisque par ailleurs $l+1=rn$. \\
\end{pf}

\begin{rem}
Si on fixe $n=2$ on retrouve bien $r=g+1$ valeur donnée par Liu (\cite{Liu} 7.4.3.). \\
\end{rem}

\subsection{Une action sur la courbe $X_0$}\label{action}

On suppose dans cette partie que $k$ contient les racines $n$-ièmes de l'unité de $\overline{\mathbb{F}_p}$, une clôture algébrique de $\mathbb{F}_p$, i.e. $\mu_n(\overline{\mathbb{F}}_p)$ est inclus dans $k$. Sous cette hypothèse, le polynôme $Q_n(X)=X^n-1$ admet $n$ racines distinctes dans $k$. Par ailleurs puisqu'on a supposé que $n$ et $p$ étaient premiers entre eux, si $Q_n(\alpha)$ est nul alors $Q'_n(\alpha)=n\alpha^{n-1}$ est non nul. Enfin puisque $k$ est un corps fini de caractéristique $p$, $k$ est isomorphe à $\mathbb{F}_q$, où $q=p^m$ pour un certain entier non nul $m$. Dans ce cas, $\mu_n(\overline{\mathbb{F}}_p)$ est inclus dans $k$ si et seulement si $n$ divise $p^m-1$. Soit donc $\mu_n(\overline{\mathbb{F}}_p)$ le groupe, cyclique, des racines $n$-ièmes de l'unité. Soit $\zeta\in \mu_n(\overline{\mathbb{F}}_p)$ une racine primitive $n$-ième de l'unité. \\

On définit $\delta$ sur $k[t,y]$ par $$\delta : (t,y)\mapsto (t,\zeta^{-1} y)$$ 
qui passe au quotient en 
\[ 
\begin{array}{c c c c}
    \delta : & \mathcal{O}_{\mathcal{U}}(\mathcal{U}) & \longrightarrow & \mathcal{O}_{\mathcal{U}}(\mathcal{U}) \\
     & t & \longmapsto & t \\
     & y & \longmapsto &  \zeta ^{-1}y
\end{array}
\]
de même on définit $\delta$ sur $k[s,z]$ par $$\delta : (s,z)\mapsto (s, \zeta^{-1} z)$$ 
qui passe au quotient en
\[ 
\begin{array}{c c c c}
    \delta : & \mathcal{O}_{\mathcal{V}}(\mathcal{V}) & \longrightarrow & \mathcal{O}_{\mathcal{V}}(\mathcal{V})\\
     & s & \longmapsto & s \\
     & z & \longmapsto & \zeta ^{-1}z
\end{array}
\]
et ces deux applications se recollent sur $\mathcal{U}\cap \mathcal{V}$. Donc $\delta$ induit un morphisme de $X_0$ dans $X_0$ qui décrit l'action de $\mu_n(\overline{\mathbb{F}}_p)$ sur $X_0$. Par fonctorialité $\delta$ induit un morphisme sur $\Omega^1_{X_0}$, noté $\delta^*$. Cette action est déduite de la précédente par différentiation, on a donc la description suivante
\[ 
\begin{array}{c c c c}
    \delta^* : & \Omega^1_{X_0}(\mathcal{U}) & \longrightarrow & \Omega^1_{X_0}(\mathcal{U})\\
     & dt & \longmapsto & dt \\
     & dy & \longmapsto &  \zeta ^{-1}dy 
\end{array}
\]
\newline
\[ 
\begin{array}{c c c c}
    \delta^* : & \Omega^1_{X_0}(\mathcal{V}) & \longrightarrow & \Omega^1_{X_0}(\mathcal{V})\\
     & ds & \longmapsto & ds \\
     & dz & \longmapsto &  \zeta ^{-1}dz
\end{array}
\]
\newline

\section{Bases des différents espaces de cohomologie}

\begin{notat*}
On appelle point à l'infini le point de $\mathcal{V}$ correspondant à l'idéal engendré par $s$ et $z$ dans $k[s,z]/(z^n-f_2(s))$. Ceci définit un diviseur que l'on notera $\infty$ \\
\end{notat*}

\begin{defi} \label{aucuneidee}
Soit $\sigma_k$ le Frobenius sur $k$ (i.e. l'élévation à la puissance $p$), alors $\sigma_k$ induit un morphisme de schémas $\Spec(k)\rightarrow \Spec(k)$. Si $Y$ est un schéma sur $\Spec(k)$ on notera $Y'_0$ le schéma $\Spec(k) \times_{\Spec(k)} Y$, où le produit fibré est pris au-dessus de $\sigma_k$. \\
\end{defi}

Dans ce qui suit nous donnons une description de différents groupes de cohomologie attachés à $X_0$ et $X'_0$. \\

\subsection{Une base de $H^0(X_0,\Omega^1_{X_0})$}

\begin{prop} \label{baseH0k}
Le $k$-espace vectoriel $H^0(X_0,\Omega_{X_0}^1)$ est de dimension $g$ et admet pour base la famille 
\begin{center}
$\omega_{i,j}=\dfrac{t^{i}dt}{y^j}$ avec $1\leq j\leq n-1$ et $0\leq i\leq rj-2$
\end{center}
où pour mémoire $r=(l+1)/n$. 
\end{prop}

\begin{pf} Comme $X_0$ est propre et lisse (proposition \ref{lisse}) alors le faisceau $\Omega_{X_0}^1$ est localement libre et $H^0(X_0,\Omega_{X_0}^1)$ est un $k$-espace vectoriel de dimension finie. Par ailleurs le fait qu'il soit de dimension $g$ est un résultat découlant directement du théorème de Riemann-Roch. Montrer que la famille $\omega_{i,j}$ considérée est une famille d'éléments de $H^0(X_0,\Omega_{X_0}^1)$ revient à montrer que pour tout $Q$ point de $X_0$, $v_Q(\omega_{i,j})\geq 0$ (avec les conventions \ref{conventions}). Par ailleurs on a vu que $\mathcal{U}=D(y)\cup D(f')$ et donc qu'un point de $X_0$ est soit un point de $D(y)$, soit un point de ramification (contenu dans $D(f')$). Il est clair que $\omega_{i,j} \in \Omega^1_{X_0}(\mathcal{U}\cap D(y))$ (sans condition sur $j$ et avec $i$ supérieur ou égal à $0$) il n'y a donc plus qu'à considérer le cas où $Q$ est un point de ramification de la courbe. Si $Q \in W(y)\subset D(f')$ est un point de ramification de $\mathcal{U}$, i.e si $\Pi(Q)=(t-\alpha_i)k[t]$ avec $\alpha\in \{ \alpha_1,...,\alpha_l \}$, alors un paramètre local en $Q$ est $y$ il faut donc commencer par exprimer $\omega_{i,j}$ selon $dy$. Or on a la relation $y^n=f(t)$ donc $ny^{n-1}dy=f'(t)dt$ et on peut exprimer $\omega_{i,j}$ sur $D(f')$ par la formule 
\[\omega_{i,j}=\dfrac{t^iy^{n-1-j}}{f'(t)}dy\]
ainsi $v_Q(\omega_{i,j})=iv_Q(t)+(n-1-j)v_Q(y)-v_Q(f'(t))=iv_Q(t)+n-1-j$. En effet $t=t-\alpha_i+\alpha_i$ donc $f'(t)$ est inversible dans $A_Q$. De plus $\alpha_i$ est non nul donc $v_Q(t)$ vaut $0$ et $v_Q(\omega_{i,j})=n-1-j$ qui est supérieur ou égal à $0$ dès que $j$ est inférieur ou égal à $n-1$. Si maintenant $Q$ est le point à l'infini défini précédemment, $\Pi(Q)=sk[s]$ et un paramètre local en $Q$ est $z$. Il faut donc exprimer $\omega_{i,j}$ selon $s$, $z$ et $dz$. On utilise les formules $t=1/s$ et $y=s^{-r}z$ où pour mémoire $r=(l+1)/n$. On a alors
\[\omega_{i,j}=-\dfrac{s^{rj-2-i}}{z^j}ds=-ns^{rj-2-i}z^{n-1-j}dz\]
et $v_Q(\omega_{i,j})=n(rj-2-i)+(n-1-j)$. Mais comme $i$ est inférieur ou égal à $rj-2$ et $j$ inférieur ou égal à $n-1$ alors $v_Q(\omega_{i,j})$ est supérieur ou égal à $0$, ce qui achève de montrer que la famille considérée est une famille de $H^0(X_0,\Omega_{X_0}^1)$. Pour montrer qu'elle est linéairement indépendante sur $k$, il suffit de montrer qu'elle l'est en restriction à $\mathcal{U}\cap D(y)$. Or $\Omega^1_{X_0}(\mathcal{U}\cap D(y))$ est un $\mathcal{O}_{X_0}(\mathcal{U}\cap D(y))$-module libre de base $dt$. Plus précisément 
\[\Omega^1_{X_0}(\mathcal{U}\cap D(y))\simeq \bigoplus_{i=0}^{l-1}k[y,y^{-1}]t^idt.\] 
Alors en restriction à $\mathcal{U}\cap D(y)$ les $\omega_{i,j}$ forment un $k$-module libre et la famille est donc libre dans $H^0(X_0,\Omega_{X_0}^1)$. Pour conclure il suffit de montrer qu'elle contient $g$ éléments pour montrer que c'est une base de $H^0(X_0,\Omega^1_{X_0})$. Or elle contient :
\[ \sum_{j=1}^{n-1} (rj-1)=r\sum_{j=1}^{n-1} j - (n-1)=r\dfrac{(n-1)n}{2}-(n-1)=g\]
éléments, c'est donc bien une base du $k$-espace vectoriel considéré. \\
\end{pf}

\begin{lem} \label{lemmepourprojeteromega}
On conserve les notations du théorème précédent. Alors \\
(i) $\displaystyle \Omega^1_{X_0}(\mathcal{U}\cap D(y))=H^0(X_0,\Omega_{X_0}^1)\oplus \left( \bigoplus_{i=0}^{l-1} k[y]t^i dt \oplus \bigoplus_{j\geq n}\bigoplus_{i=0}^{l-1} ky^{-j}t^i dt \oplus \bigoplus_{j=1}^{n-1}\bigoplus_{i=rj-1}^{l-1} ky^{-j}t^i dt\right)$ \\
(ii) $\displaystyle \Omega^1_{X_0}(\mathcal{V})= H^0(X_0,\Omega_{X_0}^1)\oplus \left( k[s]ds \oplus \bigoplus_{j=1}^{n-1}\bigoplus_{i\geq rj-1} ks^iz^{-j}ds\right) $ 
\end{lem}

\begin{pf}
A la fin de la démonstration de la proposition \ref{baseH0k} on a établi que 
\[\Omega^1_{X_0}(\mathcal{U}\cap D(y))\simeq \bigoplus_{i=0}^{l-1}k[y,y^{-1}]t^idt\] 
ce qui donne immédiatement la décomposition en somme directe du (i). En ce qui concerne le résultat sur le faisceau $\Omega^1_{{X_0}}(\mathcal{V})$ commençons par montrer qu'il est libre, engendré par 
\[-\omega_{r(n-1)-2,n-1}=-\dfrac{t^{r(n-1)-2}dt}{y^{n-1}}=\dfrac{ds}{z^n}\]
la dernière égalité étant évidente au vu des conditions de recollement. Dans $\Omega^1_{{X_0}}(\mathcal{V})$ on a la relation $nz^{n-1}dz-f'_2(s)ds=0$, de plus comme ${X_0}$ est lisse on a $\mathcal{V}=D(z)\cup D(f'_2(s))$. Il est clair que $\omega_{r(n-1)-2,n-1}$ est un élément de $\Omega^1_{{X_0}}(\mathcal{V}\cap D(z))$ et que $\Omega^1_{{X_0}}(\mathcal{V}\cap D(z))$ est engendré par 
\[ds=z^{n-1}(-\omega_{r(n-1)-2,n-1})\]
donc aussi par $(-\omega_{r(n-1)-2,n-1})$ le faisceau $\Omega^1_{{X_0}}(\mathcal{V}\cap D(f'_2(s)))$ est quant à lui engendré par 
\[dz=\dfrac{-1}{n}f'_2(s)\omega_{r(n-1)-2,n-1}\]
(où la division par $n$ est licite car $n$ est premier à $p$) et donc aussi par $(-\omega_{r(n-1)-2,n-1})$ ce qui montre notre affirmation. De plus on a 
\[\mathcal{O}_{X_0}(\mathcal{V})=k[s,z]/(z^n-f_2(s))=\bigoplus_{j=0}^{n-1} k[s]z^j\]
et comme $\Omega_{X_0}^1(\mathcal{V})$ est libre engendré par $(-\omega_{r(n-1)-2,n-1})$ on a
\[\Omega^1_{X_0}(\mathcal{V})=\bigoplus_{j=0}^{n-1} k[s]z^j\dfrac{ds}{z^{n-1}}=\bigoplus_{j=0}^{n-1} k[s]\dfrac{ds}{z^{j}}\]
mais on a montré que
\[H^0({X_0},\Omega_{X_0}^1)=\bigoplus_{j=1}^{n-1}\bigoplus_{i=0}^{rj-2} k\dfrac{t^idt}{y^j}=\bigoplus_{j=1}^{n-1}\bigoplus_{i=0}^{rj-2} k\dfrac{-s^{rj-2-i}ds}{z^j}=\bigoplus_{j=1}^{n-1}\bigoplus_{i=0}^{rj-2} k\dfrac{s^{i}ds}{z^j}\]
ce qui donne bien la décomposition en somme directe annoncée. \\
\end{pf}

\begin{prop} \label{basesprimeesH0}
Une base de $H^0(X'_0,\Omega^1_{X'_0})$ est donnée par $\omega'_{i,j}=t^iy^{-j}dt$ pour $1\leq j\leq n-1$ et $0\leq i\leq rj-2$.
\end{prop}

\begin{pf}
On dispose déjà d'une base de $H^0(X_0,\Omega^1_{X_0})$ donnée par la proposition précédente, dont les éléments sont $\omega_{i,j}=t^iy^{-j}dt$ pour $1\leq j\leq n-1$ et $0\leq i\leq rj-2$. Au vu de la démonstration et par définition de $X'_0$, il est clair que l'on obtient bien la base proposée pour $H^0(X'_0,\Omega^1_{X'_0})$. \\
\end{pf}

\subsection{Décomposition de $H^0({X_0},\Omega^1_{X_0})$ en composantes isotypiques}

\begin{theo} \label{SDDiH0}
Avec les notations précédentes, on a 
\[H^0({X_0},\Omega_{{X_0}}^1)=\bigoplus_{j=1}^{n-1} \bigoplus_{i=0}^{rj-2} k\dfrac{t^{i}dt}{y^j}.\] 
Supposons maintenant que $\mu_n(\overline{\mathbb{F}}_p)\subset k$, alors si on note $\mathcal{D}_j=H^0({X_0},\Omega_{{X_0}}^1)_j$ la $j$-ième composante isotypique associée à l'action $\delta$ décrite dans le paragraphe \ref{action} on a : 
\[\mathcal{D}_0=0\]
\[\mathcal{D}_j=\bigoplus_{i=0}^{rj-2} k\dfrac{t^{i}dt}{y^j} \ \ \forall 1\leq j\leq n-1.\]
\end{theo}

\begin{pf}
C'est une conséquence immédiate du théorème \ref{baseH0k} et de la définition de l'action $\delta$ définie à la section \ref{action}. \\
\end{pf}

\subsection{Une base de $H^1({X_0},\mathcal{O}_{X_0})$}

On considère le complexe de \v{C}ech relatif au recouvrement ${X_0} = \mathcal{U}\cup \mathcal{V}$.
$$0\longrightarrow \mathcal{O}_{X_0}(\mathcal{U})\times \mathcal{O}_{X_0}(\mathcal{V})\overset{d_{0,1}}{\longrightarrow} \mathcal{O}_{X_0}(\mathcal{U}\cap \mathcal{V}) \longrightarrow 0$$
où l'application $d_{0,1}$ est donnée par $(a,b)\mapsto (a-b)\vert_{\mathcal{U}\cap \mathcal{V}}$. Comme $H^1({X_0},\mathcal{O}_{X_0})$ est égal au premier espace de cohomologie de ce complexe de \v{C}ech, il est égal à $\mathcal{O}_{X_0}(\mathcal{U}\cap \mathcal{V})/Im(d_{0,1})$. \\

\begin{prop} \label{baseH1Ok}
Le $k$-espace vectoriel $H^1({X_0},\mathcal{O}_{{X_0}})$ est de dimension $g$ et admet pour base les classes $\overline{h_{i,j}}$ des éléments $h_{i,j}\in  \mathcal{O}_{{X_0}}(\mathcal{U} \cap \mathcal{V} )$ où  
\begin{center}
$h_{i,j}=\dfrac{y^j}{t^i}=\dfrac{z^j}{s^{rj-i}}$, $1\leq j\leq n-1$, $1\leq i\leq rj-1$
\end{center}
\end{prop}

\begin{pf}
Comme ${X_0}$ est propre on sait que $H^1({X_0},\mathcal{O}_{X_0})$ est un $k$-espace vectoriel de dimension finie. Par ailleurs par dualité de Poincaré on a $H^1({X_0},\mathcal{O}_{X_0})\simeq H^0({X_0},\Omega_{X_0}^1)$. Mais on a établi au théorème \ref{baseH0k} que ce dernier espace vectoriel était de dimension $g$ ce qui donne le résultat sur $H^1({X_0},\mathcal{O}_{X_0})$. De plus $\mathcal{O}_{X_0}(\mathcal{U})=k\left[ t,y\right] / (y^n-f(t))$ et $\mathcal{O}_{X_0}(\mathcal{V})=k\left[ s,z\right] / (z^n-f_2(t))$ et donc
\[\mathcal{O}_{X_0}(\mathcal{U}\cap \mathcal{V})=k\left[ t,y,t^{-1}\right] / (y^n-f(t))=k\left[ s,z,s^{-1}\right] / (z^n-f_2(t)).\]
Une base de $\mathcal{O}_{X_0}(\mathcal{U}\cap \mathcal{V})$ comme $k$-espace vectoriel est donc $t^iy^j$ avec $i\in \mathbb{Z}$ et $0\leq j\leq n-1$ ou encore $s^iz^j$ avec $i\in \mathbb{Z}$ et $0\leq j\leq n-1$. De même une base de $\mathcal{O}_{X_0}(\mathcal{U})$ (respectivement $\mathcal{O}_{X_0}(\mathcal{V})$) est $t^iy^j$ avec $i\geq 0$ et $0\leq j\leq n-1$ (respectivement $s^iz^j$ avec $i\geq 0$ et $0\leq j\leq n-1$). Par définition de $d_{0,1}$, tout élément $t^iy^j$ avec $i\geq 0$ et $0\leq j\leq n-1$ est dans l'image de $d_{0,1}$. De la même manière tout élément du type $s^iz^j$ avec $i\geq 0$ et $0\leq j\leq n-1$ est dans l'image de $d_{0,1}$. Or sur $\mathcal{U}\cap \mathcal{V}$ on a $s=1/t$ et $z=s^ry$ d'où $s^iz^j=t^{-(i+rj)}y^j$. Par ailleurs lorsque $j$ vaut $0$, $i+rj$ prends toutes les valeurs entières lorsque $i$ est dans $\mathbb{N}$ et lorsque $j$ est non nul, $i+rj$ prend toutes les valeurs supérieures ou égales à $rj$. Posons donc $B=\left\lbrace kt^{i}y^j, \ 0\leq j\leq n-1, i\in \mathbb{N} \right\rbrace \cup \left\lbrace kt^{-i}y^j, \ 0\leq j\leq n-1, i\geq rj \right\rbrace$, alors $Im(d)\supset B$. On sait également grâce à la description de $\mathcal{O}_{X_0}(\mathcal{U}\cap \mathcal{V})$ donnée plus haut que $H^1({X_0},\mathcal{O}_{X_0})=\sum k\overline{t^{i}y^j}$ pour $0\leq j\leq n-1$ et $i$ élément de $\mathbb{Z}$ et les considérations précédentes permettent alors de dire que 
\[ H^1({X_0},\mathcal{O}_{{X_0}})=\sum_{\substack{1\leq j\leq n-1 \\ 1\leq i\leq rj-1}} k\overline{t^{-i}y^j}.  \ \ \ (1)\]
Notons $e_{i,j}=\overline{t^{-i}y^j}$ pour $1\leq j\leq n-1$ et $1\leq i\leq rj-1$, $L=\oplus ke_{i,j}$ le $k$-espace vectoriel et $\varphi : L\rightarrow H^1({X_0},\mathcal{O}_{X_0})$ l'application qui envoie $e_{i,j}$ sur $\overline{h_{i,j}}$. Remarquons tout d'abord que $L$ est de dimension $g$ puisque le nombre d'éléments $e_{i,j}$ est 
\[\sum_{j=1}^{n-1} (rj-1)=r\sum_{j=1}^{n-1} j - (n-1)=r\dfrac{(n-1)n}{2}-(n-1)=g.\]
On sait que $\varphi$ est surjective grâce à l'égalité $(1)$, c'est donc une application $k$-linéaire surjective entre deux $k$-espaces vectoriels de même dimension, et $\varphi$ induit un isomorphisme $L\simeq H^1({X_0},\mathcal{O}_{X_0})$. En particulier les éléments $\overline{h_{i,j}}$ forment une base de $H^1({X_0},\mathcal{O}_{X_0})$. \\
\end{pf}

\begin{prop} \label{basesprimeesH1O}
Une base de $H^1(X'_0,\mathcal{O}_{X'_0})$ est donnée par $\overline{h'_{i,j}}=\overline{t^{-i}y^j}$ pour $1\leq j\leq n-1$ et $1\leq i\leq rj-1$.
\end{prop}

\begin{pf}
On dispose d'une base de $H^1(X_0,\mathcal{O}_{X_0})$ par la proposition précédente, dont les éléments sont $\overline{h_{i,j}}=t^{-i}y^j$ pour $1\leq j\leq n-1$ et $1\leq i\leq rj-1$. Au vu de la démonstration et par définition de $X'_0$, il est clair que l'on obtient bien la base proposée pour $H^1(X'_0,\mathcal{O}_{X'_0})$. \\
\end{pf}

\subsection{Décomposition de $H^1({X_0},\mathcal{O}_{X_0})$ en composantes isotypiques}

\begin{prop}\label{SDDiH1}
Avec les notations précédentes, on a
\[H^1({X_0},\mathcal{O}_{{X_0}})=\displaystyle \bigoplus_{j=1}^{n-1}\ \bigoplus_{i=1}^{rj-1} k\overline{t^{-i}y^j}=\bigoplus_{j=1}^{n-1} \ \bigoplus_{i=1}^{r(n-j)-1} k\overline{t^{-i}y^{n-j}}.\]
Supposons maintenant que $\mu_n(\overline{\mathbb{F}}_p)\subset k$, alors si on note $H^1({X_0},\mathcal{O}_{{X_0}})_j$ la $j$-ième composante isotypique associée à l'action de $\delta$ on a :
\[H^1({X_0},\mathcal{O}_{{X_0}})_0=0\]
\[H^1({X_0},\mathcal{O}_{{X_0}})_j=\bigoplus_{i=1}^{r(n-j)-1} k\overline{t^{-i}y^{n-j}} \ \ \forall 1\leq j\leq n-1.\] 
\end{prop}

\begin{pf}
C'est une conséquence immédiate de la proposition \ref{baseH1Ok} et de la définition de l'action $\delta$ définie à la section \ref{action}. \\
\end{pf}

\section{Un scindage adapté de la filtration de Hodge}\label{calculscindage}

\begin{rem} \label{identH1Ok}
Au vu de la démonstration de la proposition \ref{baseH1Ok}, il y a une identification naturelle d'un élément $\overline{h_{i,j}}$ de la base de $H^1({X_0},\mathcal{O}_{X_0})$ avec le représentant $h_{i,j}$ élément de $\mathcal{O}_{\mathcal{U}\cap \mathcal{V}}$. Par la suite on notera simplement $t^{-i}y^j$ le représentant de la classe $\overline{t^{-i}y^j}\in H^1({X_0},\mathcal{O}_{X_0})$ lorsqu'on voudra le regarder comme élément de $\mathcal{O}_{\mathcal{U}\cap \mathcal{V}}$.\\
\end{rem}

Notons $\Omega_{X_0}^{\bullet}$ le complexe de de Rham associé à $\Omega_{X_0}^1$ et $\mathcal{F}$ le recouvrement de ${X_0}$ par $\mathcal{U}$ et $\mathcal{V}$ alors on note \v{C}$(\mathcal{F},\Omega_{X_0}^{\bullet})$ le complexe de \v{C}ech associé au complexe de de Rham. Il est donné par
$$\begin{array}{ccccccc}
    & & 0 & & 0 & & \\
    & & \uparrow & & \uparrow & & \\
    0 & \rightarrow & \Omega_{X_0}^1(\mathcal{U}) \times \Omega_{X_0}^1(\mathcal{V}) & \rightarrow & \Omega_{X_0}^1(\mathcal{U}\cap \mathcal{V}) & \rightarrow & 0 \\
    & & \uparrow & & \uparrow & & \\
    0 & \rightarrow & \mathcal{O}_{X_0}(\mathcal{U}) \times  \mathcal{O}_{X_0}(\mathcal{V}) & \rightarrow & \mathcal{O}_{X_0}(\mathcal{U}\cap \mathcal{V}) & \rightarrow & 0 \\
    & & \uparrow & & \uparrow & & \\
    & & 0 & & 0 & & \\
\end{array}$$
on peut lui associer le complexe simple
\[0 \rightarrow \mathcal{O}_{X_0}(\mathcal{U}) \times  \mathcal{O}_{X_0}(\mathcal{V}) \rightarrow \Omega_{X_0}^1(\mathcal{U}) \times \Omega_{X_0}^1(\mathcal{V})\oplus \mathcal{O}_{X_0}(\mathcal{U}\cap \mathcal{V}) \rightarrow \Omega_{X_0}^1(\mathcal{U}\cap \mathcal{V}) \rightarrow 0\]
et le groupe de cohomologie $H^1_{DR}({X_0})$ s'identifie alors aux classes des éléments 
\[(\omega_\mathcal{U},\omega_\mathcal{V},h)\in \Omega^1_{X_0}(\mathcal{U})\oplus \Omega^1_{X_0}(\mathcal{V})\oplus \mathcal{O}_{X_0}(\mathcal{U}\cap \mathcal{V})\]
tels que $\omega_\mathcal{V}-\omega_\mathcal{U}+dh=0$, modulo les éléments de la forme $(dh_\mathcal{U},-dh_\mathcal{V},h)$ où $h=h_\mathcal{U}+h_\mathcal{V}$, $h_\mathcal{U}\in \mathcal{O}_{X_0}(\mathcal{U})$, $h_\mathcal{V}\in \mathcal{O}_{X_0}(\mathcal{V})$. Dans la suite on notera $H^1($\v{C}$(\mathcal{F},\Omega_{X_0}^{\bullet}))$ le premier espace de cohomologie associé au complexe simple considéré ici. \\

\begin{prop} \label{scindagemodp}
On conserve les notations précédentes, alors : \\
(i) L'image de $\omega \in H^0({X_0},\Omega^1_{X_0})$ dans $H^1($\v{C}$(\mathcal{F},\Omega_{X_0}^{\bullet}))$ est la classe de l'élément $(\omega\vert_{\mathcal{U}},\omega\vert_{\mathcal{V}},0)$. \\
(ii) Soient $1\leq j\leq n-1$ et $1\leq i\leq rj-1$. Notons $\lambda_k$ le coefficient de $t^k$ dans $f(t)$ et posons $a_j=(n-j)r-2$, 
\[\alpha_{i,j}^\mathcal{U} = \left[ \displaystyle \sum_{k=i+1}^l \dfrac{jk-in}{n}\lambda_kt^{k-(i+1)}\right] y^{j-1}\dfrac{dt}{y^{n-1}},\]
\[\alpha_{i,j}^\mathcal{V} = \left[ \sum_{k=0}^i \dfrac{in-jk}{n} \lambda_ks^{i+1+a_j-k}\right]z^{j-1}\dfrac{ds}{z^{n-1}}.\]
Alors $\alpha_{i,j}^\mathcal{U}\in \Omega^1_{X_0}(\mathcal{U})$, $\alpha_{i,j}^\mathcal{V}\in \Omega^1_{X_0}(\mathcal{V})$, $dh_{i,j}=\alpha_{i,j}^\mathcal{U}+\alpha_{i,j}^\mathcal{V}$ et on a une application $k$-linéaire :
\[ 
\begin{array}{c c c c}
    [\tau] :  & H^1({X_0},\mathcal{O}_{X_0}) & \longrightarrow & H^1($\v{C}$(\mathcal{F},\Omega_{X_0}^{\bullet}))\\
     & \overline{h_{i,j}} & \longmapsto & \left[ (\alpha_{i,j}^\mathcal{U},-\alpha_{i,j}^\mathcal{V},h_{i,j})\right] \\
\end{array}
\]
qui est un scindage de la filtration de Hodge. Les éléments $\alpha_{i,j}$ sont déterminés à une section globale $\omega \in H^0({X_0},\Omega^1_{X_0})$ près. 
\end{prop}

\begin{pf}
(i) Cela résulte de l'inclusion naturelle $H^0({X_0},\Omega^1_{X_0}) \subset H^1($\v{C}$(\mathcal{F},\Omega_{X_0}^{\bullet}))\simeq H^1_{DR}({X_0})$. \\
(ii) Rappelons tout d'abord que $dt/y^{n-1}\in H^0({X_0},\Omega^1_{X_0})$ (voir théorème \ref{baseH0k}) et de même pour $ds/z^{n-1}$. Soient $1\leq j\leq n-1$ et $1\leq i\leq rj-1$ fixés, calculons $dh_{i,j}$,
\[dh_{i,j}=\dfrac{jy^{j-1}dy}{t^i}-\dfrac{iy^jdt}{t^{i+1}} \ \ mod \ \mathcal{O}_{X_0}(\mathcal{U}\cap \mathcal{V})\left[ny^{n-1}dy-f'(t)dt \right]. \]
Par ailleurs on a la relation $y^n=f(t)$ dans $\mathcal{O}_{X_0}(\mathcal{U}\cap \mathcal{V})$ et donc
\[dh_{i,j}=\dfrac{jf'(t)y^{j-1}}{nt^i}\dfrac{dt}{y^{n-1}}-\dfrac{iy^jy^{n-1}}{t^{i+1}}\dfrac{dt}{y^{n-1}}=\left( \dfrac{jf'(t)}{nt^i}-\dfrac{if(t)}{t^{i+1}}\right) y^{j-1}\dfrac{dt}{y^{n-1}} \in \Omega^1_{X_0}(\mathcal{U}\cap \mathcal{V}).\]
On rappelle que 
\[\Omega^1_{X_0}(\mathcal{U})\supset \left( k[t,y]/(y^n-f(t))\right) \dfrac{dt}{y^{n-1}}\]
\[\Omega^1_{X_0}(\mathcal{V})\supset \left( k[s,z]/(z^n-f_2(s))\right) \dfrac{ds}{z^{n-1}}.\]
donc pour montrer le résultat voulu montrons que l'on peut écrire 
\[dh_{i,j}= R_{i,j}(t,y)\dfrac{dt}{y^{n-1}}+S_{i,j}(s,z)\dfrac{ds}{z^{n-1}}\] 
avec $R_{i,j}$ élément de $k[y,t]$ et $S_{i,j}$ élément de $k[z,s]$, $\forall 1\leq j\leq n-1$, $\forall 1\leq i\leq rj-1$. Remarquons qu'avec ces conditions $i$ est strictement inférieur à $l$. En effet
\[i\leq rj-1 \Rightarrow i\leq r(n-1)-1 =rn-1-r=l-r<l.\] 
Notons $f(t)=\sum \lambda_k t^k$ alors $f'(t)=\sum k\lambda_k t^{k-1}$ (pour mémoire $\lambda_k\in k$), alors
\[dh_{i,j}=\left[ \displaystyle \dfrac{j}{n}\dfrac{1}{t^i}\sum_{k=1}^l k\lambda_kt^{k-1}-i\dfrac{1}{t^{i+1}}\sum_{k=0}^l \lambda_kt^k\right] y^{j-1}\dfrac{dt}{y^{n-1}} \in \Omega^1_{X_0}(\mathcal{U}\cap \mathcal{V})\]
\begin{center}
$\displaystyle dh_{i,j}= \left[  \dfrac{j}{n}\sum_{k=1}^i k\lambda_kt^{k-(i+1)}-i\sum_{k=0}^i \lambda_kt^{k-(i+1)}\right] y^{j-1}\dfrac{dt}{y^{n-1}} +\left[  \sum_{k=i+1}^l \left( \dfrac{j}{n}k-i\right) \lambda_kt^{k-(i+1)} \right] y^{j-1}\dfrac{dt}{y^{n-1}}.$
\end{center}
Notons $\alpha_{i,j}^\mathcal{U}$ le second terme de cette expression. C'est un élément de $(k[y,t]/(y^n-f(t)))y^{-(n-1)}dt$. Pour conclure montrons que le premier terme, noté $\alpha_{i,j}^\mathcal{V}$, est un élément de $(k[z,s]/(z^n-f_2(s)))z^{-(n-1)}ds$.
\[y^{j-1}\dfrac{dt}{y^{n-1}}=\dfrac{-s^{-(j-1)r}z^{j-1}s^{-2}ds}{s^{-(n-1)r}z^{n-1}}=-s^{(n-j)r-2}z^{j-1}\dfrac{ds}{z^{n-1}}\in \Omega^1_{X_0}(\mathcal{U}\cap \mathcal{V}).\]
Notons $a_j=(n-j)r-2$, comme $1\leq j\leq n-1$ alors $r-2\leq a_j \leq (n-1)r-2$. Mais $r=(l+1)/n$ avec $l$ congru à $-1$ modulo $n$ et $l$ positif ou nul. Ainsi $r\geq 1$ et $-1\leq a_j \leq (n-1)r-2$.
\[\alpha_{i,j}^\mathcal{V} = \left[  -\dfrac{j}{n}\sum_{k=1}^i k\lambda_ks^{i+1+a_j-k}+i\sum_{k=0}^i \lambda_ks^{i+1+a_j-k}\right] z^{j-1}\dfrac{ds}{z^{n-1}}.\]
Comme $i-k\geq 0$ et que par les considérations précédentes $1+a_j\geq 0$, alors $\alpha_{i,j}^\mathcal{V}$ est bien dans l'espace désiré. Par linéarité cela définit une application $\tau$ qui à $h_{i,j}$ associe $(\alpha_{i,j}^\mathcal{U},-\alpha_{i,j}^\mathcal{V},h_{i,j})$ dans $\Omega_{X_0}^1(\mathcal{U})\oplus \Omega_{X_0}^1(\mathcal{V})\oplus \mathcal{O}_{X_0}(\mathcal{U}\cap \mathcal{V})$ qui donne $[\tau]$ par passage au quotient dans $H^1($\v{C}$(\mathcal{F},\Omega_{X_0}^{\bullet}))$ ce qui achève la démonstration. \\
\end{pf}

Décrivons maintenant l'action du scindage $\tau$ sur les composantes isotypiques associées à l'action $\delta$, définie à la section \ref{action}, des différents modules. \\

\begin{prop} \label{isotau}
On conserve les notations de la proposition précédente. Supposons que $\mu_n(\overline{\mathbb{F}_p})$ soit inclus dans $k$, alors $\tau$ envoie la $j$-ième composante isotypique de  $H^1({X_0},\mathcal{O}_{X_0})$ sur la $j$-ième composante isotypique de $\Omega_{X_0}^1(\mathcal{U})\oplus \Omega_{X_0}^1(\mathcal{V})\oplus \mathcal{O}_{X_0}(\mathcal{U}\cap \mathcal{V})$. De même $[\tau]$ envoie la $j$-ième composante isotypique de  $H^1({X_0},\mathcal{O}_{X_0})$ sur la $j$-ième composante isotypique de $H^1($\v{C}$(\mathcal{F},\Omega_{X_0}^{\bullet}))$. 
\end{prop}

\begin{pf}
Soit $1\leq j\leq n-1$ fixé, et soit $\tau$ comme dans la proposition précédente. Alors l'élément $\overline{h_{i,j}}$ est envoyé par $\tau$ sur un triplet du type $(\alpha_\mathcal{U}(t)y^{-(n-j)}dt,\alpha_\mathcal{V}(t)y^{-(n-j)}dt,h(t)y^j)$. Comme l'action de $\delta$ est donnée par $\delta (t)=t$ et $\delta (y)=\zeta^{-1}y$ (voir section \ref{action}), l'image de $\overline{h_{i,j}}=\overline{t^{-i}y^j}$, élément de la $(n-j)$-ième composante isotypique de $H^1({X_0},\mathcal{O}_{X_0})$, est dans la $(n-j)$-ième composante isotypique de $\Omega_{X_0}^1(\mathcal{U})\oplus \Omega_{X_0}^1(\mathcal{V})\oplus \mathcal{O}_{X_0}(\mathcal{U}\cap \mathcal{V})$ ce qui donne bien le résultat annoncé. Le résultat sur $[\tau]$ s'en déduit par passage aux classes d'équivalence. \\
\end{pf}

\section{Mise en place du calcul de Frobenius divisé}

\subsection{Rappels sur le Frobenius}\label{rappels}

\begin{defi} \label{polysigma}
Soit $\sigma$ un relèvement du Frobenius sur $\Spec(W_1)$, alors $\sigma$ induit un morphisme de schémas $\Spec(W_1)\longrightarrow \Spec(W_1)$. Si $Y_1$ est un schéma sur $\Spec(W_1)$ on notera $Y'_1$ le schéma $\Spec(W_1) \times_{\Spec(W_1)} Y$, où le produit fibré est pris au-dessus de $\sigma$. Considérons $Q(T)=\sum b_kT^k$ un polynôme à coefficients dans $W_1$, on note $Q^{\sigma}(T)=\sum \sigma(b_k)T^k$ où $\sigma$ est un relèvement du Frobenius sur $W_1$. \\
\end{defi}

Par la propriété des vecteurs de Witt, le Frobenius sur $k$ (qui est l'élévation à la puissance $p$) se relève en un isomorphisme $\sigma : W\rightarrow W$. On a le diagramme suivant de schémas sur $\Spec(k)$ :
\begin{center}
\begin{tikzcd}
X_0 \arrow[bend left]{rr}{F_{abs}}
\arrow{r}{F} & X'_0 \arrow{r}{pr} \arrow{d} & X_0 \arrow{d} \\
& \Spec(k) \arrow{r}{\overline{\sigma}} & \Spec(k)
\end{tikzcd}
\end{center}
où $F_{abs}$ est le Frobenius absolu, $F$ le Frobenius relatif factorisant $F_{abs}$, $pr$ la projection de $X'_0$ sur $X_0$ et $\overline{\sigma}:\Spec(k)\rightarrow \Spec(k)$ est le morphisme induit par le Frobenius sur $k$. Soit $\mathcal{U}'\simeq \mathcal{U}$ (comme $\mathbb{Z}$-schéma) l'analogue de $\mathcal{U}$ sur $X'_0$, on a $\mathcal{U}=\Spec(k[t,y]/(y^n-f(t)))$ et donc $\mathcal{U}'=\Spec(k[t,y]/(y^n-f^{\sigma}(t)))$ (voir définition \ref{polysigma}). En restreignant le diagramme à l'ouvert $\mathcal{U}$ on obtient
\begin{center}
\begin{tikzcd}
k[t,y]/(y^n-f(t)) \arrow[bend left]{rr}{F^{\#}_{abs}} \arrow{r}{pr^{\#}} & k[t,y]/(y^n-f^{\sigma}(t))  \arrow{r}{F^{\#}}  & k[t,y]/(y^n-f(t))  
\end{tikzcd}
\end{center}
où $F^{\#}_{abs}(h)=h^p$, $pr^{\#}(h)=h^{\sigma}$, $F^{\#}(y)=y^p$ et $F^{\#}(t)=t^p$. Par ailleurs $F^{\#}_{abs}$ est semi-linéaire par rapport à $\sigma$, $pr^{\#}$ est un isomorphisme d'anneaux également semi-linéaire par rapport à $\sigma$ et $F^{\#}$ est un morphisme de $k$-algèbres. Par fonctorialité on a alors
\begin{center}
\begin{tikzcd}
H^1_{cris}(X_0) \arrow[bend left]{rr}{F_{abs}} \arrow{r} & H^1_{cris}(X'_0)  \arrow{r}{F}  & H^1_{cris}(X_0) 
\end{tikzcd}
\end{center}
où $F_{abs}$ est un morphisme semi-linéaire et $F$ un morphisme $K$-linéaire qui factorise $F_{abs}$. De plus, Berthelot ayant construit un isomorphisme $H^1_{cris}(X_0)\simeq H^1_{DR}(X)$ (et $H^1_{cris}(X'_0)\simeq H^1_{DR}(X')$, voir \cite{Berth}), on obtient par identification
\begin{center}
\begin{tikzcd}
H^1_{DR}(X) \arrow[bend left]{rr}{F_{abs}} \arrow{r}{pr} & H^1_{DR}(X')  \arrow{r}{F}  & H^1_{DR}(X) 
\end{tikzcd}
\end{center}
avec $F_{abs}$ semi-linéaire, $pr$ isomorphisme de groupes semi-linéaire par rapport à $\sigma$ et $F$ qui factorise $F_{abs}$. \\

Supposons que l'on ait un scindage $sc:H^1(X',\mathcal{O}_{X'})\oplus H^0(X',\Omega^1_{X'})\rightarrow H^1_{DR}(X)$ ainsi que l'application $F$ décrite précédemment. On peut alors construire le Frobenius divisé modulo $p$, 
\[\varphi':H^1(X'_0,\mathcal{O}_{X'_0})\oplus H^0(X'_0,\Omega^1_{X'_0})\rightarrow H^1_{DR}(X_0)\] 
qui est $k$-linéaire. On sait par un théorème de Huyghe-Wach que $\varphi'$ est donné par le morphisme de Deligne-Illusie. Il s'agit en fait de calculer $H^1(X_0,\mathcal{O}_{X_0})\oplus H^0(X_0,\Omega^1_{X_0})\rightarrow H^1_{DR}(X_0)$, le Frobenius divisé de Mazur modulo $p$, qui décrit le $\varphi$-module filtré modulo $p$ associé à $X_0$. Cette application est la composée 
\begin{center}
\begin{tikzcd}
H^1(X_0,\mathcal{O}_{X_0})\oplus H^0(X_0,\Omega^1_{X_0})  \arrow[bend left]{rr}{\phi} \arrow{r}{pr^{\#}} & H^1(X'_0,\mathcal{O}_{X'_0})\oplus H^0(X'_0,\Omega^1_{X'_0}) \arrow{r}{\phi'} & H^1_{DR}(X_0)
\end{tikzcd}
\end{center}
avec $\phi$ semi-linéaire. \\

\begin{rem} \label{pourlecalculdelamatrice}
On dispose d'une base de $H^1(X_0,\mathcal{O}_{X_0})\oplus H^0(X_0,\Omega^1_{X_0})\simeq H^1_{DR}(X_0)$ donnée par les propositions \ref{baseH1Ok} et \ref{baseH0k}. Alors $pr^{\#}(\omega_{i,j})=\omega_{i,j}'$ et $pr^{\#}(\overline{h_{i,j}})=\overline{h_{i,j}'}$ et donc $\phi(\omega_{i,j})=\phi'(\omega_{i,j}')$ et $\phi(\overline{h_{i,j}})=\phi'(\overline{h_{i,j}'})$. Ainsi la matrice de $\phi$ est donnée par les $\phi'(\omega_{i,j}')$ et $\phi'(\overline{h_{i,j}'})$ exprimés dans les bases $(\omega_{i,j})$ et $(\overline{h_{i,j}})$. Par ailleurs la base (($\omega_{i,j}$),($\overline{h_{i,j}}$)) est adaptée à la filtration et on peut donc, de façon abusive, considérer cette base au départ et à l'arrivée de $\phi$. \\
\end{rem}

\subsection{Relèvements du Frobenius modulo $p^2$}\label{paszero}

Afin de relever le Frobenius, on va modifier les ouverts $\mathcal{U}$ et $\mathcal{V}$ qui ont été définis dans la notation \ref{defisurfibrespe}. Posons $\overline{\mathcal{U}}=\mathcal{U}\cap D(y)$ et notons $\overline{\mathcal{F}}$ le recouvrement associé aux deux ouverts $\overline{\mathcal{U}}$ et $\mathcal{V}$.

\begin{lem}
Le recouvrement $\overline{\mathcal{F}}$ est un recouvrement de ${X_0}$. 
\end{lem}

\begin{pf}
Considérons $\overline{k}$ une clôture algébrique de $k$ et montrons que ${X_0}(\overline{k})=\overline{\mathcal{U}}(\overline{k})\cup \mathcal{V}(\overline{k})$. Par définition 
\[{X_0}(\overline{k})=\left\lbrace (t,y)\in \overline{k}^2 \ tq \ y^n=f(t) \right\rbrace \cup \left\lbrace (s,z)\in\overline{k}^2 \ tq \ z^n=f_2(s) \right\rbrace\]
où $ts=1$, $z=s^ry$ et $f_2(s)=f(t)s^{rn}$. Par ailleurs 
\[\mathcal{U}(\overline{k})=\overline{\mathcal{U}}(\overline{k})\cup \left\lbrace (t,0)\in\overline{k}^2 \ tq \ 0=f(t)\right\rbrace\]
or par hypothèse, si $f(t)=0$ alors $t\neq 0$ donc 
\[\left\lbrace (t,0)\in\overline{k}^2 \ tq \ 0=f(t)\right\rbrace\subset \mathcal{V}(\overline{k})\] 
et on a bien l'égalité souhaitée. \\
\end{pf}

Dans la suite on calcule un relèvement du Frobenius sur chacun des ouverts de ce nouveau recouvrement.\\

\begin{notat*}
Soient \textbf{$f$} et \textbf{$f_2$} des polynômes de $W[t]$ et $W[s]$ respectivement, définis comme précédemment. On se permet de noter encore \textbf{$f$} et \textbf{$f_2$}, éléments de $W_1[t]$ et $W_1[s]$, les classes modulo $p^2$ de ces polynômes. On notera \textbf{$f'$} et \textbf{$f'_2$} leurs dérivées dans $W_1[t]$ et $W_1[s]$ respectivement. On définit les ouverts affines
\begin{center}
$\mathcal{U}_1=\Spec(W_1[t,y]/(y^n-$\textbf{$f$}$(t)))$ \\
$\mathcal{V}_1=\Spec(W_1[s,z]/(z^n-$\textbf{$f_2$}$(s)))$
\end{center}
ainsi que $X_1$ comme le recollement des deux schémas $\mathcal{U}_1$ et $\mathcal{V}_1$ qui est une courbe sur $\Spec(W_1)$. Enfin on note $\overline{\mathcal{U}}_1=\mathcal{U}_1\cap D(y)$. \\
\end{notat*}

Le Frobenius de ${X_0}$ admet un relèvement modulo $p^2$ sur chacun des deux ouverts $\overline{\mathcal{U}}_1$ et $\mathcal{V}_1$ par lissité de ces deux ouverts. \\

\begin{rem}
Puisque $y$ et \textbf{$f$}$(t)$ sont liés par $y^n=$\textbf{$f$}$(t)$ il est clair que $\overline{\mathcal{U}}_1=\mathcal{U}_1\cap D($\textbf{$f$}$(t))$. Par ailleurs puisque $t$ et \textbf{$f$}$(t)$ sont premiers entre eux dans $W_1[t]$, l'inverse de \textbf{$f$}$(t)$ est en fait un élément de $W_1\llbracket t \rrbracket$. En effet \textbf{$f$}$(t)=\lambda_0+tw(t)=\lambda_0(1+a_0^{-1}tw(t))$ avec $\lambda_0\neq 0$ et \textbf{$f$}$(t)^{-1}=\lambda_0^{-1}\sum (-1)^k\lambda_0^{-k}t^kw^k(t)\in W\llbracket t \rrbracket$ et on a donc 
\begin{center}
$\mathcal{O}_{\mathcal{U}_1}(\mathcal{U}_1)[$\textbf{$f$}$(t)^{-1}]=W_1[t,y,$\textbf{$f$}$(t)^{-1}]/(y^n-$\textbf{$f$}$(t))\subset \displaystyle \bigoplus_{j=1}^{n-1} W_1 \llbracket t\rrbracket y^j$ 
\end{center}
comme $W_1[t]$-module. \\
\end{rem}

\begin{rem}
La multiplication par $p$ induit une bijection entre $k[t]$ et $pW_1[t]$. Dans toute la suite lorsque $Q$ sera un polynôme de $k[t]$ on fera l'abus d'écrire $pQ(t)$ pour désigner l'image de $Q(t)$ sous cette bijection. \\
\end{rem}

\begin{prop} \label{frob1}
Il existe un relèvement du Frobenius $\overline{\mathcal{U}}_1\rightarrow \overline{\mathcal{U}'}_1$, et un polynôme $\mathcal{P}(t)$ de $k[t]$, tel que le relèvement modulo $p^2$ du Frobenius 
\begin{center}
$F_U:W_1[t,y,y^{-1}]/(y^n-$\textbf{$f$}$^{\sigma}(t))\longrightarrow W_1[t,y,y^{-1}]/(y^n-$\textbf{$f$}$(t))$ 
\end{center}
soit déterminé par 
\begin{align*}
F_U(t)=&\ t^p\\
F_U(y)=&\ y^p+py^p\left(\dfrac{\mathcal{P}(t)}{n}y^{-np}\right) \\
\end{align*}
où $\mathcal{P}(t)= \left( \right. $\textbf{$f$}$^{\sigma}(t^p)-$\textbf{$f$}$(t)^p\left. \right)/p$. 
\end{prop}

\begin{pf}
Ce résultat est inspiré du relèvement du Frobenius donné par Kedlaya dans \cite{Memoire} pour les courbes hyperelliptiques et par Gaudry-Gürel dans \cite{MR1934859} dans le cas superelliptique. On rappelle que $\mathcal{O}_{X_1}(\overline{\mathcal{U}}_1)=W_1[t,$\textbf{$f$}$(t)^{-1},y]/(y^n-$\textbf{$f$}$(t))$. On cherche un relèvement de la forme $t\mapsto F_U(t)=t^p$ et $y\mapsto F_U(y)$ tel que
\begin{center}
$(F_U(y))^n=$\textbf{$f$}$^{\sigma}(F_U(t))$.
\end{center}
Comme \textbf{$f$}$^{\sigma}(t^p)-$\textbf{$f$}$(t)^p$ est dans $p(W_1[t])$, l'élément $\mathcal{P}(t)=($\textbf{$f$}$^{\sigma}(t^p)-$\textbf{$f$}$(t)^p)/p$ est bien défini et est un polynôme de $k[t]$. On a clairement \textbf{$f$}$^{\sigma}(F_U(t))=$\textbf{$f$}$^{\sigma}(t^p)=$\textbf{$f$}$(t)^p+p\mathcal{P}(t)$. Comme $y$ est inversible sur $\overline{\mathcal{U}}_1$, on peut poser 
\[F_U(y)=y^p\left( 1+p\dfrac{\mathcal{P}(t)}{n}y^{-np}\right).\]
Alors on a
\[(F_U(y))^n=\left( y^p\left( 1+p\dfrac{\mathcal{P}(t)}{n}y^{-np}\right)\right) ^n=y^{np}\sum_{k=0}^n \binom{n}{k}\left( \dfrac{1}{n}p\mathcal{P}(t)y^{-np}\right)^k.\]
Réduisons maintenant cette expression modulo $W_1$, alors
\begin{center}
$(F_U(y))^n\equiv $\textbf{$f$}$(t)^p\left( 1+n\left( \dfrac{1}{n}p\mathcal{P}(t)y^{-np}\right) \right) \equiv $\textbf{$f$}$(t)^p(1+p\mathcal{P}(t)$\textbf{$f$}$(t)^{-p})\equiv $\textbf{$f$}$(t)^p+p\mathcal{P}(t).$
\end{center}
Et donc $(F_U(y))^n=$\textbf{$f$}$^{\sigma}(F_U(t))$. Reste à vérifier que $F_U(y)$ est inversible dans $\mathcal{O}_{X_1}(\overline{\mathcal{U}}_1)$. Or dans cette $W_1$-algèbre on a 
\[\left( 1+p\dfrac{\mathcal{P}(t)}{n}y^{-np}\right)\left( 1-p\dfrac{\mathcal{P}(t)}{n}y^{-np}\right)=1-\left( p\dfrac{\mathcal{P}(t)}{n}y^{-np}\right)^2=1\]
donc 
\[F_U(y^{-1})=F_U(y)^{-1}=y^{-p}\left( 1-p\dfrac{\mathcal{P}(t)}{n}y^{-np}\right)\in \mathcal{O}_{X_1}(\overline{\mathcal{U}}_1)\]
ce qui établit que $F_U$ défini par $F_U(t)=t^p$ et $F_U(y)=y^p+py^p\left(\mathcal{P}(t)/{n}y^{-np}\right)$ est un endomorphisme d'algèbres qui passe au quotient et montre l'énoncé. \\
\end{pf}

Décrivons maintenant l'action de $F_U$ sur les composantes isotypiques de $\mathcal{O}_{X_1}(\overline{\mathcal{U}_1})$ associées à l'action $\delta$, définie à la section \ref{action}. 

\begin{prop} \label{isotyfrob1}
On conserve les notations précédentes et on note $b$ le reste de la division euclidienne de $pj$ par $n$ pour $1\leq j\leq n-1$ fixé. Supposons que $\mu_n(\overline{\mathbb{F}_p})$ soit inclus dans $k$, alors le morphisme $F_U:\mathcal{O}_{X_1}(\overline{\mathcal{U}_1}) \rightarrow \mathcal{O}_{X_1}(\overline{\mathcal{U}_1})$ envoie la $j$-ième composante isotypique de $\mathcal{O}_{X_1}(\overline{\mathcal{U}_1})$ sur la $b$-ième composante isotypique de $\mathcal{O}_{X_1}(\overline{\mathcal{U}_1})$. De plus l'application $j\mapsto b$ est une bijection de $\left\lbrace 1,...,n-1\right\rbrace$. 
\end{prop}

\begin{pf}
Soit $1\leq j\leq n-1$ fixé. On a $F_U(y^j)=y^{pj}(1+p\mathcal{P}(t)y^{-npj}/n)$ dans $\mathcal{O}_{X_1}(\overline{\mathcal{U}_1})$ ainsi, si on effectue la division euclidienne $pj=an+b$ on a $F_U(y^j)=$\textbf{$f$}$(t)^ay^b(1+p\mathcal{P}(t)y^{-np}/n)$. Un élément de la $(n-j)$-ième composante isotypique est donc envoyé sur un élément de la $(n-b)$-ième composante isotypique puisque $t$ et $y^n$ sont invariants par cette action. Par ailleurs comme $p$ est premier à $n$ et que $j$ est strictement inférieur à $n$, l'application qui à $j$ associe $b$ est bien une bijection. \\
\end{pf}

Nous nous intéressons désormais à l'ouvert $\mathcal{V}_1$ de notre recouvrement.\\

\begin{prop} \label{frob2}
Il existe un relèvement du Frobenius $\mathcal{V}_1\rightarrow \mathcal{V}'_1$ et deux polynômes $v(s),b(s)$ de $k[s]$ tels que le relèvement modulo $p^2$ du Frobenius 
\begin{center}
$F_V:W_1[s,z]/(z^n-$\textbf{$f_2$}$^{\sigma}(s))\longrightarrow W_1[s,z]/(z^n-$\textbf{$f_2$}$(s))$ 
\end{center}
soit déterminé par 
\begin{center}
$F_V(s)=s^p+pv(s)$\\
$F_V(z)=z^p+pz^pb(s)$\\
\end{center}
où $v$ et $b$ sont reliés par la relation :
\begin{center}
\textbf{$f_2$}$(s)^p-$\textbf{$f_2$}$^{\sigma}(s^p)=p\left( \right. v(s)$\textbf{$f'_2$}$(s)^p-nb(s)$\textbf{$f_2$}$(s)^p\left. \right)$. \\
\end{center}
\end{prop}

\begin{pf}
Par définition on a \textbf{$f_2$}$^{\sigma}(s^p)\equiv $\textbf{$f_2$}$(s)^p$ modulo $p(W_1[s])$(voir définition \ref{polysigma}) donc il existe un polynôme $Q(s)\in k[s]$ tel que \textbf{$f_2$}$(s)^p-$\textbf{$f_2$}$^{\sigma}(s^p)=pQ(s)$ dans $W_1[s]$. De plus $f_2$ et $f'_2$ sont premiers entre eux dans $k[s]$ donc il existe deux polynômes, $v_0(s),b_0(s)\in k[s]$ tels que  
\[v_0(s)f'_2(s)^p+b_0(s)f_2(s)^p=1\] 
et on peut donc trouver deux polynômes $v(s),b(s)\in k[s]$ tels que 
\[Q(s)= v(s)f'_2(s)^p-nb(s)f_2(s)^p\] 
(qui est bien la relation annoncée sur $v$ et $b$) en posant $v(s)=Q(s)v_0(s)$ et $nb(s)=-Q(s)b_0(s)$. Posons $F_V(s)=s^p+pv(s)$ et $F_V(z)=z^p+pz^pb(s)$ et montrons que cela définit bien un relèvement du Frobenius sur $\mathcal{V}_1$ i.e. que $F_V(z)^n=$\textbf{$f_2$}$^{\sigma}(F_V(s))$ dans $\mathcal{O}_{X_1}(\mathcal{V}_1)$. Avec cette définition 
\begin{center}
$F_V(z)^n=z^{pn}(1+pb(s))^n=$\textbf{$f_2$}$(s)^p(1+pb(s))^n$. 
\end{center}
On utilise alors la formule du binôme de Newton pour développer le second terme et comme on travaille modulo $W_1$ on a finalement 
\begin{center}
$F_V(z)^n=$\textbf{$f_2$}$(s)^p(1+npb(s))=$\textbf{$f_2$}$(s)^p+npb(s)$\textbf{$f_2$}$(s)^p=$\textbf{$f_2$}$(s)^p+npb(s)f_2(s)^p$ 
\end{center}
car dans $W_1$ on a $p$\textbf{$f_2$}$(s)=f_2(s)$. En utilisant la relation établie précédemment entre $v$ et $b$ on a donc $F_V(z)^n=$\textbf{$f_2$}$^{\sigma}(s^p)+pv(s)f'_2(s)^p$ dans $\mathcal{O}_{X_1}(\mathcal{V}_1)$. Calculons alors \textbf{$f_2$}$^{\sigma}(F_V(s))=$\textbf{$f_2$}$^{\sigma}(s^p+pv(s))$. En utilisant la formule du binôme de Newton on montre que $(s^p+pv(s))^m=(s^p)^m+pv(s)m(s^p)^{m-1}$ dans $W_1[s]$ pour tout entier $m$ et donc dans $W_1$ on a 
\begin{center}
\textbf{$f_2$}$^{\sigma}(s^p+pv(s))=$\textbf{$f_2$}$^{\sigma}(s^p)+pv(s)$\textbf{$f'$}$^{\sigma}_2(s^p)=$\textbf{$f_2$}$^{\sigma}(s^p)+pv(s)f'^{\sigma}_2(s^p)$. 
\end{center}
Comme par ailleurs $f'^{\sigma}_2(s^p)= f'_2(s)^p$ dans $k[s]$ on a finalement 
\begin{center}
\textbf{$f_2$}$^{\sigma}(F_V(s))=$\textbf{$f_2$}$^{\sigma}(s^p)+pv(s)f'_2(s)^p=F_V(z)^n\in \mathcal{O}_{X_1}(\mathcal{V}_1)$
\end{center}
ce qui achève la démonstration. \\
\end{pf}

Décrivons maintenant l'action de $F_V$ sur les composantes isotypiques de $\mathcal{O}_{X_1}(\mathcal{V}_1)$ associées à l'action $\delta$, définie à la section \ref{action}. 

\begin{prop} \label{isotyfrob2}
On conserve les notations précédentes et on note $b$ le reste de la division euclidienne de $pj$ par $n$ pour $1\leq j\leq n-1$ fixé. Supposons que $\mu_n(\overline{\mathbb{F}_p})$ soit inclus dans $k$ alors le morphisme $F_V:\mathcal{O}_{X_1}(\mathcal{V}_1) \rightarrow \mathcal{O}_{X_1}(\mathcal{V}_1)$ envoie la $j$-ième composante isotypique de $\mathcal{O}_{X_1}(\mathcal{V}_1)$ sur la $b$-ième composante isotypique de $\mathcal{O}_{X_1}(\mathcal{V}_1)$. De plus l'application $j\mapsto b$ est une bijection de $\left\lbrace 1,...,n-1\right\rbrace$. 
\end{prop}

\begin{pf}
Soit $1\leq j\leq n-1$ fixé. On a $F_V(z^j)=z^{pj}(1+pb(s))$ et donc si on effectue la division euclidienne $pj=an+b$ on a $F_V(z^j)=$\textbf{$f_2$}$(s)^az^b(1+pb(s))$. Comme $s$ est invariant par $\delta$ on a bien le résultat annoncé. Par ailleurs comme $p$ est premier à $n$ et que $j$ est strictement inférieur à $n$, l'application qui à $j$ associe $b$ est bien une bijection. \\
\end{pf}

Nous aurons besoin ultérieurement (voir remarque \ref{DI1explicite}) de la

\begin{prop} \label{divisible}
Considérons $v(s)$ tel que défini dans la proposition \ref{frob2}, alors
\[f_2(s)^{p-1}\vert (v'(s)+s^{p-1})\]
dans $k[s]$. 
\end{prop}

\begin{pf}
On a $pQ(s)=$\textbf{$f_2$}$(s)^p-$\textbf{$f_2$}$^{\sigma}(s^p)$ donc $pQ'(s)=p$\textbf{$f'_2$}$(s)$\textbf{$f_2$}$(s)^{p-1}-ps^{p-1}$\textbf{$f'$}$^{\sigma}_2(s^p)$ dans $W_1[s]$. Soit après simplification par $p$, $Q'(s)=f'_2(s)f_2(s)^{p-1}-s^{p-1}f'^{\sigma}_2(s^p)$ dans $k[s]$. D'autre part dans $k[s]$, $Q'(s)=v'(s)f'_2(s)^p-nb'(s)f_2(s)^p$ et $f'^{\sigma}_2(s^p)=f'_2(s^p)=f'_2(s)^p$, on obtient donc
$$f_2(s)^{p-1}(f'_2(s)+nb'(s)f_2(s))=f'_2(s)^p(v'(s)+s^{p-1}) \in k[s]$$
mais $f_2$ et $f_2'$ sont premiers entre eux dans $k[s]$ et donc $f_2(s)^{p-1}\vert (v'(s)+s^{p-1})$ dans $k[s]$. \\
\end{pf}

D'après $3.8(c)$ de \cite{PT} on peut introduire sur $\overline{\mathcal{U}}$ (respectivement $\mathcal{V}$) l'application 
\begin{center}
$f_U:\Omega^1_{\overline{\mathcal{U}'}}\rightarrow F_*\Omega^1_{\overline{\mathcal{U}}}$ (respectivement $f_V:\Omega^1_{\mathcal{V}'}\rightarrow F_*\Omega^1_{\mathcal{V}}$) 
\end{center}
définie par 
\begin{center}
$f_U(dt)=\dfrac{F_U(dt)}{p}$ (respectivement $f_V(ds)=\dfrac{F_V(ds)}{p}$)
\end{center}
qu'on appelle par la suite Frobenius divisé différentiel. \\

\subsection{Calcul du morphisme de Deligne-Illusie}\label{Complexe}

\begin{notat*}
Dans la suite si $\mathcal{E}$ est un faisceau sur ${X_0}$, si $x$ est une section locale de $\mathcal{E}$ sur un ouvert $U$ de ${X_0}$ on notera $x\in \mathcal{E}$. \\
\end{notat*}

On rappelle que $\overline{\mathcal{F}}$ est le recouvrement de la courbe ${X_0}$ par les deux ouverts $\overline{\mathcal{U}}$ et $\mathcal{V}$. Alors le complexe de \v{C}ech associé est :
\[\begin{array}{ccccccc}
    & & 0 & & 0 & & \\
    & & \uparrow & & \uparrow & & \\
    0 & \rightarrow & \Omega_{{X_0}}^1(\overline{\mathcal{U}}) \times \Omega_{{X_0}}^1(\mathcal{V}) & \rightarrow & \Omega_{{X_0}}^1(\overline{\mathcal{U}}\cap \mathcal{V}) & \rightarrow & 0 \\
    & & \uparrow & & \uparrow & & \\
    0 & \rightarrow & \mathcal{O}_{{X_0}}(\overline{\mathcal{U}}) \times  \mathcal{O}_{{X_0}}(\mathcal{V}) & \rightarrow & \mathcal{O}_{{X_0}}(\overline{\mathcal{U}}\cap \mathcal{V}) & \rightarrow & 0 \\
    & & \uparrow & & \uparrow & & \\
    & & 0 & & 0 & & \\
\end{array}\]
auquel on peut associer le complexe simple suivant : 
\[0 \rightarrow \mathcal{O}_{{X_0}}(\overline{\mathcal{U}}) \times  \mathcal{O}_{{X_0}}(\mathcal{V}) \rightarrow \Omega_{{X_0}}^1(\overline{\mathcal{U}}) \times \Omega_{{X_0}}^1(\mathcal{V})\oplus \mathcal{O}_{{X_0}}(\overline{\mathcal{U}}\cap \mathcal{V}) \rightarrow \Omega_{{X_0}}^1(\overline{\mathcal{U}}\cap \mathcal{V}) \rightarrow 0.\]
De la même manière on a un complexe double et un complexe simple associés au recouvrement $\overline{\mathcal{F}}_1$ de $X_1$  par $\overline{\mathcal{U}}_1$ et $\mathcal{V}_1$. On se place ici sur la courbe ${X_0}$. Ainsi on considère les ouverts $\overline{\mathcal{U}}$ et $\mathcal{V}$ et on travaille donc modulo $p$. \\ 

Le choix d'un relèvement de Frobenius sur les deux ouverts affines $\overline{\mathcal{U}}_1$ et $\mathcal{V}_1$ permet de construire un morphisme $\mathcal{O}_{X'_0}$-linéaire
\begin{center}
$\displaystyle DI : \bigoplus _{i=0}^1 \Omega_{X'_0}^i [-i] \rightarrow F_*$\v{C}$(\overline{\mathcal{F}},\Omega_{{X_0}}^{\bullet})$
\end{center}
en suivant la méthode donnée par Illusie dans \cite{PT} au paragraphe 5.1. On notera $\Phi$ le morphisme induit 
\[\Phi : H^1(X'_0,\mathcal{O}_{X'_0})\oplus H^0(X'_0,\Omega^1_{X'_0})\rightarrow H^1_{DR}({X_0}).\] 

\begin{rem}\label{conditiondiff}
Dans notre cas on dispose d'un recouvrement à deux ouverts de $X_1$ par $\overline{\mathcal{U}}_1$ et $\mathcal{V}_1$ sur lesquels nous avons déjà calculé un relèvement du Frobenius (voir propositions \ref{frob1} et \ref{frob2}). Il n'y a alors qu'à calculer $f_U$ et $f_V$ le Frobenius divisé différentiel sur $\overline{\mathcal{U}}$ respectivement $\mathcal{V}$ et un homomorphisme $h:\Omega^1_{X'_0}\rightarrow F_*\mathcal{O}_{\overline{\mathcal{U}}\cap \mathcal{V}}$ devant vérifier la condition
\[ f_V(\omega\vert_{\mathcal{V}})-f_U(\omega\vert_{\mathcal{U}})=dh_{ij}(\omega)\vert_{\mathcal{U}\cap \mathcal{V}}\]
sur $\overline{\mathcal{U}}\cap \mathcal{V}$ pour tout élément $\omega$ de $\Omega^1_{X'_0}$. \\
\end{rem}

Rappelons la définition de $DI$ et donnons des expressions explicites de $f_U$, $f_V$ et $h$ dans notre cadre de travail. \\

\begin{prop} \label{morphiDI}
En degré 0, $DI_0:\mathcal{O}_{X'_0}\rightarrow F_*\mathcal{O}_{\overline{\mathcal{U}}}\oplus F_*\mathcal{O}_{\mathcal{V}}$ est donné par le Frobenius relatif $F$ i.e. $DI_0(\alpha)=(F(\alpha)\vert_{\overline{\mathcal{U}}},F(\alpha)\vert_{\mathcal{V}})$. \\
En degré 1, 
$$\begin{array}{cccc}
    DI_1 : & \Omega^1_{X'_0} & \longrightarrow & F_*\Omega_{\overline{\mathcal{U}}}^1 \oplus  F_*\Omega_{\mathcal{V}}^1 \oplus F_*\mathcal{O}_{\overline{\mathcal{U}}\cap \mathcal{V}} \\
    & \omega = \alpha_0dt=\alpha_1ds & \longmapsto &  (\alpha_0^pf_U(dt) , \alpha_1^pf_V(ds) , \alpha_0^ph(dt)) \\
\end{array}$$
où $f_U$ est le Frobenius divisé différentiel sur $\overline{\mathcal{U}}$, $f_V$ est le Frobenius divisé différentiel sur $\mathcal{V}$ correspondants aux relèvement du Frobenius donnés par les propositions \ref{frob1} et \ref{frob2} et $h$ vérifie la condition donnée à la remarque \ref{conditiondiff}. On remarque que 
\begin{align*}   f_U(dt) =&  \ t^{p-1}dt \\
    f_V(ds) =& \  (s^{p-1}+v'(s))ds \\
    h(dt) =& \  t^{2p}v(s)   
\end{align*} 
où $v$ est le polynôme du relèvement du Frobenius sur $\mathcal{V}_1$ (voir proposition \ref{frob2}). 
\end{prop}

\begin{pf}
En degré $0$ il n'y a rien à montrer. En degré 1,  $f_U$ (resp. $f_V$) est simplement le Frobenius divisé associé au relèvement du Frobenius sur $\overline{\mathcal{U}}_1$ (resp. $\mathcal{V}_1$). On rappelle que sur $\overline{\mathcal{U}}_1$ le relèvement du Frobenius (modulo $p^2$) choisi vérifie $F_U(t)=t^p$. Alors $F_U(dt)=dF_U(t)=d(t^p)=pt^{p-1}dt$. Or par définition 
\[f_U(dt)=\dfrac{F_U(dt)}{p}=\dfrac{pt^{p-1}dt}{p}=t^{p-1}dt\in \Omega^1_{{X_0}}.\]
Sur $\mathcal{V}_1$, $F_V(s)=s^p+pv(s)$. Alors $F_V(ds)=dF_V(s)=d(s^p+pv(s))=p(s^{p-1}+v'(s))ds$. D'où 
\[f_V(ds)=\dfrac{F_V(ds)}{p}=\dfrac{p(s^{p-1}+v'(s))ds}{p}=(s^{p-1}+v'(s))ds\in \Omega^1_{{X_0}}.\]
Par définition $h$ est une application définie sur $\overline{\mathcal{U}}\cap \mathcal{V}$ vérifiant $f_U(dt)-f_V(dt)=dh(dt)$. On rappelle que $t$ est un paramètre sur $\overline{\mathcal{U}}\cap \mathcal{V}$. Calculons $f_V(dt)$. Comme $t=s^{-1}$ alors $dt=-s^{-2}ds$, d'où 
\[f_V(dt)=f_V(-s^{-2}ds)=(-s^{-2})^pf_V(ds)=-s^{-2p}(s^{p-1}+v'(s))ds.\]
On a donc $dh(dt)=f_U(dt)-f_V(dt)=t^{p-1}dt-(-s^{-2p}(s^{p-1}+v'(s))ds)$. Utilisons à nouveau les relations entre $s$ et $t$ pour simplifier cette expression : 
\[dh(dt)=t^{p-1}dt-(-t^{2p}(t^{1-p}+v'(1/t)))(-t^{-2}dt)=-t^{2p-2}v'(1/t)dt=t^{2p}(-1/t^2)v'(1/t)dt=t^{2p}(v(1/t))'.\]
Enfin ${X_0}$ étant un $k$-schéma, on a $(t^{2p})'=2pt^{2p-1}dt=0$ d'où $(t^{2p}v(s))'=t^{2p}(v(1/t))'$ sur $\overline{\mathcal{U}}\cap \mathcal{V}$ et donc $h(dt)=t^{2p}v(s)$. Par définition il est clair que sur $\overline{\mathcal{U}}\cap \mathcal{V}$ la relation $f_V-f_U+dh=0$ est vérifiée. Enfin la construction du morphisme de Deligne et Illusie garantit que $f_U(\omega)\in F_*\Omega_{{X_0}}^1(\overline{\mathcal{U}})$, $f_V(\omega)\in F_*\Omega_{{X_0}}^1(\mathcal{V})$ et $h(\omega)\in F_*\mathcal{O}_{{X_0}}(\overline{\mathcal{U}}\cap \mathcal{V})$ pour tout élément $\omega$ de $\Omega_{X'_0}^1$. \\
\end{pf}

\begin{rem} \label{DI1explicite}
Une base de $H^0(X'_0,\Omega_{X'_0}^1)$ est donnée par $\omega'_{i,j}=t^iy^{-j}dt$ avec $1\leq j\leq n-1$ et $0\leq i\leq rj-2$ (voir proposition \ref{basesprimeesH0}). On a donc 
\begin{center}
$f_U(\omega'_{i,j})=\dfrac{t^{p(i+1)-1}}{y^{pj}}dt\in \Omega^1_{{X_0}}(\overline{\mathcal{U}})$ \\
$f_V(\omega'_{i,j})=-s^{p(rj-(i+2))}\dfrac{(s^{p-1}+v'(s))}{z^{pj}}ds\in \Omega^1_{{X_0}}(\mathcal{V})$ \\
$h(\omega'_{i,j})=\dfrac{t^{p(i+2)}}{y^{pj}}v(s)\in \mathcal{O}_{{X_0}}(\overline{\mathcal{V}}\cap \mathcal{V}).$ \\
\end{center}
Soit encore, en notant $a$ le quotient de la division euclidienne de $pj$ par $n$ et $b$ son reste (non nul puisque $n$ et $p$ sont premiers entre eux et $1\leq j\leq n-1$)
\begin{center}
$f_U(\omega'_{i,j})=\dfrac{t^{p(i+1)-1}}{f(t)^a}\dfrac{dt}{y^b}$ \\
$f_V(\omega'_{i,j})=-s^{p(rj-(i+2))}\dfrac{(s^{p-1}+v'(s))}{f_2(s)^a}\dfrac{ds}{z^b}$ \\
$h(\omega'_{i,j})=\dfrac{t^{p(i+2)}}{f(t)^{a+1}}v(s)y^{n-b}.$ \\
\end{center}
Enfin grâce à la propriété de divisibilité énoncée à la proposition \ref{divisible} on constate qu'on a bien 
\[f_V(\omega'_{i,j}) \in k[s]\dfrac{ds}{z^b}\subset \Omega_{\mathcal{V}}^1=\bigoplus_{i=0}^n k[s]\dfrac{ds}{y^j}.\]
En effet $a<p$ car si on avait $a\geq p$ alors $pj=an+b\geq pn+b\geq pn>p(n-1)\geq pj$ ce qui est absurde. Ainsi $f_2(s)^a\vert f_2(s)^{p-1}$ dans $k[s]$ et on conclut grâce à la proposition \ref{divisible}.  \\
\end{rem}

\begin{lem}\label{pourledi0}
Notons $\mathcal{F}'$ le recouvrement de $X'_0$ par les deux ouverts $\overline{\mathcal{U}}'$ et $\mathcal{V}'$ et \v{C}$(\mathcal{F}',\mathcal{O}_{X'_0})$ le complexe de \v{C}ech de $\mathcal{O}_{X'_0}$ relativement à $\mathcal{F}'$. \\
(i) Le morphisme $DI_0$ induit par fonctorialité un morphisme de complexes 
\begin{center}
$\alpha : $\v{C}$(\mathcal{F}',\mathcal{O}_{X'_0})\rightarrow $\v{C}$(\overline{\mathcal{F}},\Omega_{{X_0}}^{\bullet})$
\end{center}
donné par \\
\adjustbox{scale=1,center}{%
\begin{tikzcd}
0 \arrow{r} & \mathcal{O}_{X'_0}(\overline{\mathcal{U}}') \times  \mathcal{O}_{X'_0}(\mathcal{V}') \arrow{r} \arrow{d}{\alpha_0} & \mathcal{O}_{X'_0}(\overline{\mathcal{U}}'\cap \mathcal{V}') \arrow{r} \arrow{d}{\alpha_1} & 0 \arrow{d} \\
0 \arrow{r} & \mathcal{O}_{{X_0}}(\overline{\mathcal{U}}) \times  \mathcal{O}_{{X_0}}(\mathcal{V}) \arrow{r} & \Omega_{{X_0}}^1(\overline{\mathcal{U}}) \times \Omega_{{X_0}}^1(\mathcal{V})\oplus \mathcal{O}_{{X_0}}(\overline{\mathcal{U}}\cap \mathcal{V}) \arrow{r} & \Omega_{{X_0}}^1(\overline{\mathcal{U}}\cap \mathcal{V}) \arrow{r} & 0
\end{tikzcd}
}
avec $\alpha_0 (h_1,h_2)=(F^{\#}(h_1),F^{\#}(h_2))$ et $\alpha_1(h)=(0,0,F^{\#}(h))$. \\
(ii) Le morphisme induit par $DI_0$, $H^1(X'_0,\mathcal{O}_{X'_0})\rightarrow H^1_{DR}({X_0})$ est donc induit par le morphisme de complexes $\alpha$ précédent.
\end{lem}

\begin{pf}
Il est clair que $\alpha$ ainsi défini est bien un morphisme de complexes. \\
\end{pf}

Décrivons enfin comme le morphisme $DI$ agit sur les composantes isotypiques associées à l'action $\delta$ décrite à la section \ref{action}.\\

\begin{prop} \label{lapropositionquifaut}
On conserve les notations précédentes. Soit $1\leq j\leq n-1$ fixé, notons $b$ le reste de la division euclidienne de $pj$ par $n$. Enfin supposons que $\mu_n(\overline{\mathbb{F}_p})$ soit inclus dans $k$, alors le morphisme de Deligne-Illusie induit un isomorphisme 
\[(H^0({X_0},\Omega_{{X_0}}^1))_j\oplus (H^1({X_0},\mathcal{O}_{{X_0}}))_j\simeq (H^1_{DR}({X_0}))_{\lfloor pj/n \rfloor}\] 
où l'indexation par $j$ désigne le sous-espace correspondant à la $j$-ième composante isotypique. 
\end{prop}

\begin{pf}
Les arguments de la démonstration sont à nouveau sensiblement les mêmes que pour les propositions \ref{isotau}, \ref{isotyfrob1} et \ref{isotyfrob2}. Par définition du Frobenius relatif, il est évident qu'en degré $0$ un élément de la $j$-ième composante isotypique de $\mathcal{O}_{X'_0}$ est envoyé sur un élément de la $b$-ième composante isotypique de $F_*\mathcal{O}_{\overline{\mathcal{U}}}\oplus F_*\mathcal{O}_{\mathcal{V}}$ où $b$ est le reste de la division euclidienne de $pj$ par $n$. En degré 1, grâce aux formules explicites données dans la remarque précédente on peut dire également qu'un élément de la $j$-ième composante isotypique de $H^0({X_0},\Omega_{{X_0}}^1)$ est envoyé sur un élément de la $b$-ième composante isotypique de $F_*\Omega_{\overline{\mathcal{U}}}^1 \oplus  F_*\Omega_{\mathcal{V}}^1 \oplus F_*\mathcal{O}_{\overline{\mathcal{U}}\cap \mathcal{V}}$ avec toujours $b$ reste de la division euclidienne de $pj$ par $n$. Les mêmes arguments que précédemment permettent d'affirmer que $j\mapsto b$ est une permutation. Enfin puisque le morphisme de Deligne-Illusie est un isomorphisme, cela donne un isomorphisme au niveau de composantes isotypiques. \\
\end{pf}

\begin{rem} \label{laremarquequifaut}
Plus généralement même lorsque la condition $\mu_n(\overline{\mathbb{F}_p})\subset k$ n'est pas vérifiée, définissons $(H^0({X_0},\Omega_{{X_0}}^1))_j$ comme le sous-espace de $H^0({X_0},\Omega_{{X_0}}^1)$ engendré par $t^iy^{-j}dt$ pour $0\leq i\leq rj-2$ et $(H^1({X_0},\mathcal{O}_{{X_0}}))_j$ comme le sous-espace de $H^1({X_0},\mathcal{O}_{{X_0}})$ engendré par $t^{-i}y^{n-j}$ pour $1\leq i\leq r(n-j)-1$, c'est-à-dire ce que l'on pourrait appeler par abus de langage la $j$-ième composante isotypique de ces espaces respectifs. Alors comme précédemment on peut montrer que le morphisme de Deligne-Illusie induit un isomorphisme 
\[(H^0({X_0},\Omega_{{X_0}}^1))_j\oplus (H^1({X_0},\mathcal{O}_{{X_0}}))_j\simeq (H^1_{DR}({X_0}))_{\lfloor pj/n \rfloor}\] 
où $(H^1_{DR}({X_0}))_j$ désigne le sous-espace qui correspondrait à la $j$-ième composante isotypique. Cela nous permet de conclure qu'il y a une permutation des "composantes isotypiques" par le morphisme de Deligne-Illusie et donc que la matrice du Frobenius présentera des blocs de zéros et $n-1$ blocs inversibles de taille $rn-2=l-1$ (qui correspond bien à une matrice de taille $(n-1)(l-1)=2g$). Il est donc possible de choisir une présentation de la matrice à la fois adaptée à cette permutation et à la filtration. Pour des raisons de simplicité d'écriture on choisira cependant une présentation uniquement adaptée à la filtration (voir section \ref{presentationmatrice}). \\
\end{rem}

\subsection{Calcul théorique de la matrice du Frobenius divisé sur $H^1_{DR}({X_0})$}

On a la suite exacte de $k$-espaces vectoriels $0\rightarrow H^0({X_0},\Omega_{{X_0}}^1) \rightarrow H^1_{DR}({X_0}) \rightarrow H^1({X_0},\mathcal{O}_{{X_0}}) \rightarrow 0$. Par ailleurs nous avons déjà calculé un scindage de cette suite exacte dans la proposition \ref{scindagemodp}. Dans la suite, on calcule $\Phi$ le Frobenius divisé, qui a été défini à la section \ref{rappels}, sur $H^0(X'_0,\Omega_{X'_0}^1)\oplus H^1(X'_0,\mathcal{O}_{X'_0})$. On utilise le scindage $[\tau]$ construit à la section \ref{calculscindage} pour identifier $H^1_{DR}({X_0})$, qui est isomorphe à $H^1($\v{C}$(\mathcal{F},\Omega_{{X_0}}^{\bullet}))$, à $H^1({X_0},\mathcal{O}_{{X_0}})\oplus H^0({X_0},\Omega^1_{{X_0}})$. \\

\begin{defi} \label{lesfonctionsquivontbien}
On rappelle quelques définitions et notations qui serviront dans toute cette section. L'identification de $H^1({X_0},\mathcal{O}_{{X_0}})$ à $\mathcal{O}_{{X_0}}(\overline{\mathcal{U}}\cap \mathcal{V})/Im(d_{0,1})$ induit (voir remarque \ref{identH1Ok})
\begin{center}
$\iota:H^1({X_0},\mathcal{O}_{{X_0}})\rightarrow \mathcal{O}_{{X_0}}(\overline{\mathcal{U}}\cap \mathcal{V})$
\end{center}
ainsi que 
\begin{center}
$\overline{...\vphantom{t}}:\mathcal{O}_{{X_0}}(\overline{\mathcal{U}}\cap \mathcal{V})\rightarrow H^1({X_0},\mathcal{O}_{{X_0}})$.
\end{center}
Par ailleurs on a construit à la proposition \ref{scindagemodp}
\begin{center}
$\tau : H^1({X_0},\mathcal{O}_{{X_0}})\rightarrow \Omega_{{X_0}}^1(\overline{\mathcal{U}}) \times\Omega_{{X_0}}^1(\mathcal{V})\times \mathcal{O}_{{X_0}}(\overline{\mathcal{U}}\cap \mathcal{V})$
\end{center}
et on note
\begin{center}
$\left[ \vphantom{t}...\right] : \Omega_{{X_0}}^1(\overline{\mathcal{U}}) \times\Omega_{{X_0}}^1(\mathcal{V})\times \mathcal{O}_{{X_0}}(\overline{\mathcal{U}}\cap \mathcal{V})\rightarrow H^1($\v{C}$(\mathcal{F},\Omega_{{X_0}}^{\bullet}))$. \\
\end{center}
On notera enfin $P_1$ (resp. $P_2$) la projection de $H^0({X_0},\Omega^1_{{X_0}}) \oplus H^1({X_0},\mathcal{O}_{{X_0}})$ sur le second (resp. le premier) facteur. \\
\end{defi}

\begin{prop} \label{premcol}
Soit $\omega_{i,j}=t^iy^{-j}dt$ un élément de la base de $H^0({X_0},\Omega^1_{{X_0}})$ donnée à la proposition \ref{baseH0k}, posons $h=h(\omega'_{i,j})-\iota(\overline{h(\omega'_{i,j})}) \in \mathcal{O}_{{X_0}}(\overline{\mathcal{U}}\cap \mathcal{V})$. Il existe $h_U$ dans $\mathcal{O}_{{X_0}}(\overline{\mathcal{U}})$ et $h_V$ dans $\mathcal{O}_{{X_0}}(\mathcal{V})$ tels que $h=h_U+h_V$. Enfin notons $(\beta_U^{ij},-\beta_V^{ij},\iota(\overline{h(\omega'_{i,j})}))=\tau (\overline{h(\omega'_{i,j})})$, alors 
\[P_1\left( \Phi(\omega_{i,j})\right) =\overline{h(\omega'_{i,j})}\]
\[P_2\left( \Phi(\omega_{i,j})\right) =f_V(\omega'_{i,j})+\beta_V^{ij}+dh_V.\]
\end{prop}

\begin{pf} On a vu à la section \ref{rappels} que $\Phi(\omega_{i,j})=\Phi'(\omega'_{ij})$ et $\Phi(\overline{h_{i,j}})=\Phi'(\overline{h'_{i,j}})$. Or sur $H^0(X'_0,\Omega_{X'_0}^1)$ on a $\Phi'(\omega)=\left[ DI_1(\omega)\right] $ donc $P_1(\Phi(\omega_{i,j}))=P_1(\Phi'(\omega'_{i,j}))=\overline{h(\omega'_{i,j})}$. En ce qui concerne le second point il s'agit de calculer la classe de $DI_1(\omega'_{i,j})$ dans $H^1($\v{C}$(\mathcal{F},\Omega_{{X_0}}^{\bullet}))$. En utilisant les notations de l'énoncé on a 
\begin{align*}
DI_1(\omega'_{i,j})=&(f_U(\omega'_{ij}),f_V(\omega'_{ij}),h(\omega'_{ij}))= (\beta_U^{ij},-\beta_V^{ij},\iota(\overline{h(\omega'_{i,j})})) + (dh_U,-dh_V,h_U+h_V) \\
& + (f_U(\omega'_{ij})-\beta_U^{ij}-dh_U,f_V(\omega'_{ij})+\beta_V^{ij}+dh_V,0)
\end{align*}
Or le triplet $(dh_U,-dh_V,h_U+h_V)$ est clairement un bord du complexe de \v{C}ech puisqu'il est l'image de $(h_U,-h_V)$ par l'application 
\begin{align*}
\mathcal{O}_{{X_0}}(\overline{\mathcal{U}})\oplus \mathcal{O}_{{X_0}}(\mathcal{V})\rightarrow &\Omega^1_{{X_0}}(\overline{\mathcal{U}})\oplus\Omega^1_{{X_0}}(\mathcal{V})\oplus\mathcal{O}_{{X_0}}(\overline{\mathcal{U}}\cap \mathcal{V})\\
(a,b)\mapsto &\left( da,db,(a-b)\vert_{\mathcal{O}_{{X_0}}(\overline{\mathcal{U}}\cap \mathcal{V})}\right) 
\end{align*}
et est donc de classe nulle dans $H^1($\v{C}$(\mathcal{F},\Omega_{{X_0}}^{\bullet}))$. En ce qui concerne le dernier triplet on a 
\[f_U(\omega'_{ij})-\beta_U^{ij}-dh_U-(f_V(\omega'_{ij})+\beta_V^{ij}+dh_V)=0\] 
ce qui signifie que $f_V(\omega'_{ij})+\beta_V^{ij}+dh_V$ (resp. $f_U(\omega'_{ij})-\beta_U^{ij}-dh_U$) définit une section globale de $H^0({X_0},\Omega_{{X_0}}^1)$ qui est donc l'image recherchée de la projection de $\Phi(\omega_{i,j})$ sur $H^0({X_0},\Omega_{{X_0}}^1)$. Ce qui achève la démonstration. Par ailleurs on a choisi d'exprimer cette section globale sur l'ouvert $\mathcal{V}$ pour des questions de calcul pratique qui seront détaillées plus loin. \\
\end{pf}

\begin{defi} \label{lesprojecteursquivontbien}
Par définition $\overline{\mathcal{U}}=\mathcal{U}\cap D(y)$ donc la décomposition en somme directe de $\Omega^1_{{X_0}}(\overline{\mathcal{U}})$ et $\Omega^1_{{X_0}}(\mathcal{V})$ donnée par le lemme \ref{lemmepourprojeteromega} permet de définir les projecteurs suivants :
\[q_U : \Omega^1_{{X_0}}(\overline{\mathcal{U}}) \rightarrow H^0({X_0},\Omega^1_{{X_0}})\]
\[q_V : \Omega^1_{{X_0}}(\mathcal{V}) \rightarrow H^0({X_0},\Omega^1_{{X_0}})\]
qui vont servir par la suite. \\
\end{defi}

\begin{rem} \label{remarquepourtronquersurW}
On sait que $f_V(\omega'_{i,j})+\beta_V^{ij}+dh_V$ est une section globale de $H^0({X_0},\Omega_{{X_0}}^1)$ par la démonstration de la proposition précédente. Ce n'est pas le cas de $f_V(\omega'_{i,j})$, $\beta_V^{ij}$ et $dh_V$ qui sont seulement des éléments de $\Omega_{{X_0}}^1(\mathcal{V})$. Cependant on a 
\[f_V(\omega'_{i,j})+\beta_V^{ij}+dh_V=q_V(f_V(\omega'_{i,j})+\beta_V^{ij}+dh_V)=q_V(f_V(\omega'_{i,j}))+q_V(\beta_V^{ij})+q_V(dh_V)\]
ce qui signifie que pour les calculs effectifs il suffit de calculer les projections de $f_V(\omega'_{i,j})$, $\beta_V^{ij}$ et $dh_V$ sur $H^0({X_0},\Omega^1_{{X_0}})$. \\
\end{rem}

\begin{prop} \label{derncol}
Soit $\overline{h_{i,j}}=\overline{t^{-i}y^j}$ un élément de la base de $H^1({X_0},\mathcal{O}_{{X_0}})$ donnée à la proposition \ref{baseH1Ok}, posons $h=h_{i,j}^p-\iota(\overline{h_{i,j}^p}) \in \mathcal{O}_{{X_0}}(\overline{\mathcal{U}}\cap \mathcal{V})$. Il existe $h_U$ dans $\mathcal{O}_{{X_0}}(\overline{\mathcal{U}})$ et $h_V$ dans $\mathcal{O}_{{X_0}}(\mathcal{V})$ tels que $h=h_U+h_V$. Enfin notons $(\beta_U^{ij},-\beta_V^{ij},\iota(\overline{h_{i,j}^p}))=\tau (\overline{h_{i,j}^p})$, alors 
\[P_1\left( \Phi(\overline{h_{i,j}})\right) =\overline{h_{i,j}^p}\]
\[P_2\left( \Phi(\overline{h_{i,j}})\right) =-\beta_U^{ij}-dh_U.\]
\end{prop}

\begin{pf}
On utilise à nouveau le fait que $\Phi(\omega_{i,j})=\Phi'(\omega'_{ij})$ et $\Phi(\overline{h_{i,j}})=\Phi'(\overline{h'_{i,j}})$. Par le lemme \ref{pourledi0}, sur $H^1({X_0},\mathcal{O}_{{X_0}})$ on a $\Phi'(\overline{h_{i,j}'})=\left[ (0,0,h_{i,j}^p)\right] $ donc $P_1\left( \Phi(\overline{h_{i,j}})\right) =\overline{h_{i,j}^p}$. En ce qui concerne le second point il s'agit à nouveau de calculer la classe de $DI_0(\overline{h'_{i,j}})$ dans $H^1($\v{C}$(\mathcal{F},\Omega_{{X_0}}^{\bullet}))$. Toujours à l'aide du lemme \ref{pourledi0} et en utilisant les notation de l'énoncé on a donc : 
\begin{align*}
(0,0,h_{i,j}^p)
= (\beta_U^{ij},-\beta_V^{ij},\iota(\overline{h_{i,j}^p})) + (dh_U,-dh_V,h_U+h_V) 
+(-\beta_U^{ij}-dh_U,\beta_V^{ij}+dh_V,0).
\end{align*}
A nouveau le triplet $(dh_U,-dh_V,h_U+h_V)$ est un bord du complexe de \v{C}ech et est donc de classe nulle dans $H^1($\v{C}$(\mathcal{F},\Omega_{{X_0}}^{\bullet}))$. En ce qui concerne le dernier triplet on a $-\beta_U^{ij}-dh_U-(\beta_V^{ij}+dh_V)=0$ ce qui signifie donc que $-\beta_U^{ij}-dh_U$ (resp. $\beta_V^{ij}+dh_V$) définit une section globale de $H^0({X_0},\Omega_{{X_0}}^1)$ qui est donc l'image recherchée de la projection de $\Phi(\overline{h_{i,j}})$ sur $H^0({X_0},\Omega_{{X_0}}^1)$. Ce qui achève la démonstration. \\
\end{pf}

\begin{rem} \label{remarquepourtronquersurU}
On sait que $-\beta_U^{ij}-dh_U$ est une section globale de $H^0({X_0},\Omega_{{X_0}}^1)$ grâce à la démonstration de la proposition précédente. Ce n'est pas le cas de $\beta_U^{ij}$ et $dh_U$ qui sont seulement des éléments de $\Omega_{{X_0}}^1(\overline{\mathcal{U}})$. Mais 
\[-\beta_U^{ij}-dh_U=q_U(-\beta_U^{ij}-dh_U)=-q_U(\beta_U^{ij})-q_U(dh_U)\]
et pour les calculs effectifs il suffit de calculer les projections de $\beta_U^{ij}$ et $dh_U$ sur $H^0({X_0},\Omega_{{X_0}}^1)$. \\
\end{rem}

\begin{rem}
Pour pouvoir effectivement calculer la matrice de $\Phi$ il ne reste plus qu'à expliciter pour chaque élément de $\mathcal{O}_{{X_0}} (\overline{\mathcal{U}}\cap \mathcal{V})$ sa classe d'équivalence dans $H^1({X_0},\mathcal{O}_{{X_0}})$. \\
\end{rem}

\subsection{Calcul des classes d'éléments de $\mathcal{O}_{{X_0}} (\overline{\mathcal{U}}\cap \mathcal{V})$}

On considère l'application
\[ 
\begin{array}{c c c c}
      d : & \mathcal{O}_{{X_0}}(\overline{\mathcal{U}})\times \mathcal{O}_{{X_0}}(\mathcal{V}) & \rightarrow & \mathcal{O}_{{X_0}} (\overline{\mathcal{U}}\cap \mathcal{V}) \\
      & (a,b) & \mapsto & (a-b)\vert_{\overline{\mathcal{U}}\cap \mathcal{V}} \\
\end{array}
\]
alors $H^1({X_0},\mathcal{O}_{{X_0}})=\mathcal{O}_{{X_0}} (\overline{\mathcal{U}}\cap \mathcal{V}) /Im(d)$. On commence par expliciter $d$ sur les éléments de la base de $\mathcal{O}_{{X_0}}(\overline{\mathcal{U}})\times \mathcal{O}_{{X_0}}(\mathcal{V})$. Pour cela montrons que toutes les algèbres considérées se plongent dans l'algèbre $A=k\llbracket t\rrbracket [t^{-1}][y]/(y^n-f(t))$ (où $k\llbracket t\rrbracket [t^{-1}]$ est l'algèbre des séries de Laurent). \\

\begin{rem}
Comme $k\llbracket t\rrbracket$-module, \[A=\bigoplus_{j=0}^{n-1} k\llbracket t\rrbracket[t^{-1}]y^j.\] \\
\end{rem}

\begin{lem}
On a $\mathcal{O}_{{X_0}}(\overline{\mathcal{U}})\subset A$ et plus précisément
\[\mathcal{O}_{{X_0}}(\overline{\mathcal{U}}) \subset \bigoplus_{j=0}^{n-1} k\llbracket t \rrbracket y^j\] 
comme $k\llbracket t\rrbracket$-module. 
\end{lem}

\begin{pf}
Comme $\overline{\mathcal{U}}=\mathcal{U}\cap D(y)=\mathcal{U}\cap D(f(t))$, alors $\overline{\mathcal{U}}=\Spec \left( k[t,y,f(t)^{-1}]/y^n-f(t)\right) $. Par ailleurs, $f(t)$ admet un inverse dans $k\llbracket t \rrbracket$ (car $t$ et $f(t)$ sont premiers dans $k[t]$ par hypothèse) et on a donc 
\[\mathcal{O}_{{X_0}} (\overline{\mathcal{U}}) \subset k\llbracket t \rrbracket [y]/(y^n-f(t))\] 
ce qui donne clairement le résultat annoncé (et donc l'inclusion dans $A$). \\
\end{pf}

\begin{lem} \label{serieWk}
On a $\mathcal{O}_{{X_0}}(\mathcal{V})\subset A$ et plus précisément
\[\mathcal{O}_{{X_0}}(\mathcal{V}) =\bigoplus_{j=0}^{n-1} k[t^{-1}]t^{-rj}y^j \]
et
\[\mathcal{O}_{{X_0}}(\mathcal{V})=\bigoplus_{j=0}^{n-1} \bigoplus_{\alpha\geq rj}kt^{-\alpha}y^j\] 
comme $k$-espace vectoriel.
\end{lem}

\begin{pf}
Par définition $\mathcal{V}=\Spec\left( k[s,z]/(z^n-f_2(s))\right) $ d'où
\[\mathcal{O}_{{X_0}}(\mathcal{V})=\bigoplus_{j=0}^{n-1} k[s]z^j.\]
Par ailleurs sur $\overline{\mathcal{U}}\cap \mathcal{V}$ on a les relations $s=t^{-1}$, $z=t^{-r}y$ et alors par abus de notations on peut écrire
\[\mathcal{O}_{{X_0}}(\mathcal{V})=\bigoplus_{j=0}^{n-1} k[t^{-1}]t^{-rj}y^j\]
ce qui donne la première égalité, la seconde en étant une conséquence immédiate. L'inclusion dans $A$ est alors évidente. \\
\end{pf}

\begin{lem} \label{jenaibesoin}
On a $\mathcal{O}_{{X_0}} (\overline{\mathcal{U}}\cap \mathcal{V})\subset A$ et plus précisément comme $k$-espace vectoriel on a 
\[\mathcal{O}_{{X_0}} (\overline{\mathcal{U}}\cap \mathcal{V})= \mathcal{O}_{{X_0}}  (\overline{\mathcal{U}}) \oplus \left(\bigoplus_{\alpha>0}kt^{-\alpha} \oplus \bigoplus_{j=1}^{n-1} \ \bigoplus_{\alpha=1}^{rj-1}kt^{-\alpha}y^j \oplus \ \bigoplus_{j=1}^{n-1} \bigoplus_{\alpha\geq rj}kt^{-\alpha}y^j\right) \]
ou encore 
\[\mathcal{O}_{{X_0}}  (\overline{\mathcal{U}}\cap \mathcal{V})= \mathcal{O}_{{X_0}}  (\overline{\mathcal{U}}) \oplus \widehat{\mathcal{O}_{{X_0}}(\mathcal{V})} \oplus  \bigoplus_{j=1}^{n-1} \bigoplus_{\alpha=1}^{rj-1}kt^{-\alpha}y^j \]
où l'on a posé 
\[\widehat{\mathcal{O}_{{X_0}}(\mathcal{V})}=\bigoplus_{\alpha>0}kt^{-\alpha} \oplus \bigoplus_{j=1}^{n-1} \bigoplus_{\alpha\geq rj}kt^{-\alpha}y^j.\]
De plus $\mathcal{O}_{{X_0}}(\mathcal{V}) = k\oplus \widehat{\mathcal{O}_{{X_0}}(\mathcal{V})}$. 
\end{lem}

\begin{pf}
Par définition de ouverts affines $\overline{\mathcal{U}}$ et $\mathcal{V}$ on a $\overline{\mathcal{U}}\cap \mathcal{V}=\Spec(k[t,t^{-1},f(t)^{-1},y]/(y^n-f(t)))$ donc $\mathcal{O}_{{X_0}}(\overline{\mathcal{U}}\cap \mathcal{V})=k[t,t^{-1},f(t)^{-1},y]/(y^n-f(t))$ et 
\[\mathcal{O}_{{X_0}}(\overline{\mathcal{U}}\cap \mathcal{V})\simeq \bigoplus_{j=0}^{n-1} k[t][f(t)^{-1}][t^{-1}]y^j\] 
comme $k[t]$-module. Mais $t$ et $f(t)$ sont premiers dans $k[t]$ donc $f(t)^{-1}$ est dans $k\llbracket t \rrbracket$ et $\mathcal{O}_{{X_0}}(\overline{\mathcal{U}}\cap \mathcal{V})$ est inclus dans $A$. En ce qui concerne les trois égalités, la dernière est une évidence au regard du lemme \ref{serieWk} et la seconde n'est qu'une réécriture de la première. Par ailleurs il est clair par définition des différents espaces que le membre de droite de la première égalité est bien inclus dans $\mathcal{O}_{{X_0}}(\overline{\mathcal{U}}\cap \mathcal{V})$ il reste donc à prouver l'inclusion réciproque. Posons $B=k\llbracket t \rrbracket [y]/(y^n-f(t))$, alors 
\[B\simeq \bigoplus_{j=0}^{n-1} k\llbracket t\rrbracket y^j\] 
comme $k\llbracket t\rrbracket$-module. Posons également $E=k[t,y]/(y^n-f(t))$, alors 
\[E\simeq \bigoplus_{j=0}^{n-1}k[t]y^j\] 
comme $k[t]$-module. Le complété de $E$ pour la topologie $t$-adique est $k\llbracket t \rrbracket \otimes_{k[t]} E$ car $k[t]$ est noethérien. Or $k\llbracket t \rrbracket$ est un $k[t]$-module plat et $E$ est séparé pour la topologie $t$-adique donc $E\hookrightarrow k\llbracket t \rrbracket \otimes_{k[t]} E\simeq B$. Par définition $\mathcal{O}_{{X_0}}(\overline{\mathcal{U}})=E[y^{-1}]=E[f(t)^{-1}]$ et alors $\mathcal{O}_{{X_0}}(\overline{\mathcal{U}})\hookrightarrow B[f(t)^{-1}]$ car la localisation est plate. Par ailleurs $B[f(t)^{-1}]\simeq B$ car $f(t)^{-1}$ est dans $k\llbracket t\rrbracket$ et $\mathcal{O}_{{X_0}}(\overline{\mathcal{U}})\hookrightarrow B$ par platitude de la localisation. Par définition $\mathcal{O}_{{X_0}}(\overline{\mathcal{U}}\cap \mathcal{V})=\mathcal{O}_{{X_0}}(\overline{\mathcal{U}})[t^{-1}]=E[f(t)^{-1}][t^{-1}]$ et $\mathcal{O}_{{X_0}}(\overline{\mathcal{U}}\cap \mathcal{V})\hookrightarrow B[t^{-1}]$.
Pour conclure montrons le lemme suivant : 
\begin{changemargin}{0.5cm}{0.5cm}
\begin{lem}$B\cap \mathcal{O}_{{X_0}}(\overline{\mathcal{U}} \cap \mathcal{V})=\mathcal{O}_{{X_0}}(\overline{\mathcal{U}})$. 
\end{lem} 

\begin{pf} Rappelons tout d'abord que l'on a les inclusions d'anneaux compatibles à l'action de $\delta$
\begin{center}
$\mathcal{O}_{{X_0}}(\overline{\mathcal{U}})\overset{i_1}{\subset} B \overset{i_2}{\subset} A$ \\
$\mathcal{O}_{{X_0}}(\overline{\mathcal{U}} \cap \mathcal{V})= \mathcal{O}_{{X_0}}(\overline{\mathcal{U}})[f^{-1}]\overset{i_3}{\subset} A$
\end{center}
notons $B_j=k\llbracket t\rrbracket y^j$ et $A_j=k\llbracket t\rrbracket [t^{-1}]y^j$ pour $1\leq j\leq n-1$, les sous-espaces de $A$ et $B$ tels que $B=\oplus B_j$ et $A=\oplus A_j$. Alors $i_1\left( \mathcal{O}_{{X_0}}(\overline{\mathcal{U}})_j\right) \subset B_j$, $i_2(B_j)\subset A_j$ et $i_3\left( \mathcal{O}_{{X_0}}(\overline{\mathcal{U}} \cap \mathcal{V})_j\right)\subset A_j$. Ainsi il suffit de montrer que pour tout $0\leq j\leq n-1$ on a $B_j \cap \mathcal{O}_{{X_0}}(\overline{\mathcal{U}} \cap \mathcal{V})_j= \mathcal{O}_{{X_0}}(\overline{\mathcal{U}})_j$ c'est-à-dire $k\llbracket t\rrbracket y^j \cap k[t,t^{-1},f^{-1}]y^j=k[t,f^{-1}]y^j$ ce qui revient à montrer que 
\begin{center}
$k\llbracket t\rrbracket \cap k[t,t^{-1},f^{-1}]=k[t,f^{-1}]$.
\end{center}
Or l'algèbre $k\llbracket t\rrbracket$ est munie de la valuation $t$-adique, $v_t$, qui s'étend à $k\llbracket t\rrbracket[t^{-1}]$. De plus $v_t(f)$ est non nulle car $f(0)$ est différent de $0$. Soit $u$ un élément de $k[t,f^{-1}]$, alors il s'écrit $u=af^{-n}$ avec $a$ dans $k[t]$ et $v_t(u)=v_t(af^{-n})=v_t(a)\geq 0$. Soit $w$ dans $k\llbracket t\rrbracket \cap k[t,t^{-1},f^{-1}]$, alors $w=af^{-n}$ avec $a$ dans $k[t,t^{-1}]$. Comme $w$ est un élément de $k\llbracket t\rrbracket$ alors $v_t(w)=v_t(a)\leq 0$ et donc $a$ est un élément de $k[t]$ ce qui signifie que $w$ est dans $k[t,f^{-1}]$ et montre donc le lemme. \\
\end{pf}
\end{changemargin}
Montrons maintenant l'inclusion 
\[\mathcal{O}_{{X_0}} (\overline{\mathcal{U}}\cap \mathcal{V})\subset \mathcal{O}_{{X_0}}  (\overline{\mathcal{U}}) \oplus \left(\bigoplus_{\alpha>0}kt^{-\alpha}  \oplus \bigoplus_{j=1}^{n-1} \ \bigoplus_{\alpha=1}^{rj-1}kt^{-\alpha}y^j \oplus \ \bigoplus_{j=1}^{n-1} \bigoplus_{\alpha\geq rj}kt^{-\alpha}y^j\right). \]
Considérons $h$ élément de $\mathcal{O}_{{X_0}}(\overline{\mathcal{U}}\cap \mathcal{V})$ et regardons-le comme élément de $A$ (puisque $\mathcal{O}_{{X_0}}(\overline{\mathcal{U}}\cap \mathcal{V})\subset A$), alors $\forall 0\leq j\leq n-1$, il existe $h_j$ élément de $k\llbracket t \rrbracket [t^{-1}]$ tel que $h=\sum h_jy^j$. Par définition, pour $j$ fixé on peut écrire 
\[h_j = \sum_{s=-a}^0 \lambda_s^jt^s + \sum_{s\geq 0} \lambda_s^jt^s\] 
avec $\lambda^j_i$ dans $k$. Posons 
\[H_j=\sum_{s=-a}^0 \lambda_s^jt^s\] 
\[K_j=\sum_{s\geq 0} \lambda_s^jt^s\] 
alors $h=\sum (H_j+K_j)y^j$. Or $\sum H_jy^j$ est un élément de $\mathcal{O}_{{X_0}}(\overline{\mathcal{U}}\cap \mathcal{V})$ et comme c'est également le cas de $h$ on en déduit que $\sum K_jy^j$ est dans $\mathcal{O}_{{X_0}}(\overline{\mathcal{U}}\cap \mathcal{V})$ et d'après le lemme $\sum K_jy^j$ est un élément de $\mathcal{O}_{{X_0}}(\overline{\mathcal{U}})$. On a donc décomposé un élément de $\mathcal{O}_{{X_0}}(\overline{\mathcal{U}}\cap \mathcal{V})$ selon le second membre de la première égalité de l'énoncé ce qui prouve la seconde inclusion voulue et achève donc la démonstration. \\
\end{pf}

\begin{rem} 
A l'aide de ce dernier lemme on retrouve clairement le résultat déjà établi à la proposition \ref{baseH1Ok} qui est l'expression d'une base de $H^1({X_0},\mathcal{O}_{{X_0}})$. Cependant on voit que pour obtenir la classe d'un élément du type $f(t)^{-a}t^by^j$ avec $a\in \mathbb{N}$, $b\in\mathbb{Z}$ et $1\leq j\leq n-1$ il est nécessaire de connaître l'expression de l'inverse de $f$ dans $k\llbracket t \rrbracket$. Cela pose plusieurs problèmes et il est notamment peu raisonnable d'espérer gérer une telle donnée algorithmiquement. Cependant il existe une astuce, due à Xavier Caruso que l'on détaille dans le lemme suivant. \\
\end{rem}

\begin{lem} \label{Bezout}
Rappelons que l'on note $\overline{...\vphantom{t}}$ les classes des éléments de $\mathcal{O}_{{X_0}}(\overline{\mathcal{U}}\cap \mathcal{V})$ dans $H^1({X_0},\mathcal{O}_{{X_0}})$, soient $a,b\in \mathbb{N}$, $1\leq j\leq n-1$ alors :
\[\overline{f(t)^{-a}t^by^j}=0\] 
\[\overline{f(t)^{-a}t^{-b}y^j}=\overline{I_b^at^{-b}y^j}\]
où $I_b$ est l'inverse de $f$ à l'ordre $b+1$ (donc dans $k[t]/t^{b+1}$) i.e. $f(t)I_b(t)=1+t^{b+1}P(t)$ avec $P\in k[t]$. 
\end{lem}

\begin{pf}
Il est évident au vu de ce qui précède que $f(t)^{-a}t^by^j$ est un élément de $\mathcal{O}_{{X_0}} (\overline{\mathcal{U}})$ ce qui donne le premier résultat. En ce qui concerne le second. Soit $I$ l'inverse de $f$ dans $k\llbracket t \rrbracket$ et soit $I_{b}$ sa troncature à l'ordre $b$. On a donc $I=I_{b} + t^{b+1}R(t)$ où $R(t)\in k\llbracket t \rrbracket$. Alors 
\[f(t)^{-a}t^{-b}y^j=I^at^{-b}y^j=(I_{b} + t^{b+1}R(t))^at^{-b}y^j.\] 
En utilisant la formule du binôme de Newton on obtient finalement 
\[f(t)^{-a}t^{-b}y^j=(I_b^a+t^{b+1}\tilde{R}(t))t^{-b}y^j=I_b^at^{-b}y^j+t\tilde{R}(t)y^j.\] 
Le second terme de la somme étant clairement dans $Im(d)$, la classe de l'élément considéré est la classe de $I_b^at^{-b}y^j$ d'où le résultat voulu. \\
\end{pf}

\section{Calcul de la matrice du Frobenius divisé de la fibre spéciale} \label{algorithmematrice}

Ce qui précède permet le calcul de la matrice du Frobenius divisé $$\Phi : H^0(X_0,\Omega_{X_0}^1)\oplus H^1(X_0,\mathcal{O}_{X_0})\rightarrow H^1_{DR}(X_0)$$ défini à la section \ref{rappels}. On a donné dans les propositions \ref{baseH0k} et \ref{baseH1Ok} des bases de $H^0(X_0,\Omega_{X_0}^1)$ et $H^1(X_0,\mathcal{O}_{X_0})$ qui par abus de notation donnent une base de $H^1_{DR}(X_0)$ (voir remarque \ref{pourlecalculdelamatrice}). La matrice donnée par les formules présentées plus loin est exprimée dans cette base, on donne donc la matrice de 
\[\Phi : H^0(X_0,\Omega_{X_0}^1)\oplus H^1(X_0,\mathcal{O}_{X_0})\rightarrow H^0(X_0,\Omega_{X_0}^1)\oplus H^1(X_0,\mathcal{O}_{X_0}).\]

\subsection{Présentation de la matrice} \label{presentationmatrice}

On fait le choix de présentation suivant, l'image des éléments de la base de $H^0(X_0,\Omega_{X_0}^1)$ se lit dans les $g$ premières colonnes de la matrice et l'image des éléments de la base de $H^1(X_0,\mathcal{O}_{X_0})$ dans les $g$ dernières colonnes. Cela permet de respecter la filtration. Soit maintenant un élément de la base de $H^0(X_0,\Omega_{X_0}^1)\oplus H^1(X_0,\mathcal{O}_{X_0})$ la projection de son image sur $H^0(X_0,\Omega_{X_0}^1)$ se lit dans les $g$ premières lignes et la projection de son image sur $H^1(X_0,\mathcal{O}_{X_0})$ dans les $g$ dernières lignes. Ceci ne fait pas apparaître le découpage par composantes isotypiques comme présenté dans la proposition \ref{lapropositionquifaut} et la remarque \ref{laremarquequifaut} mais si on découpe la matrice en $4$ blocs de taille $g\times g$ alors le bloc supérieur gauche correspond à la matrice de Cartier et le bloc inférieur droit à la matrice de Hasse-Witt. Par ailleurs pour chaque bloc on ordonne les éléments dans l'ordre lexicographique en prenant d'abord en compte la puissance de $y$ puis celle de $t$ afin de respecter les composantes isotypiques au sein de chaque bloc. En d'autres termes la matrice se présente de la manière suivante \\

\[
\begin{blockarray}{ccccccc}
   & y^{-1}dt & \dots & t^{r(n-1)-2}y^{-(n-1)}dt & t^{-1}y & \dots & t^{-(r(n-1)-1)}y^{n-1} \\
  \begin{block}{c(cccccc)}
  y^{-1}dt & * & \dots  & * & * & \dots  & *  \\
  \vdots  & * & \dots   & * & * & \dots  & *  \\
  t^{r-2}y^{-1}dt & * & \dots   & * & * & \dots  & *  \\
  y^{-2}dt  & * & \dots   & * & * & \dots  & *  \\
  \vdots & * & \dots   & * & * & \dots  & *  \\
  t^{r(n-1)-2}y^{-(n-1)}dt & * & \dots  & * & * & \dots  & *  \\
  t^{-1}y  & * & \dots   & * & * & \dots  & *  \\
  \vdots & * & \dots   & * & * & \dots  & *  \\
  t^{r-1}y & * & \dots  & * & * & \dots  & *  \\
  t^{-1}y^2 & * & \dots   & * & * & \dots  & *  \\
  \vdots & * & \dots  & * & * & \dots  & *  \\
  t^{-(r(n-1)-1)}y^{n-1} &  * & \dots & * & * & \dots  & *  \\
  \end{block}
\end{blockarray}
\] \\
\begin{defi*} 
On rappelle qu'on utilise les notations suivantes dans tout cette section. L'identification de $H^1({X_0},\mathcal{O}_{{X_0}})$ à $\mathcal{O}_{{X_0}}(\overline{\mathcal{U}}\cap \mathcal{V})/Im(d_{0,1})$ induit (voir remarque \ref{identH1Ok})
\begin{center}
$\iota:H^1({X_0},\mathcal{O}_{{X_0}})\rightarrow \mathcal{O}_{{X_0}}(\overline{\mathcal{U}}\cap \mathcal{V})$
\end{center}
ainsi que 
\begin{center}
$\overline{...\vphantom{t}}:\mathcal{O}_{{X_0}}(\overline{\mathcal{U}}\cap \mathcal{V})\rightarrow H^1({X_0},\mathcal{O}_{{X_0}})$.
\end{center}
Par ailleurs on a construit à la proposition \ref{scindagemodp}
\begin{center}
$\tau : H^1({X_0},\mathcal{O}_{{X_0}})\rightarrow \Omega_{{X_0}}^1(\overline{\mathcal{U}}) \times\Omega_{{X_0}}^1(\mathcal{V})\times \mathcal{O}_{{X_0}}(\overline{\mathcal{U}}\cap \mathcal{V})$
\end{center}
et on note
\begin{center}
$\left[ \vphantom{t}...\right] : \Omega_{{X_0}}^1(\overline{\mathcal{U}}) \times\Omega_{{X_0}}^1(\mathcal{V})\times \mathcal{O}_{{X_0}}(\overline{\mathcal{U}}\cap \mathcal{V})\rightarrow H^1($\v{C}$(\mathcal{F},\Omega_{{X_0}}^{\bullet}))$. \\
\end{center}
On a également défini les projecteurs suivants (définition \ref{lesprojecteursquivontbien})
$$q_U : \Omega^1_{X_0}(\overline{\mathcal{U}}) \rightarrow H^0(X_0,\Omega^1_{X_0})$$
$$q_V : \Omega^1_{X_0}(\mathcal{V}) \rightarrow H^0(X_0,\Omega^1_{X_0})$$
et on notera enfin $P_1$ (resp. $P_2$) la projection de $H^0({X_0},\Omega^1_{{X_0}})\oplus H^1({X_0},\mathcal{O}_{{X_0}})$ sur le second (resp. le premier) facteur. \\
\end{defi*}

\subsection{Calcul de la matrice de Hasse-Witt}

La matrice de Hasse-Witt correspond au quart inférieur droit de la matrice du Frobenius divisé que nous sommes en train de calculer, c'est en fait $P_1\circ \Phi\vert_{H^1(X_0,\mathcal{O}_{X_0})} :H^1(X_0,\mathcal{O}_{X_0})\rightarrow H^1(X_0,\mathcal{O}_{X_0})$. En vertu de la proposition \ref{derncol} on l'obtient tout simplement en calculant la classe de $h_{i,j}^p$ dans $H^1(X_0,\mathcal{O}_{X_0})$. \\

\begin{prop} \label{HW}
Soit $\overline{h_{i,j}}$ un élément de la base de $H^1(X_0,\mathcal{O}_{X_0})$ (voir proposition \ref{baseH1Ok}). Soient $a$ le quotient et $b$ le reste de la division euclidienne de $pj$ par $n$. Notons $F(t)=f(t)^a$ et $F[k]$ le coefficient de $t^k$ de $F$. Alors $1\leq b\leq n-1$ et
\[P_1\left( \Phi\left( \overline{h_{i,j}}\right) \right)=\overline{h_{i,j}^p}=\sum_{k=1}^{rb-1}F[pi-k]\overline{t^{-k}y^b}=\sum_{k=1}^{rb-1}F[pi-k]\overline{h_{k,b}} \]
dans $H^1(X_0,\mathcal{O}_{X_0})$. 
\end{prop}

\begin{pf}
Par définition $0\leq b\leq n-1$. Par ailleurs $b$ est non nul car $p$ est premier à $n$ et $1\leq j\leq n-1$ donc $n$ ne divise pas le produit $pj$. Soient $1\leq j\leq n-1$, $1\leq i\leq rj-1$, par abus de notation notons $DI_0(t^{-i}y^j)$ l'élément $(t^{-i}y^j)^p$ de $\mathcal{O}_{X_0}(\overline{\mathcal{U}})$ alors
\[DI_0(t^{-i}y^j)=t^{-pi}y^{pj}=t^{-pi}f(t)^ay^b=t^{-pi}F(t)y^b=\sum_{k=0}^d F[k]t^{-(pi-k)}y^b\]
où $d$ désigne le degré de $F$ en tant que polynôme en $t$. Réalisons le changement de variable $k'=pi-k$, alors
\[DI_0(t^{-i}y^j)=\sum_{k'=pi-d}^{pi} F[pi-k']t^{-k'}y^b.\]
Or pour tout $\alpha$ dans $\mathbb{N}$, $\overline{t^{\alpha}y^b}=0$ de plus $\overline{t^{-\alpha}y^b}=\overline{h_{\alpha,b}}$ lorsque $1\leq \alpha \leq rb-1$ et vaut $0$ sinon. D'où
\[\overline{DI_0(t^{-i}y^j)}=\sum_{k'=1}^{rb-1} F[pi-k']\overline{t^{-k'}y^b}\]
ce qui est le résultat annoncé. Par ailleurs notons que si $pi-d>1$ ou $pi<rb-1$ la somme contient à priori plus de termes qu'elle ne devrait. Dans le premier cas un tel terme $k$ vérifie $pi-k>d$ et donc $F[pi-k]=0$ puisque $d=\deg(F)$. Dans le second cas $k$ vérifie $pi-k<0$ et donc $F[pi-k]=0$ puisque $F$ est un polynôme. \\
\end{pf}

\begin{rem}
Dans le cas hyperelliptique, on a $j=1$ et comme $p$ est impair alors $a=(p-1)/2$ et $b=1$, de plus $g=r-1$. Ainsi la formule devient
\[P_1\left( \Phi(\overline{t^{-i}y})\right) =\sum_{k=1}^{g}f^{(p-1)/2}[pi-k]\overline{t^{-k}y}\]
dans $H^1(X_0,\mathcal{O}_{X_0})$ et le coefficient en $\overline{t^{-j}y}$ vaut donc $f^{(p-1)/2}[pi-j]$. On retrouve bien le résultat cité dans \cite{MR3240808} (à transposition près). \\
\end{rem}

\subsection{Calcul du bloc supérieur droit}

Ce bloc correspond à $P_2\circ\Phi\vert_{H^1(X_0,\mathcal{O}_{X_0})} :H^1(X_0,\mathcal{O}_{X_0})\rightarrow H^0(X_0,\Omega_{X_0}^1)$. On utilise la proposition \ref{derncol} pour calculer $P_2\left( \Phi\left( \overline{h_{i,j}}\right) \right)=-\beta_U^{ij}-dh_U$. Comme $q_U\left( P_2\left( \Phi\left( \overline{h_{i,j}}\right) \right) \right)=P_2\left( \Phi\left( \overline{h_{i,j}}\right) \right)$, alors 
\[P_2\left( \Phi\left( \overline{h_{i,j}}\right) \right)=-q_U\left( \beta_U^{ij}\right) -q_U\left( dh_U\right).\] 
Dans ce qui suit nous procédons au calcul de ces deux termes. \\

\begin{lem} \label{betaU}
Soit $\overline{h_{i,j}}$ un élément de la base de $H^1(X_0,\mathcal{O}_{X_0})$ (voir proposition \ref{baseH1Ok}). Soient $a$ le quotient et $b$ le reste de la division euclidienne de $pj$ par $n$. Notons $F(t)=f(t)^a$, $F[k]$ le coefficient de $t^k$ de $F$ et $\beta_{ij}^U\in \Omega_{X_0}^1(\overline{\mathcal{U}})$ la première composante de $\tau( \overline{h_{i,j}^p})$ (voir proposition \ref{scindagemodp}). Alors $1\leq b\leq n-1$ et
\[q_U\left( \beta_{ij}^U \right) = \sum_{m=0}^{r(n-b)-2}\left( \sum_{k=1}^{rb-1} \dfrac{b(m+k+1)-kn}{n}F[pi-k]f[m+k+1]\right) \dfrac{t^mdt}{y^{n-b}}\]
dans $H^0(X_0,\Omega_{X_0}^1)$. 
\end{lem}

\begin{pf}
En ce qui concerne la condition sur $b$ l'argument est le même que dans la démonstration précédente. Par ailleurs la formule proposée vient tout simplement de la linéarité du scindage $\tau$ et de la formule donnée à la proposition \ref{HW}. En effet le terme que nous cherchons à calculer est le $\beta_{ij}^U$ de la proposition \ref{derncol} qui est donné par la première composante de $\tau(\overline{t^{-k}y^b})$ et qui vaut 
\[\sum_{m=k+1}^l  \dfrac{bm-kn}{n}f[m]t^{m-(k+1)}\dfrac{dt}{y^{n-b}}\]
par la proposition \ref{scindagemodp}. On a donc 
\begin{align*}
\beta_{ij}^U =&\sum_{k=1}^{rb-1} \sum_{m=k+1}^l \dfrac{bm-kn}{n}F[pi-k]f[m] \dfrac{t^{m-(k+1)}dt}{y^{n-b}} \in \Omega_{X_0}^1(\overline{\mathcal{U}})
\end{align*}
et en réalisant le changement de variable $m'=m-(k+1)$ on obtient
\begin{align*}
\beta_{ij}^U =&\sum_{k=1}^{rb-1} \sum_{m'=0}^{l-(k+1)} \dfrac{b(m'+k+1)-kn}{n}F[pi-k]f[m'+k+1] \dfrac{t^{m}dt}{y^{n-b}}.
\end{align*}
Par la remarque \ref{remarquepourtronquersurU} il suffit de calculer $q_U\left( \beta_{ij}^U \right)$ et pour cela on considère uniquement les termes de cette somme tels que $0\leq m'\leq r(n-b)-2$ (voir lemme \ref{lemmepourprojeteromega}) ce qui achève la démonstration. Encore une fois on pourrait penser que si $r(n-b)-2>l-(k+1)$ la somme contient plus de termes qu'elle ne devrait. Cependant un tel terme vérifie $m'>l-(k+1)$ donc $m'+k+1>l$ et $f[m'+k+1]=0$ puisque $\deg(f)=l$. \\
\end{pf}

\begin{lem} \label{dhU}
Soit $\overline{h_{i,j}}$ un élément de la base de $H^1(X_0,\mathcal{O}_{X_0})$ (voir proposition \ref{baseH1Ok}). Soient $a$ le quotient et $b$ le reste de la division euclidienne de $pj$ par $n$. Notons $F(t)=f(t)^a$, $F[k]$ le coefficient de $t^k$ de $F$ et $d=\deg(F)$. Soit $h=h^p_{i,j}-\iota(\overline{h^p_{i,j}})$ et notons $h_U\in\mathcal{O}_{X_0}(\overline{\mathcal{U}})$ et $h_V\in\mathcal{O}_{X_0}(\mathcal{V})$ tels que $h=h_U+h_V$ (voir proposition \ref{derncol}). Alors $1\leq b\leq n-1$ et
\[q_U\left( dh_U\right)  = \sum_{m=0}^{r(n-b)-2}\left( \sum_{k=pi-d}^{0} \dfrac{b(m+k+1)-kn}{n}F[pi-k]f[m+k+1]\right) \dfrac{t^mdt}{y^{n-b}}\]
dans $H^0(X_0,\Omega_{X_0}^1)$. Où par convention la somme est nulle si $pi-d>0$. 
\end{lem}

\begin{pf}
Pour la condition sur $b$ l'argument est toujours le même. Dans la démonstration de la proposition \ref{HW} on a établi la formule suivante 
\[DI_0(t^{-i}y^j)=\sum_{k=pi-d}^{pi} F[pi-k]t^{-k}y^b\]
et on a donc
\[h=\sum_{k=pi-d}^{0} F[pi-k]t^{-k}y^b+\sum_{k=rb}^{pi} F[pi-k]t^{-k}y^b\]
et clairement $h_U$ correspond à la première somme et $h_V$ à la seconde grâce au lemme \ref{jenaibesoin}. Dans le cas où $pi-d>0$ la première somme est vide donc $h_U=0$ et $dh_U=0$. Dans le cas contraire on obtient par différentiation
\[dh_U=\sum_{k=pi-d}^{0} F[pi-k]\left( -kt^{-k-1}y^bdt+t^{-k}by^{b-1}dy\right) \in \Omega_{X_0}^1(\overline{\mathcal{U}})\]
cette relation étant valable modulo $ny^{n-1}dy-f'(t)dt$ puisque sur $\overline{\mathcal{U}}$ on a la relation $y^n=f(t)$. Ainsi
\begin{align*}
dh_U=&\sum_{k=pi-d}^{0} F[pi-k]\left( -kt^{-k-1}y^ny^{b-n}dt+t^{-k}by^{b-1}\dfrac{f'(t)dt}{ny^{n-1}}\right) \\
 =&\sum_{k=pi-d}^{0} F[pi-k]\left( -kt^{-k-1}f(t)+\dfrac{bt^{-k}f'(t)}{n}\right) \dfrac{dt}{y^{n-b}}.
\end{align*}
On remplace alors $f(t)$ et $f'(t)$ par leur expression en fonction des $f[k]$ et on obtient
\begin{align*}
dh_U=&\sum_{k=pi-d}^{0} F[pi-k]\left( -kt^{-k-1}\sum_{m=0}^l f[m]t^m+\dfrac{bt^{-k}}{n}\sum_{m=0}^l mf[m]t^{m-1}\right) \dfrac{dt}{y^{n-b}} \\
=& \sum_{k=pi-d}^{0} F[pi-k]\left( -k\sum_{m=0}^l f[m]t^{m-k-1}+\dfrac{b}{n}\sum_{m=0}^l mf[m]t^{m-1-k}\right) \dfrac{dt}{y^{n-b}}\\
=&\sum_{k=pi-d}^{0} \sum_{m=0}^l \dfrac{bm-kn}{n}F[pi-k]f[m]\dfrac{t^{m-k-1}dt}{y^{n-b}}.
\end{align*}
On réalise le changement de variable $m'=m-k-1$, alors
\[dh_U=\sum_{k=pi-d}^{0} \sum_{m'=-(k+1)}^{l-(k+1)} \dfrac{b(m'+k+1)-kn}{n}F[pi-k]f[m'+k+1]\dfrac{t^{m'}dt}{y^{n-b}}.\]
Par la remarque \ref{remarquepourtronquersurU} il suffit de calculer $q_U\left( dh_U\right)$ et pour cela on considère uniquement les termes de cette somme tels que $0\leq m'\leq r(n-b)-2$ (voir lemme \ref{lemmepourprojeteromega}) ce qui achève la démonstration. De plus si $r(n-b)-2>l-(k+1)$ la somme ne contient pas plus de termes qu'elle ne devrait puisqu'on aurait $m'>l-(k+1)$ et donc $f[m'+k+1]=0$ car $\deg(f)=l$. \\
\end{pf}

\begin{prop}\label{supdroit}
Soit $\overline{h_{i,j}}$ un élément de la base de $H^1(X_0,\mathcal{O}_{X_0})$ (voir proposition \ref{baseH1Ok}). Soient $a$ le quotient et $b$ le reste de la division euclidienne de $pj$ par $n$. Notons alors $F(t)=f(t)^a$, $F[k]$ le coefficient de $t^k$ de $F$ et $d=\deg(F)$. Alors $1\leq b\leq n-1$ et
\[P_2\left( \Phi\left( \overline{h_{i,j}}\right) \right)=\sum_{m=0}^{r(n-b)-2}\left( \sum_{k=pi-d}^{rb-1}\dfrac{kn-b(m+k+1)}{n}F[pi-k]f[m+k+1]\right) \dfrac{t^mdt}{y^{n-b}}\]
dans $H^0(X_0,\Omega_{X_0}^1)$. 
\end{prop}

\begin{pf}
Cette proposition est tout simplement une conséquence de la proposition \ref{derncol} et des lemmes \ref{betaU} et \ref{dhU}. Remarquons que dans le lemme \ref{dhU} on avait dû considérer séparément les cas où $pi-d\leq 0$ et $pi-d>0$. Dans le premier cas il est clair qu'on a la somme proposée, dans le second cas on a $dh_U=0$ mais en fait cela signifie également que la formule donnée dans le lemme \ref{betaU} présente des termes nuls jusqu'à ce que $k=pi-d$ ce qui fait qu'on obtient bien la même formule dans les deux cas. \\
\end{pf}

\subsection{Calcul du bloc inférieur gauche}

Ce bloc correspond à $P_1\circ\Phi\vert_{H^0(X_0,\Omega_{X_0}^1)} :H^0(X_0,\Omega_{X_0}^1)\rightarrow H^1(X_0,\mathcal{O}_{X_0})$. On utilise donc la proposition \ref{premcol}. Calculons la classe de $h(\omega'_{i,j})$ dans $H^1(X_0,\mathcal{O}_{X_0})$ pour $1\leq j\leq n-1$ et $0\leq i\leq rj-2$. \\

\begin{lem} \label{approx}
Soit $\omega'_{i,j}$, $1\leq j\leq n-1$, $0\leq i\leq rj-2$ un élément de la base de $H^0(X'_0,\Omega_{X'_0}^1)$ (voir proposition \ref{basesprimeesH0}). Soit $v(s)$ le polynôme intervenant dans le relèvement du Frobenius modulo $p^2$ sur $\mathcal{V}_1$ (voir proposition \ref{frob2}). Notons $d=max\left\lbrace 0,\deg(v)-2p\right\rbrace $. Alors 
\[\overline{h(\omega'_{i,j})}=\overline{t^{p(i+2)}v(s)I_d(t)^{a+1}y^{n-b}}\]
dans $H^1(X_0,\mathcal{O}_{X_0})$ où $I_d$ est l'inverse de $f$ à l'ordre $d+1$ (voir lemme \ref{Bezout}). 
\end{lem}

\begin{pf}
Tout d'abord on a vu à la remarque \ref{DI1explicite} que 
\[h(\omega'_{i,j})=\dfrac{t^{p(i+2)}}{f(t)^{a+1}}v(s)y^{n-b}.\]
Par ailleurs en vertu du lemme \ref{Bezout} on sait qu'il suffit de connaître la puissance négative maximale en $t$ pouvant apparaître dans cette expression et que celle-ci nous donne l'ordre d'inversion de $f$ suffisant pour calculer la classe de $h(\omega'_{i,j})$.
Comme $s=t^{-1}$ la puissance négative maximale en $t$ apparaissant dans $t^{p(i+2)}v(s)$ est donc $\deg(v)-p(i+2)$. Or $0\leq i\leq rj-1$ donc $\deg(v)-p(i+2)\geq \deg(v)-2p$, $\forall 0\leq i\leq rj-1$. Ainsi la puissance maximale en $t^{-1}$ est $\deg(v)-2p$ et ce quel que soit $j$. Soit donc $d=\deg(v)-2p$ alors connaître l'inverse de $f$ à l'ordre $d+1$ est suffisant pour obtenir la classe cherchée. \\
\end{pf}

\begin{rem}
On a déjà vu que puisqu'il était nécessaire de ne considérer qu'une estimation de l'inverse de $f$ il fallait mener les calculs non pas sur $\overline{\mathcal{U}}$ mais sur $\mathcal{V}$. On utilise donc les relations entre $(t,y)$ et $(s,z)$ qui sont pour mémoire $s=t^{-1}$ et $z=t^{-r}y$ où $r=(l+1)/n$. Clairement $\forall 1\leq j\leq n-1$, $\forall 1\leq i\leq rj-1$, $t^{-i}y^j=s^{-(rj-i)}z^j$ et $\forall 1\leq j\leq n-1$, $\forall 0\leq i\leq rj-2$, $t^iy^{-j}dt=-s^{rj-(i+2)}z^{-j}ds$. Donc la base de $H^1(X_0,\mathcal{O}_{X_0})$ peut s'écrire $s^{-i}z^j$ avec $1\leq j\leq n-1$ et $1\leq i\leq rj-1$ et celle de $H^0(X_0,\Omega_{X_0}^1)$ peut s'écrire $s^iz^{-j}ds$ avec $1\leq j\leq n-1$ et $0\leq i\leq rj-2$. \\
\end{rem}

\begin{prop} \label{DI1H1}
Soit $\omega_{i,j}\in H^0(X_0,\Omega_{X_0}^1)$ un élément de la base (voir proposition \ref{baseH0k}). Soient $v(s)$ le polynôme du relèvement du Frobenius modulo $p^2$ sur $\mathcal{V}_1$ (voir proposition \ref{frob2}), $d=max\left\lbrace 0,\deg(v)-2p\right\rbrace $, $I_d$ l'inverse de $f$ à l'ordre $d+1$ (voir lemme \ref{Bezout}), $a$ le quotient et $b$ le reste de la division euclidienne de $pj$ par $n$. Notons $IP(t)=I_d(t)^{a+1}$, $d^{\circ}IP=\deg(IP)$, $IP[k]$ le coefficient du monôme $t^k$ dans $IP(t)$ et $v[m]$ le coefficient du monôme $s^m$ dans $v(s)$. Alors $1\leq b\leq n-1$ et
\[P_1\left( \Phi(\omega_{i,j})\right) =\overline{h(\omega'_{i,j})}=\sum_{m=1}^{r(n-b)-1} \left( \sum_{k=0}^{d^{\circ}IP}IP[k]\cdot v[p(i+2)+k+m]\right) \overline{t^{-m}y^{n-b}}\]
dans $H^1(X_0,\mathcal{O}_{X_0})$. 
\end{prop}

\begin{pf}
Tout d'abord on rappelle qu'en vertu du lemme \ref{approx} on sait que 
\[\overline{h(\omega'_{i,j})}=\overline{t^{p(i+2)}v(s)I_d(t)^{a+1}y^{n-b}}=\overline{t^{p(i+2)}v(s)IP(t)y^{n-b}}.\]
Par ailleurs
\begin{center}
$\displaystyle v(s)=\sum_{m=0}^{d+2p}v[m]s^m$ \ \ \ et  \ \ \ $\displaystyle IP(t)=\sum_{k=0}^{d^{\circ}IP}IP[k]t^k=\sum_{k=0}^{d^{\circ}IP}IP[k]s^{-k}$
\end{center}
D'où
\[t^{p(i+2)}IP(t)v(s)y^{n-b}=\sum_{m=0}^{d+2p}\sum_{k=0}^{d^{\circ}IP} v[m]\cdot IP[k]\cdot s^{m-k-p(i+2)}s^{-r(n-b)}z^{n-b} \in \mathcal{O}_{X_0}(\overline{\mathcal{U}}\cap \mathcal{V})\]
et en faisant le changement de variable $m'=p(i+2)+r(n-b)+k-m$ on a alors
\[t^{p(i+2)}IP(t)v(s)y^{n-b}=\sum_{k=0}^{d^{\circ}IP} \sum_{m'=pi+r(n-b)+k-d}^{p(i+2)+r(n-b)+k}IP[k]\cdot v[p(i+2)+r(n-b)+k-m']\cdot s^{-m'}z^{n-b}.\]
Or tous les éléments $s^{-m'}z^{n-b}$ ne vérifiant pas $1\leq m'\leq r(n-b)-1$ ont une classe nulle dans $H^1(X_0,\mathcal{O}_{X_0})$. On obtient donc 
\[\overline{h(\omega'_{i,j})}=\sum_{k=0}^{d^{\circ}IP} \sum_{m'=1}^{r(n-b)-1}IP[k]\cdot v[p(i+2)+r(n-b)+k-m']\cdot \overline{s^{-m'}z^{n-b}}.\]
Comme pour les cas précédents, la somme ne contient pas plus de termes qu'elle ne devrait et ce même si on a $pi+r(n-b)+k-d>1$ ou $p(i+2)+r(n-b)+k<r(n-b)-1$. Dans le premier cas un tel terme vérifierait $p(i+2)+r(n-b)+k-m>2p+d=\deg(v)$ et son coefficient est donc nul. Dans le second cas on a $p(i+2)+r(n-b)+k-m<0$ et le coefficient est également nul puisque $v$ est un polynôme. On utilise alors la condition de recollement pour exprimer cette égalité selon $t$ et $y$ et on a 
\[\overline{h(\omega'_{i,j})}=\sum_{k=0}^{d^{\circ}IP} \sum_{m'=1}^{r(n-b)-1}IP[k]\cdot v[p(i+2)+r(n-b)+k-m']\cdot \overline{t^{m'}t^{-r(n-b)}y^{n-b}}\]
puis avec le changement de variable $m=r(n-b)-m'$ on obtient finalement
\[\overline{h(\omega'_{i,j})}=\sum_{k=0}^{d^{\circ}IP} \sum_{m=1}^{r(n-b)-1}IP[k]\cdot v[p(i+2)+k+m]\cdot \overline{t^{-m}y^{n-b}}\]
ce qui achève la démonstration. \\
\end{pf}

\subsection{Calcul de la matrice de Cartier}

La matrice de Cartier correspond au quart supérieur gauche de la matrice du Frobenius divisé que nous sommes en train de calculer, c'est en fait $P_2\circ\Phi\vert_{H^0(X_0,\Omega_{X_0}^1)} :H^0(X_0,\Omega_{X_0}^1)\rightarrow H^0(X_0,\Omega_{X_0}^1)$. On utilise donc la proposition \ref{premcol} pour calculer $P_2(\Phi(\omega_{i,j}))=f_V(\omega_{i,j}')+\beta_{ij}^V+dh_V$. Cependant on sait que $q_V\left( P_2(\Phi(\omega_{i,j}))\right) = P_2(\Phi(\omega_{i,j}))$ et donc
\[P_2(\Phi(\omega_{i,j}))=q_V\left( f_V(\omega_{i,j}')\right) +q_V\left( \beta_{ij}^V\right) +q_V\left( dh_V\right).\]
Dans ce qui suit on va donc procéder au calcul de ces trois termes. \\

\begin{lem} \label{fW}
Soit $\omega'_{i,j}\in H^0(X'_0,\Omega_{X'_0}^1)$ un élément de la base (voir proposition \ref{basesprimeesH0}). Soient $a$ le quotient et $b$ le reste de la division euclidienne de $pj$ par $n$. Soient $v(s)$ le polynôme du relèvement du Frobenius modulo $p^2$ sur $\mathcal{V}_1$ (voir proposition \ref{frob2}) et $Q(s)$ la fraction rationnelle $(s^{p-1}+v'(s))/f_2(s)^a$. Alors $1\leq b\leq n-1$, $Q$ est un polynôme en $s$ et si l'on note $Q[k]$ le coefficient du monôme $s^k$ dans $Q$ et $d^{\circ}Q$ le degré de $Q$, alors
\[q_V\left( f_V(\omega'_{i,j})\right) =\sum_{\mu=0}^{rb-2} Q[rb-2+p(i+2-rj)-\mu]\dfrac{t^{\mu}dt}{y^b}\]
dans $H^0(X_0,\Omega_{X_0}^1)$. Pour mémoire $f_V(\omega'_{i,j})$ désigne la deuxième composante de $DI_1(\omega'_{i,j})$ (voir proposition \ref{morphiDI}). 
\end{lem}

\begin{pf}
L'argument pour la condition sur $b$ est toujours le même. De plus on a vu à la remarque \ref{DI1explicite} que l'entier $a$ était toujours strictement inférieur à $p$ et donc en vertu de la proposition \ref{divisible} $Q$ est bien un polynôme en $s$.
On utilise la formule de la remarque \ref{DI1explicite} pour $f_V$ et on a
\[f_V(\omega'_{i,j})=-s^{p(rj-(i+2))}Q(s)\dfrac{ds}{z^b}=\sum_{\mu=0}^{d^{\circ}Q} -Q[\mu]s^{\mu}s^{p(rj-(i+2))}\dfrac{ds}{z^b}\]
ce qui donne par le changement de variable $\mu '=p(rj-(i+2))+\mu$ l'expression
\[f_V(\omega'_{i,j})=\sum_{\mu '=p(rj-(i+2))}^{p(rj-(i+2))+rb-2} -Q[\mu '-p(rj-(i+2))]s^{\mu '}\dfrac{ds}{z^b}.\]
Grâce au lemme \ref{lemmepourprojeteromega} on sait que pour calculer $q_V\left( f_V(\omega'_{i,j})\right)$ on ne garde que les valeurs de $\mu '$ pour lesquelles $0\leq \mu' \leq rb-2$ et on obtient ainsi
\[q_V\left( f_V(\omega'_{i,j})\right) =\sum_{\mu '=0}^{rb-2} -Q[\mu '-p(rj-(i+2))]s^{\mu '}\dfrac{ds}{z^b}\]
et cette somme ne contient pas de termes supplémentaires puisque $p(rj-(i+2))$ est positif ou nul et donc $p(rj-(i+2))+rb-2\geq rb-2$ de plus pour $\mu '$ compris entre $0$ et $p(rj-(i+2))$ on a $\mu '-p(rj-(i+2))<0$ et comme $Q$ est un polynôme $Q[\mu '-p(rj-(i+2))]=0$. En utilisant la condition de recollement on a 
\[q_V\left( f_V(\omega'_{i,j})\right)=\sum_{\mu '=0}^{rb-2} -Q[\mu '-p(rj-(i+2))]t^{-\mu '}\dfrac{-t^{-2}dt}{t^{-rb}y^b}=\sum_{\mu '=0}^{rb-2} Q[\mu '-p(rj-(i+2))]t^{rb-2-\mu '}\dfrac{dt}{y^b}\]
soit par le changement de variable $\mu = rb-2-\mu '$
\[q_V\left( f_V(\omega'_{i,j})\right)=\sum_{\mu =0}^{rb-2} Q[rb-2-p(rj-(i+2))-\mu ]\dfrac{t^{\mu}dt}{y^b}\]
ce qui achève la démonstration. \\
\end{pf}

\begin{lem} \label{betaW}
Soit $\omega'_{i,j}\in H^0(X'_0,\Omega_{X'_0}^1)$ un élément de la base (voir proposition \ref{basesprimeesH0}). Notons $\beta_{ij}^V$ la deuxième composante de $\tau( \overline{h(\omega'_{i,j}}))$ (voir proposition \ref{scindagemodp}). Alors
$$q_V\left( \beta_{ij}^V\right) = 0$$ 
dans $H^0(X_0,\Omega_{X_0}^1)$. 
\end{lem}

\begin{pf}
On sait par la proposition \ref{DI1H1} que dans $H^1(X_0,\mathcal{O}_{X_0})$ on a l'égalité
\[\overline{h(\omega'_{i,j})}=\sum_{m=1}^{r(n-b)-1} \left( \sum_{k=0}^{d^{\circ}IP}IP[k]v[p(i+2)+k+m]\right) \overline{t^{-m}y^{n-b}}\]
et alors par linéarité du scindage il suffit de calculer la deuxième composante de $\tau(\overline{t^{-m}y^{n-b}})$ qui est égale à $-\beta_{ij}^V$ où
\[\beta_{ij}^V=\sum_{\mu =0}^{m}\dfrac{mn-(n-b)\mu }{n}f[\mu ]s^{m+rb-1-\mu }\dfrac{ds}{z^b},\]
(voir proposition \ref{scindagemodp}) et donc
\[\beta_{ij}^V =\sum_{m=1}^{r(n-b)-1} \sum_{k=0}^{d^{\circ}IP}IP[k]v[p(i+2)+k+m]\sum_{\mu =0}^{m}\dfrac{mn-(n-b)\mu }{n}f[\mu ]s^{m+rb-1-\mu }\dfrac{ds}{z^b}.\]
On constate que puisque $\mu \leq m$ alors $m+rb-1-\mu \geq rb-1$ et donc dans cette somme la puissance en $s$ est toujours strictement supérieure à $rb-2$. Mais grâce au lemme \ref{lemmepourprojeteromega} on sait que pour calculer  $q_V\left( \beta_{ij}^V\right)$ on ne conserve que les termes dont la puissance en $s$ est comprise entre $0$ et $rb-2$ et on obtient donc le résultat annoncé. \\
\end{pf}

\begin{lem} \label{dhW}
Soit $\omega'_{i,j}\in H^0(X'_0,\Omega_{X'_0}^1)$ un élément de la base (voir proposition \ref{basesprimeesH0}). Soit $v(s)$ le polynôme intervenant dans le relèvement du Frobenius modulo $p^2$ sur $\mathcal{V}_1$ (voir proposition \ref{frob2}), notons $d=max\left\lbrace 0,\deg(v)-2p\right\rbrace $ et $I_d$ l'inverse de $f$ à l'ordre $d+1$ (voir lemme \ref{Bezout}). Soient $a$ le quotient et $b$ le reste de la division euclidienne de $pj$ par $n$, notons $IP(t)=I_d(t)^{a+1}$, $d^{\circ}IP=\deg(IP)$, $IP[k]$ le coefficient de $t^k$ dans $IP(t)$ et $v[m]$ le coefficient de $s^m$ dans $v(s)$. Soit $h=h(\omega'_{i,j})-\iota(\overline{h(\omega'_{i,j})})$ et notons $h_U \in \mathcal{O}_{\overline{\mathcal{U}}}$ et $h_V \in \mathcal{O}_{\mathcal{V}}$ tels que $h=h_U+h_V$ (voir proposition \ref{premcol}). Posons $\nu_{m,k}=p(i+2)+k+m$, alors $1\leq b\leq n-1$ et
\[q_V\left( dh_V\right) =\sum_{\mu =0}^{rb-2}\left( \sum_{k=0}^{d^{\circ}IP}\sum_{m=r(n-b)}^{d-(pi+k)} \dfrac{n(\mu +1)-b(\mu +m+1)}{n}IP[k]\cdot f[\mu +m+1]\cdot v[\nu_{m,k}]\right)\dfrac{t^{\mu}dt}{y^b} \]
dans $H^0(X_0,\Omega_{X_0}^1)$. Où par convention la somme est nulle si $r(n-b)>d-(pi+k)$. 
\end{lem}

\begin{pf}
Dans la démonstration de la proposition \ref{DI1H1} on a obtenu l'égalité suivante 
\[t^{p(i+2)}IP(t)v(s)y^{n-b}=\sum_{k=0}^{d^{\circ}IP} \sum_{m=pi+r(n-b)+k-d}^{p(i+2)+r(n-b)+k}IP[k]\cdot v[p(i+2)+r(n-b)+k-m]\cdot s^{-m}z^{n-b}.\]
Clairement si $pi+r(n-b)+k-d>0$ i.e $r(n-b)>d-(pi+k)$ cette somme ne contient aucun terme polynômial en $s$ et donc $h_V=0$ ce qui implique que $dh_V=0$. Dans le cas contraire notons on a alors 
\[h_V=\sum_{k=0}^{d^{\circ}IP} \sum_{m=pi+r(n-b)+k-d}^{0}IP[k]\cdot v[p(i+2)+r(n-b)+k-m]\cdot s^{-m}z^{n-b}\]
et en utilisant les conditions de recollement
\[h_V=\sum_{k=0}^{d^{\circ}IP} \sum_{m=pi+r(n-b)+k-d}^{0}IP[k]\cdot v[p(i+2)+r(n-b)+k-m]\cdot t^{m}t^{-r(n-b)}y^{n-b}.\]
Posons $m'=r(n-b)-m$ on a alors
\[h_V=\sum_{k=0}^{d^{\circ}IP} \sum_{m'=r(n-b)}^{d-(pi+k)}IP[k]\cdot v[p(i+2)+k+m']\cdot t^{-m'}y^{n-b}.\]
Par différentiation on obtient
\[dh_V=\sum_{k=0}^{d^{\circ}IP} \sum_{m'=d-(pi+k)}^{r(n-b)}IP[k]\cdot v[p(i+2)+k+m']\left(-m't^{-m'-1}y^{n-b}dt+(n-b)t^{-m'}y^{n-b-1}dy\right)\]
qui est donc un élément de $\Omega^1_{X_0}(\mathcal{V})$,
alors comme $y^n=f(t)$ et $ny^{n-1}dy=f'(t)dt$
\[dh_V=\sum_{k=0}^{d^{\circ}IP} \sum_{m'=d-(pi+k)}^{r(n-b)}IP[k]\cdot v[p(i+2)+k+m']\left( -m't^{-m'-1}f(t)+(n-b)t^{-m'}\dfrac{f'(t)}{n}\right) \dfrac{dt}{y^b}.\]
Or on a 
\begin{align*}
-m't^{-m'-1}f(t)+(n-b)t^{-m'}\dfrac{f'(t)}{n}=&\sum_{\mu =0}^l -m'f[\mu]t^{\mu-m'-1}+\sum_{\mu =0}^l \dfrac{n-b}{n}\mu f[\mu]t^{\mu-1-m'}\\
=&\sum_{\mu =0}^l \dfrac{\mu (n-b)-m'n}{n}f[\mu] t^{\mu-m'-1}
\end{align*}
d'où
\[dh_V=\sum_{k=0}^{d^{\circ}IP} \sum_{m'=d-(pi+k)}^{r(n-b)}\sum_{\mu =0}^l \dfrac{\mu (n-b)-m'n}{n}IP[k]\cdot v[p(i+2)+k+m']\cdot f[\mu] \dfrac{t^{\mu-m'-1} dt}{y^b}.\]
On réalise le changement de variable $\mu'=\mu-m'-1$ et on pose $\nu_{m',k}=p(i+2)+k+m'$, alors
\[dh_V=\sum_{k=0}^{d^{\circ}IP} \sum_{m'=d-(pi+k)}^{r(n-b)}\sum_{\mu' =-(m'+1)}^{l-(m'+1)} \dfrac{n(\mu' +1)-b(\mu' +m'+1)}{n}IP[k]\cdot v[\nu_{m,k}]\cdot f[\mu'+m'+1] \dfrac{t^{\mu'} dt}{y^b}.\]
Enfin grâce au lemme \ref{lemmepourprojeteromega} on sait que pour calculer $q_V\left( dh_V\right)$ on ne conserve que les termes tels que $0\leq \mu'\leq rb-2$ ce qui donne la formule annoncée. Par ailleurs il n'y a pas de termes supplémentaires car si $-(m'+1)>0$ ou $l-(m'+1)<rb-2$ alors dans le premier cas on aura un terme avec $\mu' <-(m'+1)$ et donc $\mu'+m'+1<0$ et le coefficient sera nul puisque $f$ est un polynôme. Dans le second cas $\mu'+m'+1>l$ et comme $\deg(f)=l$ le coefficient sera à nouveau nul. \\
\end{pf}

\begin{prop} \label{projection}
Soit $\omega_{i,j}\in H^0(X_0,\Omega_{X_0}^1)$ un élément de la base (voir proposition \ref{baseH0k}). Soient $v(s)$ le polynôme intervenant dans le relèvement du Frobenius modulo $p^2$ sur $\mathcal{V}_1$ (voir proposition \ref{frob2}), $d=max\left\lbrace 0,\deg(v)-2p\right\rbrace $ et $I_d$ l'inverse de $f$ à l'ordre $d+1$ (voir lemme \ref{Bezout}). Soient $a$ le quotient et $b$ le reste de la division euclidienne de $pj$ par $n$, notons $IP(t)=I_d(t)^{a+1}$, $d^{\circ}IP=\deg(IP)$, $IP[k]$ le coefficient du monôme $t^k$ dans $IP(t)$ et $v[m]$ le coefficient du monôme $s^m$ dans $v(s)$. Soit $Q(s)$ le polynôme $(s^{p-1}+v'(s))/f_2(s)^a$ (voir lemme \ref{fW}), notons $Q[k]$ le coefficient du monôme $s^k$ dans $Q$ et $d^{\circ}Q$ le degré de $Q$. Alors $1\leq b\leq n-1$ et si on note $\nu_{m,k}=p(i+2)+k+m$,
\begin{align*}
\begin{split}
P_2\left( \Phi(\omega_{i,j})\right) =\sum_{\mu =0}^{rb-2} & \left( \sum_{k=0}^{d^{\circ}IP}\sum_{m=r(n-b)}^{d-(pi+k)} \dfrac{n(\mu +1)-b(\mu +m+1)}{n}IP[k]\cdot f[\mu +m+1]\cdot v[\nu_{m,k}]  \right. \\
& \left.  + Q[rb-2+p(i+2-rj)-\mu]  \vphantom{\sum_{k=0}^{dPI}} \right) \dfrac{t^{\mu}dt}{y^b} 
\end{split}
\end{align*}
dans $H^0(X_0,\Omega_{X_0}^1)$. Où par convention la double somme est nulle si $r(n-b)>d-(pi+k)$. 
\end{prop}

\begin{pf}
Cette proposition est tout simplement une conséquence de la proposition \ref{premcol} et des lemmes \ref{fW}, \ref{betaW} et \ref{dhW}. \\
\end{pf}

\section{Complexité des formules}

Dans ce qui suit on va donner une approximation de coût du calcul des formules exposées dans la section précédente. Par ailleurs on s'intéresse au calcul de la complexité algébrique on ne comptera donc le nombre d'opérations dans $k$ ou dans $W_1$. \\

\subsection{Complexité pour la matrice de Hasse-Witt}

On utilise la formule de la proposition \ref{HW}. Tout d'abord remarquons que à $j$ fixé et pour $i$ quelconque entre $1$ et $rj-1$ le calcul de l'image de $\overline{h_{i,j}}$ par $P_1\circ \Phi$ fait seulement appel à un coefficient du polynôme $F$ qui ne dépend que de $j$. Une fois ce polynôme calculé on peut compléter les $rj-1$ colonnes de la matrice sans sur-coût. \\

En ce qui concerne le calcul de $F$ il faut tout d'abord calculer $a$ qui est égal à la partie entière de $pj/n$ puis mettre à la puissance $a$ (inférieur ou égal à $p$) un polynôme de degré $l$. Ce calcul a donc un coût en $O(pl)$.\\

Grâce à la toute première remarque il n'y a donc qu'à calculer le polynôme $F$ pour chaque valeurs de $j$ ce qui fait donc $n-1$ calculs à effectuer. Le coût total du calcul est donc en $O(pnl)$ enfin comme le genre de la courbe est donnée par la formule $g=(n-1)(l-1)/2$ on peut finalement dire que le coût de la matrice de Hasse-Witt est en $O(pg)$. \\

\begin{rem}
Comme on n'a pas besoin de tous les coefficients de $F$ mais seulement de certains un calcul à la Harvey-Sutherland (voir \cite{MR3240808} pourrait être intéressant. \\
\end{rem}

\subsection{Complexité pour le bloc supérieur droit}

On utilise la formule de la proposition \ref{supdroit}. Il faut tout d'abord commencer par calculer le polynôme $F$ à $j$ fixé. Nous avons déjà vu que ceci avait un coût en $O(pl)$. \\

Ensuite pour chaque valeur de $i$ compris entre $1$ et $rj-1$ il faut calculer la somme. Chaque terme de celle-ci comprend $8$ opérations. Le nombre maximum de termes est donné par la plus grande valeur de $rb-1$ et la plus petite valeur de $pi-d$ possible. Comme $i$ vaut au moins $1$ et $d$ au plus $pl$ alors $pi-d$ vaut au moins $p(1-l)$. En ce qui concerne $rb-1$ puisque $r$ vaut $(l+1)/n$ il vaut au plus $l$ il y a donc au maximum $l-p(1-l)+1$ termes. Le nombre d'opérations pour calculer la somme est donc en $O(pl)$. A $j$ fixé on calcule donc $F$ en $O(pl)$ puis tous les termes pour $i$ entre $1$ et $rj-1$ en $O(pl^2)$. \\

Pour calculer tous le bloc il faut encore faire varier $j$ qui prend $n-1$ valeurs différentes, on a donc un coût de calcul en $O(pnl^2)$ i.e. $O(pgl)$. \\

\subsection{Complexité pour le bloc inférieur gauche}

On utilise la formule de la proposition \ref{DI1H1}. Tout d'abord le calcul du relèvement du Frobenius modulo $p^2$ ainsi que le calcul de l'inverse de $f$ à l'ordre $d+1$ ne dépendent ni de $j$ ni de $i$. \\

En ce qui concerne le calcul de $v$ il s'obtient tout d'abord en calculant la relation de Bezout $v_0 f'_2+b_0 f_2=1$ modulo $p^2$ qui a un coût en le maximum des degrés des polynômes $f'_2$ et $f_2$ i.e en $l+1$ opérations dans $W_1$. Il s'agit ensuite de mettre $v_0$ à la puissance $p$ qui a un coût dépendant du degré de $v_0$ i.e. en $(l+1)p$. Enfin il s'agit de multiplier $v_0^p$ par le polynôme $\mathcal{Q}(s)=(f_2(s^p)-f_2^{\sigma}(s)^p/p$ qui a un coût en le maximum des degrés des deux polynômes i.e. inférieur ou égale à $(l+1)p$. En ce qui concerne le calcul du polynôme $\mathcal{Q}(s)$ il faut tout d'abord calculer $f_2^{\sigma}(s)$ qui a un coût en $l+1$ puis de le mettre à la puissance $p$ ce qui a un coût en $p(l+1)$. Il faut ensuite faire la différence entre deux polynômes de degré inférieur ou égal à $p(l+1)$ et diviser par $p$ ce qui dans $W_1$ a un coût $1$. Le calcul de $\mathcal{Q}$ a donc un coût en $O(pl)$ et le calcul de $v$ a un coût en $O(pl)$.\\

Pour le calcul de l'inverse, en utilisant un algorithme de Newton le coût est en le degré auquel on souhaite calculer l'inverse. Or le degré voulu est inférieur à $deg(v)-2p$ et le degré de $v$ étant au plus $2p(l+1)$ on a donc un coût du calcul de l'inverse en $O(pl)$. \\

En ce qui concerne les calculs dépendants de $i$ et $j$. A $j$ fixé il faut calculer $IP$ qui est le polynôme inverse à la puissance $a+1$ qui vaut au plus $p$. Il a donc un coût qui dépend du degré de l'inverse dont on vient de voir qu'il vaut au plus $2pl$. Le coût du calcul de $IP$ est donc en $O(p^2l)$. \\

Pour chaque valeur de $i$ il s'agit ensuite de calculer la somme dont les termes ne nécessitent qu'une seule opération. Le nombre maximum de termes de cette somme est $deg(IP)+1$ il est donc inférieur ou égale à $2p^2l+1$ et le coût de calcul de la somme est donc en $O(p^2l)$. A $j$ fixé on a donc à calculer $IP$ en $O(p^2l)$ puis tous les termes pour chaque valeur de $i$ entre $1$ et $rj-1$ ce qui a un coût en $O(p^2 l^2)$. \\

Enfin pour calculer tout le bloc il faut faire varier $j$ et le calcul est donc en $O(p^2l^2n)$ ceci domine clairement tous les autres calculs et le coût du calcul total est en $O(p^2gl)$. \\

\subsection{Complexité pour la matrice de Cartier}

On utilise la formule de la proposition \ref{projection}. On sait qu'il s'agit tout d'abord de calculer $v$ ainsi que l'inverse de $f$ à l'ordre $d+1$, calculs ne dépendant pas de $i$ ni de $j$ et dont le coût, qui a déjà été calculé à la section précédente, est en $O(pl)$ pour les deux. Il s'agit ensuite, à $j$ fixé de calculer $IP$ ce qui a un coût en $O(p^2l)$. \\

Enfin il s'agit de calculer le polynôme $Q$ toujours à $j$ fixé. Tout d'abord il faut dériver le polynôme $v$ ce qui a un coût en le degré de $v$ c'est-à-dire en $O(pl)$ et l'ajout du monôme $s^{p+1}$ a un coût de $1$. Il faut ensuite calculer la puissance $a$ du polynôme $f_2$ qui a un coût en $O(pl)$ comme nous l'avons déjà vu. Enfin il faut diviser deux polynômes qui a un coût en le maximum des degrés i.e. en $O(pl)$. Le calcul de $Q$ a donc un coût en $O(pl)$. \\

Pour calculer un terme de la somme à $j$ fixé il faut effectuer $11$ opérations. Par ailleurs dans cette double somme il y a $(d^{\circ}IP+1)(d-(pi+k)-r(n-b))$ termes. Or on sait que le degré de $IP$ est inférieur ou égal à $(2p(l-1)+p)p$, que $r(n-b)$ est positif ou nul que $d$ vaut au plus $2p(l-1)+p$ et que $pi+k$ vaut au plus $0$, il y a donc au maximum $p(2p(l-1)+p)(2p(l-1)+p)$ termes dans cette somme et le calcul de la somme a un coût en $O(p^3l^2)$ et comme il doit être effectué pour chaque valeur de $i$ le coût est en $O(p^3l^3)$. \\

Enfin pour calculer tout le bloc il faut faire varier $j$ et le calcul est donc en $O(p^3l^3n)$ ceci domine clairement tous les autres calculs et le coût du calcul total est en $O(p^3gl^2)$. \\

\section{Quelques exemples pratiques}

Plusieurs exemples de courbes ont été traités. On n'exposera ici qu'une poignée de ces exemples. Un cas qui nous a semblé intéressant a été de fixer le polynôme $f$ mais de faire varier le nombre premier $p$ considéré. Nous présentons ici le cas $y^3=(t-1)(t-2)(t-3)(t-4)(t-5)$ et les valeurs de $p$ suivantes, $p=17$, $p=31$ et $p=41$. Dans chacun des ces trois cas, la matrice du Frobenius divisé, présentée de manière à ce que le bloc supérieur gauche soit la matrice de Cartier et le bloc inférieur droit la matrice de Hasse-Witt, a la forme suivante :
\begin{center}
$p= 17 \ \ \ det(A4)= 0$ \hspace{120pt} $p= 31 \ \ \ det(A4)= 8$
$$\left(\begin{array}{m{0,4cm}m{0,4cm}m{0,4cm}m{0,4cm}|m{0,4cm}m{0,4cm}m{0,4cm}m{0,4cm}}
 0 &  12 &  16 &  6 &  6 &  0 &  0 &  0 \\
 8 &  0 &  0 &  0 &  0 &  15 &  1 &  11 \\
 9 &  0 &  0 &  0 &  0 &  7 &  0 &  12 \\
 7 &  0 &  0 &  0 &  0 &  1 &  4 &  15 \\
\hline
 10 &  0 &  0 &  0 &  0 &  10 &  8 &  0 \\
 0 &  13 &  5 &  8 &  0 &  0 &  0 &  0 \\
 0 &  5 &  10 &  11 &  13 &  0 &  0 &  0 \\
 0 &  8 &  11 &  0 &  10 &  0 &  0 &  0
\end{array}\right)
\ \ \ \ \ \ \ \ \ \ 
\left(\begin{array}{m{0,4cm}m{0,4cm}m{0,4cm}m{0,4cm}|m{0,4cm}m{0,4cm}m{0,4cm}m{0,4cm}}
 21 &  0 &  0 &  0 &  0 &  8 &  2 &  4 \\
 0 &  17 &  19 &  25 &  8 &  0 &  0 &  0 \\
 0 &  14 &  10 &  16 &  30 &  0 &  0 &  0 \\
 0 &  8 &  8 &  14 &  20 &  0 &  0 &  0 \\
\hline
 0 &  0 &  16 &  3 &  11 &  0 &  0 &  0 \\
 0 &  0 &  0 &  0 &  0 &  9 &  11 &  24 \\
 28 &  0 &  0 &  0 &  0 &  8 &  29 &  10 \\
 13 &  0 &  0 &  0 &  0 &  29 &  15 &  9
\end{array}\right)
$$
\end{center}
\begin{center}
$p= 41 \ \ \ det(A4)= 0$
$$\left(\begin{array}{m{0,4cm}m{0,4cm}m{0,4cm}m{0,4cm}|m{0,4cm}m{0,4cm}m{0,4cm}m{0,4cm}}
 0 &  1 &  31 &  30 &  36 &  0 &  0 &  0 \\
 5 &  0 &  0 &  0 &  0 &  37 &  27 &  14 \\
 34 &  0 &  0 &  0 &  0 &  27 &  9 &  6 \\
 8 &  0 &  0 &  0 &  0 &  14 &  6 &  40 \\
\hline
 33 &  0 &  0 &  0 &  0 &  0 &  0 &  0 \\
 0 &  36 &  26 &  4 &  0 &  0 &  0 &  0 \\
 0 &  26 &  6 &  31 &  0 &  0 &  0 &  0 \\
 0 &  4 &  31 &  38 &  0 &  0 &  0 &  0
\end{array}\right)$$
\end{center}
où $A4$ désigne la matrice de Hasse-Witt. En dernier exemple nous présentons le cas où $p$ vaut $13$, $n$ vaut $4$ et $f(t)=(t-1)(t-2)(t-3)(t-4)(t-5)(t-6)(t-7)$. Alors la matrice du Frobenius est de la forme :
$$\left(\begin{array}{m{0,4cm}m{0,4cm}m{0,4cm}m{0,4cm}m{0,4cm}m{0,4cm}m{0,4cm}m{0,4cm}m{0,4cm}|m{0,4cm}m{0,4cm}m{0,4cm}m{0,4cm}m{0,4cm}m{0,4cm}m{0,4cm}m{0,4cm}m{0,4cm}}
10 & 0 & 0 & 0 & 0 & 0 & 0 & 0 & 0 & 0 & 0 & 0 & 0 & 2 & 8 & 1 & 8 & 6 \\
0 & 4 & 9 & 6 & 0 & 0 & 0 & 0 & 0 & 0 & 6 & 8 & 10 & 0 & 0 & 0 & 0 & 0 \\
0 & 10 & 4 & 0 & 0 & 0 & 0 & 0 & 0 & 0 & 5 & 0 & 8 & 0 & 0 & 0 & 0 & 0 \\
0 & 1 & 8 & 8 & 0 & 0 & 0 & 0 & 0 & 0 & 8 & 9 & 10 & 0 & 0 & 0 & 0 & 0 \\
0 & 0 & 0 & 0 & 9 & 5 & 6 & 2 & 12 & 6 & 0 & 0 & 0 & 0 & 0 & 0 & 0 & 0 \\
0 & 0 & 0 & 0 & 11 & 3 & 9 & 6 & 12 & 7 & 0 & 0 & 0 & 0 & 0 & 0 & 0 & 0 \\
0 & 0 & 0 & 0 & 2 & 1 & 1 & 2 & 0 & 12 & 0 & 0 & 0 & 0 & 0 & 0 & 0 & 0 \\
0 & 0 & 0 & 0 & 6 & 12 & 0 & 8 & 2 & 11 & 0 & 0 & 0 & 0 & 0 & 0 & 0 & 0 \\
0 & 0 & 0 & 0 & 8 & 10 & 6 & 7 & 10 & 5 & 0 & 0 & 0 & 0 & 0 & 0 & 0 & 0 \\
\hline
 0 & 0 & 0 & 0 & 5 & 7 & 0 & 10 & 8 & 0 & 0 & 0 & 0 & 0 & 0 & 0 & 0 & 0 \\
0 & 0 & 9 & 7 & 0 & 0 & 0 & 0 & 0 & 0 & 9 & 11 & 10 & 0 & 0 & 0 & 0 & 0 \\
0 & 5 & 4 & 7 & 0 & 0 & 0 & 0 & 0 & 0 & 0 & 1 & 1 & 0 & 0 & 0 & 0 & 0 \\
0 & 1 & 4 & 1 & 0 & 0 & 0 & 0 & 0 & 0 & 3 & 12 & 1 & 0 & 0 & 0 & 0 & 0 \\
4 & 0 & 0 & 0 & 0 & 0 & 0 & 0 & 0 & 0 & 0 & 0 & 0 & 10 & 7 & 7 & 7 & 0 \\
3 & 0 & 0 & 0 & 0 & 0 & 0 & 0 & 0 & 0 & 0 & 0 & 0 & 10 & 8 & 10 & 12 & 1 \\
10 & 0 & 0 & 0 & 0 & 0 & 0 & 0 & 0 & 0 & 0 & 0 & 0 & 12 & 10 & 1 & 4 & 8 \\
11 & 0 & 0 & 0 & 0 & 0 & 0 & 0 & 0 & 0 & 0 & 0 & 0 & 11 & 1 & 7 & 12 & 0 \\
5 & 0 & 0 & 0 & 0 & 0 & 0 & 0 & 0 & 0 & 0 & 0 & 0 & 12 & 1 & 4 & 9 & 8
\end{array}\right)$$
et la matrice de Hasse-Witt est de déterminant nul sans être la matrice nulle. Pour finir rappelons que la remarque \ref{laremarquequifaut} nous permettait d'affirmer que les matrices calculées présentent des blocs de zéros et $n-1$ blocs inversibles. Il nous semble intéressant de faire apparaître ces blocs pour les quatre exemples que nous venons de présenter. Ainsi une présentation de ces matrices adaptée à la fois à la permutation et à la filtration donne les matrices suivantes (présentées dans le même ordre que précédemment) :
\begin{center}
$n=3 \ \ \ p=17$ \hspace{160pt} $n=3 \ \ \ p=31$
$$\left(\begin{array}{m{0,4cm}m{0,4cm}m{0,4cm}m{0,4cm}|m{0,4cm}m{0,4cm}m{0,4cm}m{0,4cm}}
 0 &  0 &  0 &  0 &  12 &  16 &  6 &  6 \\
 0  &  0 &  0 &  0 &  13 &  5 &  8 &  0\\
 0  &  0 &  0 &  0&  5 &  10 &  11 &  13 \\
 0  &  0 &  0 &  0&  8 &  11 &  0 &  10 \\
\hline
 8 &  15 &  1 &  11&  0 &  0 &  0 &  0  \\
 9 &  7 &  0 &  12&  0 &  0 &  0 &  0  \\
 7  &  1 &  4 &  15&  0 &  0 &  0 &  0 \\
 10 &  10 &  8 &  0 &  0 &  0 &  0 &  0  
\end{array}\right)
\ \ \ \ \ \ \ \ \ \ 
\left(\begin{array}{m{0,4cm}m{0,4cm}m{0,4cm}m{0,4cm}|m{0,4cm}m{0,4cm}m{0,4cm}m{0,4cm}}
 21 &  8 &  2 &  4&  0 &  0 &  0 &  0  \\
  0 &  9 &  11 &  24&  0 &  0 &  0 &  0  \\
 28 &  8 &  29 &  10&  0 &  0 &  0 &  0  \\
 13 &  29 &  15 &  9&  0 &  0 &  0 &  0 \\
 
\hline
 0  &  0 &  0 &  0 &  17 &  19 &  25 &  8\\
 0 &  0 &  0 &  0&  14 &  10 &  16 &  30  \\
 0  &  0 &  0 &  0&  8 &  8 &  14 &  20 \\
 0 &  0 &  0 &  0 &  0 &  16 &  3 &  11 
\end{array}\right)
$$
\end{center}
\begin{center}
$n=3 \ \ \ p=41$
$$\left(\begin{array}{m{0,4cm}m{0,4cm}m{0,4cm}m{0,4cm}|m{0,4cm}m{0,4cm}m{0,4cm}m{0,4cm}}
 0  &  0 &  0 &  0 &  1 &  31 &  30 &  36\\
  0  &  0 &  0 &  0 &  36 &  26 &  4 &  0\\
 0  &  0 &  0 &  0 &  26 &  6 &  31 &  0\\
 0  &  0 &  0 &  0&  4 &  31 &  38 &  0\\
\hline
 5&  37 &  27 &  14 &  0 &  0 &  0 &  0  \\
 34&  27 &  9 &  6 &  0 &  0 &  0 &  0  \\
 8 &  14 &  6 &  40&  0 &  0 &  0 &  0  \\
 33 &  0 &  0 &  0 &  0 &  0 &  0 &  0 
\end{array}\right)$$
\end{center}
\begin{center}
$n=4 \ \ \ p=13$
$$\left(\begin{array}{cccccc|cccccc|cccccc}
10 & 2 & 8 & 1 & 8 & 6 & 0 & 0 & 0 & 0 & 0 & 0 & 0 & 0 & 0 & 0 & 0 & 0 \\
4& 10 & 7 & 7 & 7 & 0 & 0 & 0 & 0 & 0 & 0 & 0 & 0 & 0 & 0 & 0 & 0 & 0  \\
3 & 10 & 8 & 10 & 12 & 1& 0 & 0 & 0 & 0 & 0 & 0 & 0 & 0 & 0 & 0 & 0 & 0  \\
10& 12 & 10 & 1 & 4 & 8 & 0 & 0 & 0 & 0 & 0 & 0 & 0 & 0 & 0 & 0 & 0 & 0  \\
11 & 11 & 1 & 7 & 12 & 0& 0 & 0 & 0 & 0 & 0 & 0 & 0 & 0 & 0 & 0 & 0 & 0  \\
5 & 12 & 1 & 4 & 9 & 8 & 0 & 0 & 0 & 0 & 0 & 0 & 0 & 0 & 0 & 0 & 0 & 0  \\
\hline
0 &  0 & 0 & 0 & 0 & 0 & 4 & 9 & 6  & 6 & 8 & 10 & 0 & 0 & 0 & 0 & 0 &0\\
0 &  0 & 0 & 0 & 0 & 0 & 10 & 4 & 0 & 5 & 0 & 8 & 0 & 0 & 0 & 0 & 0 &0\\
0 &  0 & 0 & 0 & 0 & 0 & 1 & 8 & 8  & 8 & 9 & 10 & 0 & 0 & 0 & 0 & 0 &0 \\
0 &  0 & 0 & 0 & 0 & 0 & 0 & 9 & 7 & 9 & 11 & 10 & 0 & 0 & 0 & 0 & 0&0 \\
0 &  0 & 0 & 0 & 0 & 0 & 5 & 4 & 7 & 0 & 1 & 1 & 0 & 0 & 0 & 0 & 0 &0\\
0 &  0 & 0 & 0 & 0 & 0 & 1 & 4 & 1 & 3 & 12 & 1 & 0 & 0 & 0 & 0 & 0 &0\\
\hline
0 & 0 & 0 & 0 & 0 & 0 & 0 & 0 & 0 & 0 & 0 & 0 & 9 & 5 & 6 & 2 & 12 & 6 \\
0 & 0 & 0 & 0  & 0 & 0 & 0 & 0 & 0 & 0 & 0 & 0 & 11 & 3 & 9 & 6 & 12 & 7\\
0 & 0 & 0 & 0  & 0 & 0 & 0 & 0 & 0 & 0 & 0 & 0 & 2 & 1 & 1 & 2 & 0 & 12\\
0 & 0 & 0 & 0 & 0 & 0 & 0 & 0 & 0 & 0 & 0 & 0 & 6 & 12 & 0 & 8 & 2 & 11 \\
0 & 0 & 0 & 0  & 0 & 0 & 0 & 0 & 0 & 0 & 0 & 0 & 8 & 10 & 6 & 7 & 10 & 5\\
 0 & 0 & 0 & 0  & 0 & 0 & 0 & 0 & 0 & 0 & 0 & 0 & 5 & 7 & 0 & 10 & 8 & 0
\end{array}\right)$$
\end{center}
qui par contre ne permettent plus de lire la matrice de Cartier ni la matrice de Hasse-Witt pour chacune des courbes considérées.

\newpage

\bibliographystyle{plain}
\nocite{*}
\bibliography{Z-Ref}

\end{document}